\documentclass[reqno,11pt]{amsart}
\pdfoutput=1

\usepackage{eucal}
\usepackage{outlines}
\usepackage{bbm}
\usepackage{amssymb,amsfonts,braket}
\usepackage{stmaryrd}

\usepackage{autobreak}
\usepackage{enumerate}
\usepackage{enumitem}
\usepackage{mathrsfs}
\usepackage{tikz-cd}
\usepackage{mathtools}
\usepackage{yhmath}
\usepackage{quiver}
\usepackage{bbm}
\usepackage{comment}
\usepackage{bm}
\usepackage{quiver} 
\usepackage{scalerel}
\usepackage{stackengine,wasysym}
\usepackage{csquotes}
\usepackage{xcolor}
\usepackage{amsmath,amsthm}
\usepackage{etoolbox}
\usepackage{calligra,mathrsfs}
\usepackage[style=alphabetic,
minalphanames=1,
maxalphanames=20,
maxbibnames=99]{biblatex}
\usepackage[colorlinks,allcolors=blue]{hyperref}
\usepackage{zref-clever}
\zcsetup{cap,nameinlink}

\newcommand{\cref}[1]{\zcref{#1}}
\newcommand{\Cref}[1]{\zcref[S]{#1}}

\newcommand{\labelcref}[1]{\zcref[noname]{#1}}

\ExplSyntaxOn
\newcommand{\zcrefglobalstringname}[2]{
	\__zrefclever_opt_varname_lang_type:nnnn{\languagename}{#1}{#2}{tl}
}
\newcommand{\zcreflocalstringname}[2]{
	\__zrefclever_opt_varname_type:een{#1}{#2}{tl}
}
\newcommand{\zcrefgetstring}[2]{
	\__zrefclever_provide_langfile:e { \languagename }
	\__zrefclever_process_language_settings:
	\ifcsvoid{\zcreflocalstringname{#1}{#2}}
	{\csname \zcrefglobalstringname{#1}{#2} \endcsname}
	{\csname \zcreflocalstringname{#1}{#2} \endcsname}
}
\ExplSyntaxOff
\newcommand{\zcreftypename}[1]{\zcrefgetstring{#1}{Name-sg}}

\newtheorem{sharedtheoremcounter}{}[section]
\undef\sharedtheoremcounter % prevent dummy theorem from being used

\NewDocumentCommand{\newzctheorem}{mO{#1}O{\zcreftypename{#2}}}{%
    \newtheorem{#1}[sharedtheoremcounter]{#3}%

    \zcsetup{countertype={#1=#2}}%

    \AtBeginEnvironment{#1}{%
        \zcsetup{countertype={sharedtheoremcounter=#2}}%
    }%
}
\zcRefTypeSetup{convention}{
	Name-sg = Convention ,
	name-sg = convention ,
	Name-pl = Conventions ,
	name-pl = conventions ,
}
\zcRefTypeSetup{keythm}{
	Name-sg = Key Theorem ,
	name-sg = key theorem ,
	Name-pl = Key Theorems ,
	name-pl = key theorems ,
}
\zcRefTypeSetup{mainthm}{
	Name-sg = Main Theorem ,
	name-sg = main theorem ,
	Name-pl = Main Theorems ,
	name-pl = main theorems ,
}
\zcRefTypeSetup{exmp}{
	Name-sg = Example ,
	name-sg = example ,
	Name-pl = Examples ,
	name-pl = examples ,
}
\zcRefTypeSetup{exmps}{
	Name-sg = Examples ,
	name-sg = examples ,
	Name-pl = Examples ,
	name-pl = examples ,
}
\zcRefTypeSetup{notn}{
	Name-sg = Notation ,
	name-sg = notation ,
	Name-pl = Notations ,
	name-pl = notations ,
}
\zcRefTypeSetup{notns}{
	Name-sg = Notations ,
	name-sg = notations ,
	Name-pl = Notations ,
	name-pl = notations ,
}
\zcRefTypeSetup{disc}{
	Name-sg = Discussion ,
	name-sg = discussion ,
	Name-pl = Discussions ,
	name-pl = discussions ,
}
\zcRefTypeSetup{addm}{
	Name-sg = Addendum ,
	name-sg = addendum ,
	Name-pl = Addenda ,
	name-pl = addenda ,
}
\zcRefTypeSetup{exer}{
	Name-sg = Exercise ,
	name-sg = exercise ,
	Name-pl = Exercises ,
	name-pl = exercises ,
}
\zcRefTypeSetup{property}{
	Name-sg = Property ,
	name-sg = property ,
	Name-pl = Properties ,
	name-pl = properties ,
}
\zcRefTypeSetup{properties}{
	Name-sg = Properties ,
	name-sg = properties ,
	Name-pl = Properties ,
	name-pl = properties ,
}
\zcRefTypeSetup{warn}{
	Name-sg = Warning ,
	name-sg = warning ,
	Name-pl = Warnings ,
	name-pl = warnings ,
}
\zcRefTypeSetup{conv}{
	Name-sg = Convention ,
	name-sg = convention ,
	Name-pl = Conventions ,
	name-pl = conventions ,
}
\zcRefTypeSetup{assumption}{
	Name-sg = Assumption ,
	name-sg = assumption ,
	Name-pl = Assumptions ,
	name-pl = assumptions ,
}
\zcRefTypeSetup{interlude}{
	Name-sg = Interlude ,
	name-sg = interlude ,
	Name-pl = Interludes ,
	name-pl = interludes ,
}
\zcRefTypeSetup{importantrem}{
	Name-sg = Important Remark ,
	name-sg = important remark ,
	Name-pl = Important Remarks ,
	name-pl = important remarks ,
}
\zcRefTypeSetup{rems}{
	Name-sg = Remarks ,
	name-sg = remarks ,
	Name-pl = Remarks ,
	name-pl = remarks ,
}
\zcRefTypeSetup{warning}{
	Name-sg = Warning ,
	name-sg = warning ,
	Name-pl = Warnings ,
	name-pl = warnings ,
}
\zcRefTypeSetup{sch}{
	Name-sg = Scholium ,
	name-sg = scholium ,
	Name-pl = Scholia ,
	name-pl = scholia ,
}
\zcRefTypeSetup{dig}{
	Name-sg = Digression ,
	name-sg = digression ,
	Name-pl = Digressions ,
	name-pl = digressions ,
}
\zcRefTypeSetup{con}{
	Name-sg = Construction ,
	name-sg = construction ,
	Name-pl = Constructions ,
	name-pl = constructions ,
}
\zcRefTypeSetup{irem}{
	Name-sg = Important Remark ,
	name-sg = important remark ,
	Name-pl = Important Remarks ,
	name-pl = important remarks ,
}
\zcRefTypeSetup{tempnotation}{
	Name-sg = Temporary Notation ,
	name-sg = temporary notation ,
	Name-pl = Temporary Notations ,
	name-pl = temporary notations ,
}
\zcRefTypeSetup{conventions}{
	Name-sg = Conventions ,
	name-sg = conventions ,
	Name-pl = Conventions ,
	name-pl = conventions ,
}

\zcRefTypeSetup{defns}{
	Name-sg = Definitions ,
	name-sg = definitions ,
	Name-pl = Definitions ,
	name-pl = definitions ,
}

\zcRefTypeSetup{claim}{
	Name-sg = Claim ,
	name-sg = claim ,
	Name-pl = Claims ,
	name-pl = claims ,
}

\zcRefTypeSetup{notation}{
	Name-sg = Notation ,
	name-sg = notation ,
	Name-pl = Notations ,
	name-pl = notations ,
}

\zcRefTypeSetup{notnn}{
	Name-sg = Notation ,
	name-sg = notation ,
	Name-pl = Notations ,
	name-pl = notations ,
}

\zcRefTypeSetup{warning}{
	Name-sg = Warning ,
	name-sg =  warning ,
	Name-pl = Warnings ,
	name-pl = warnings ,
}

\zcRefTypeSetup{quest}{
	Name-sg = Question ,
	name-sg = question ,
	Name-pl = Questions ,
	name-pl = questions ,
}

\theoremstyle{plain}
\newzctheorem{theorem}
\newzctheorem{keythm}[keythm][Key Theorem]
\newzctheorem{lem}[lemma]
\newzctheorem{thm}[theorem]
\newzctheorem{thmm}[theorem]
\newzctheorem{cor}[corollary]
\newzctheorem{prop}[proposition]
\newzctheorem{conj}[conjecture]
\newzctheorem{quest}[quest]
\newzctheorem{claim}[claim][Claim]
\newzctheorem{mainthm}[mainthm][Main Theorem]

\theoremstyle{definition}
\newzctheorem{defi}[definition]
\newzctheorem{convention}
\newzctheorem{defn}[definition][Definition]	
\newzctheorem{defns}[defns][Definitions]
\newzctheorem{exmp}[exmp][Example]
\newzctheorem{exmps}[exmps][Examples]
\newzctheorem{notn}[notn][Notation]
\newzctheorem{notns}[notns][Notations]
\newzctheorem{disc}[disc][Discussion]
\newzctheorem{addm}[addm][Addendum]
\newzctheorem{exer}[exer][Exercise]
\newzctheorem{property}[property][Property]
\newzctheorem{properties}[properties][Properties]
\newzctheorem{warn}[warn][Warning]
\newzctheorem{remark}[remark][Remark]
\newzctheorem{assumption}[assumption][Assumption]
\newzctheorem{interlude}[interlude][Interlude]
\newzctheorem{importantrem}[importantrem][Important Remark]
\newzctheorem{con}[con][Construction]

\theoremstyle{remark}
\newzctheorem{rems}[rems][Remarks]
\newzctheorem{warning}[warning][Warning]
\newzctheorem{sch}[sch][Scholium]	
\newzctheorem{rem}[remark][Remark]
\newzctheorem{dig}[dig][Digression]
\newzctheorem{notation}[notation][Notation]
\newzctheorem{irem}[irem][Important Remark]
\newzctheorem{tempnotation}[tempnotation][Temporary Notation]
\newzctheorem{conventions}[conventions][Conventions]
\newzctheorem{notnn}[notnn][Notation]
\newzctheorem{conv}[conv][Convention]

\makeatletter
\let\c@equation\c@sharedtheoremcounter

\makeatother

\newcommand{\beqn}{\begin{equation}}
	\newcommand{\eeqn}{\end{equation}}
\newcommand{\bclaim}{\begin{claim}}
	\newcommand{\eclaim}{\end{claim}}
\newcommand{\blem}{\begin{lem}}
	\newcommand{\elem}{\end{lem}}
\newcommand{\bproof}{\begin{proof}}
	\newcommand{\eproof}{\end{proof}}
\newcommand{\bdef}{\begin{defn}}
	\newcommand{\edefn}{\end{defn}}
\newcommand{\bprop}{\begin{prop}}
	\newcommand{\eprop}{\end{prop}}
\newcommand{\bthm}{\begin{thm}}
	\newcommand{\ethm}{\end{thm}}
\newcommand{\brem}{\begin{rem}}
	\newcommand{\erem}{\end{rem}}
\newcommand{\bcor}{\begin{cor}}
	\newcommand{\ecor}{\end{cor}}

\newcommand{\calD}{\mathcal{D}}

\newcommand{\calF}{\mathcal{F}}

\newcommand{\calH}{\mathcal{H}}

\newcommand{\calL}{\mathcal{L}}
\newcommand{\calM}{\mathcal{M}}

\newcommand{\calO}{\mathcal{O}}
\newcommand{\calP}{\mathcal{P}}

\newcommand{\scrC}{\mathscr{C}}
\newcommand{\scrD}{\mathscr{D}}

\newcommand{\scrR}{\mathscr{R}}

\renewcommand{\b}{\beta}

\newcommand{\on}{\operatorname}

\newcommand{\wih}{\widehat}

\newcommand{\hook}{\hookrightarrow}

\newcommand{\kb}{k[\b^{\pm 1}]}

\newcommand{\mathendash}{\text{\textendash}}
\renewcommand{\mod}{\mathendash\mathrm{mod}}

\mathchardef\mhyphen="2D

\renewcommand{\bmod}{\mathendash\mathbf{mod}}

\newcommand{\heart}{\heartsuit}

\newcommand{\oblv}{\mathrm{oblv}}

\newcommand{\Cat}{\mathcal{C}\mathrm{at}}
\newcommand{\id}{\mathrm{id}}
\newcommand{\Map}{\calM\mathrm{ap}}

\newcommand{\sm}{\mathrm{sm}}

\newcommand{\dgcat}{\mathbf{dgcat}}
\newcommand{\bCat}{\mathbf{Cat}}
\newcommand{\bDGCat}{\mathbf{DGCat}}
\newcommand{\bPr}{\mathbf{Pr}}
\newcommand{\bCoh}{\mathbf{Coh}}
\newcommand{\bReg}{\mathbf{Reg}}
\newcommand{\bSpan}{\mathbf{Span}}

\newcommand{\QCoh}{\on{\mathsf{QCoh}}}
\newcommand{\IndCoh}{\on{\mathsf{IndCoh}}}
\newcommand{\Vect}{\mathsf{Vect}}

\newcommand{\Coh}{\on{\mathsf{Coh}}}
\newcommand{\Perf}{\on{\mathsf{Perf}}}

\newcommand{\Ind}{\on{\mathsf{Ind}}}

\newcommand{\Reg}{\mathsf{Reg}}

\newcommand{\Fun}{\on{Fun}}
\newcommand{\bFun}{\mathbf{Fun}}

\newcommand{\op}{\on{op}}
\newcommand{\pt}{\on{pt}}

\newcommand{\st}{\on{st}}

\newcommand{\DGCat}{\mathrm{DGCat}}

\DeclareMathOperator{\colim}{colim}

\DeclareMathOperator{\Mod}{Mod}

\DeclareMathOperator{\Spec}{Spec}
\DeclareMathOperator{\Hom}{Hom}

\DeclareMathOperator{\sHom}{\mathscr{H}\text{\kern -3pt {\calligra\large om}}\,}

\renewcommand{\Pr}{\calP\mathrm{r}}

\newcommand{\Ch}{\mathsf{Ch}}

\newcommand{\im}{\on{im}}

\newcommand{\ev}{\mathrm{ev}}

\newcommand{\Mor}{\calM\mathrm{or}}

\renewcommand{\sHom}{\calH\mathrm{om}}

\newcommand{\ku}{k[u^{\pm 1}]}

\newcommand{\Stk}{\mathsf{Stk}}

\newcommand{\Grpd}{\mathsf{Grpd}}

\newcommand{\laft}{\mathrm{laft}}

\newcommand{\Corr}{\mathrm{Corr}}

\newcommand{\tT}{t\mathrm{Tate}}
\newcommand{\Pot}{\calD^{\mathrm{pot}}}

\newcommand{\bbone}{\bm{1}}

\newcommand{\triv}{\mathrm{triv}}
\newcommand{\Ob}{\calO\mathrm{b}}

\newcommand{\sfMap}{\mathsf{Map}}

\DeclareFontFamily{U}{mathx}{\hyphenchar\font45}
\DeclareFontShape{U}{mathx}{m}{n}{
	<5> <6> <7> <8> <9> <10>
	<10.95> <12> <14.4> <17.28> <20.74> <24.88>
	mathx10
}{}
\DeclareSymbolFont{mathx}{U}{mathx}{m}{n}
\DeclareFontSubstitution{U}{mathx}{m}{n}
\DeclareMathAccent{\widecheck}{0}{mathx}{"71}
\DeclareMathAccent{\wideparen}{0}{mathx}{"75}

\renewcommand{\colim}{\operatornamewithlimits{colim}}

\makeatletter
\newcommand*\bigcdot{\mathpalette\bigcdot@{.5}}
\newcommand*\bigcdot@[2]{\mathbin{\vcenter{\hbox{\scalebox{#2}{$\m@th#1\bullet$}}}}}
\makeatother

\usepackage{geometry}
\makeatletter
\renewcommand{\widehat}{%
  \mathaccent"0\hexnumber@\symyhlargesymbols62\relax
}
\makeatother

\title{A six-functor formalism for potent $\calD$-modules}

\author{Jacob Erlikhman and Kendric Schefers}

\begin{document}
\raggedbottom
\begin{abstract}
On stacks there is a richer version of $\calD$-modules---called \textit{potent $\calD$-modules}---which decompletes at equivariant parameters, in the same sense that genuine equivariant $K$-theory decompletes equivariant cohomology. In this paper, we construct---for a broad class of geometric stacks and their formal completions---a six-functor formalism for potent $\calD$-modules. We also show that on finite type separated schemes, this six-functor formalism is equivalent to the Gaitsgory-Rozenblyum six-functor formalism for 2-periodic $\calD$-modules.
\end{abstract}

\maketitle

\tableofcontents

\pagebreak
\raggedbottom

%--------------------------------------------------------------------------------------------------%
%--------------------------------------INTRODUCTION-------------------------------------------%
%--------------------------------------------------------------------------------------------------%

\section{Introduction}
\subsection{Statement of main results}
Let $X$ be a smooth scheme of finite type over a charactersitic 0 field. As explained in \cite{GRCrystals}, one can understand the derived category of (left) $\mathcal{D}$-modules on $X$ as the derived category of quasi-coherent sheaves on the de Rham prestack $X_{\mathrm{dR}}$ of $X$, called the category of \textit{(left) crystals} on $X$. The de Rham prestack is a space whose points are the infinitesimally identified points of $X$.\footnote{More precisely, for a prestack $X$, the de Rham prestack is defined as the functor which to an affine scheme $ S$ assigns $X(S_{\mathrm{red}})$, where $S_{\mathrm{red}}$ is the underlying classical reduced scheme of $S$.} More precisely, two $S$-points $x,y\in X(S)$ are \textit{infinitesimally close} if they map to the same point under the canonical map $X(S)\to X(S_{\mathrm{red}})$. Quasi-coherent sheaves on $X_{\mathrm{dR}}$ are, by definition, quasi-coherent sheaves $\mathcal{F}$ on $X$ such that for any two infinitesimally close points $x,y$ of $X$, $x^*\mathcal{F}\simeq y^*\mathcal{F}$, together with a ($\infty$-categorical) cocycle condition. The reason that this recovers (left) $\mathcal{D}$-modules is that a flat connection on a vector bundle is equivalent to an identification of nearby fibers of that vector bundle (by parallel transport). In algebraic geometry, points on a scheme can only be ``nearby'' if they are infinitesimally close, whence the definition. For technical reasons (explained in \textit{loc. cit.}), the category of right $\mathcal{D}$-modules has better properties than that of left $\mathcal{D}$-modules; it is equivalent to the category of \textit{right crystals}: ind-coherent sheaves on the de Rham prestack. We denote this category by
\begin{align*}
	\mathsf{Crys}^r(X)\coloneq \mathsf{IndCoh}(X_{\mathrm{dR}}).
\end{align*} 

In \cite{BZNLoops}, Ben-Zvi--Nadler explained that one can recover (2-periodicized) $\mathcal{D}$-modules on $X$ by another construction. Namely, they prove that Koszul duality and a version of the Hochschild-Kostant-Rosenberg isomorphism identify 2-periodic $\mathcal{D}$-modules with cyclic sheaves on the derived free loop stack of $X$\footnote{This stack is the derived stack $\mathcal{L}X\coloneq \Map(S^1,X)\simeq X\times _{X\times X}X$, where $S^1$ is the algebraic stack $B\mathbb{Z} $, or, equivalently, the constant derived stack with value $S^ 1\in \mathsf{Spc}$.}:
\begin{align}\label{eq: potent D-modules are D-modules}
	\mathsf{IndCoh}(\mathcal{L}X)^{\tT}\simeq \mathcal{D}_X\mod_{\mathbb{Z} /2}\simeq \mathsf{Crys}^r(X)_{\mathbb{Z} /2},\footnotemark
\end{align}
\footnotetext{The category $\mathsf{IndCoh}(\mathcal{L}X)^{\tT}$ is the regularized Tate invariants, defined in \zcref{sssec:Tate invariants} below, with respect to the $S^1$-action on $\mathsf{IndCoh}(\mathcal{L}X)$ induced by the loop rotation action on $\mathcal{L}X\coloneq \Map(S^1,X)$. Concretely, it is the ind-completion of $\mathsf{Coh}(\mathcal{L}X)^{tS^1}$, the ordinary Tate invariants (with respect to the loop rotation action) of the derived category of coherent sheaves on $\mathcal{L}X$.}
where $\mathcal{D}_X\mod_{\mathbb{Z} /2}$ denotes the category of (2-perodicized) $\mathcal{D}$-modules on $X$. This category of cyclic sheaves, called the category of \textit{potent $\mathcal{D}$-modules} in the ICM address \cite{BZPotent}, was the main object of study of \cite{ES-Kashiwara}. We refer the reader to that paper for motivation for and instructive examples of potent $\mathcal{D}$-modules.

In \cite[Theorem 1.3.5]{Preygel}---which also extended \zcref[noname]{eq: potent D-modules are D-modules} to the non-smooth case, the category of potent $\mathcal{D}$-modules on a finite type scheme was shown to be equivalent to (the 2-periodicization of) $\mathsf{Crys}^r(X)$. However, as explained in \cite{ES-Kashiwara}, these categories differ when we allow $X$ to be a stack (with affine diagonal)---more precisely, $\mathcal{D}$-modules are a full subcategory of potent $\mathcal{D}$-modules.

In the generality of prestacks locally almost of finite type (laft), the assignment $X\leadsto \mathsf{Crys}^r(X)$ was shown to give rise to a six-functor formalism in \cite[Chapter 4 Section 2]{GR17II}. In this paper, we show that the assignment $X\leadsto \Pot(X)$ also gives rise to a six-functor formalism. Before stating the theorem, we caution that we were not able to construct this formalism in the generality of laft prestacks.\footnote{Consider the laft prestack $\coprod_{n\ge 0}\mathbb{A}^n$, and let $f$ be its structure morphism to $\Spec k$. On the $n^{\text{th}}$ component, $f^!k$ is $\Omega^n_{\mathbb{A}^n}[n]$, so it has no global lower bound. We require such lower bounds on $f^!$ in order for $\tT$ invariants to define a functor. See \zcref{sssec: intro tTate functoriality} for details.} Rather, we are only able to construct this formalism for \textit{ind-finite type stacks}---roughly, these are (possibly formal) stacks locally (almost) of finite type which are filtered colimits of algebraic stacks (almost) of finite type and with affine diagonal. Denote the category of such by $\mathsf{indGStk}_{\mathrm{ft}}$.
\begin{thm}[\zcref{thm:potent-D-modules-is-a-6FF}]\label{intro: potent D-modules is a 6FF}
	There is a lax symmetric monoidal $(\infty,2)$-functor
	\begin{align*}
		\Pot:\bSpan(\mathsf{indGStk}_{\mathrm{ft}},\mathrm{ind}\text{-}\mathrm{inf}\text{-}\mathrm{schem\&proper})_{\mathrm{equiv,ind\text{-}\mathrm{inf}\text{-}schem\&proper}}\to \bDGCat_{\ku},
	\end{align*} 
	where for an $\infty$-category $\mathscr{C}$ and suitable\footnote{See \cite[Convention 3.1]{CLL_six-functor-formalisms} for the precise conditions on $I,P,E,\mathscr{C}$.} classes of morphisms $I,P\subset E$ in $\mathscr{C}$, $\bSpan(\mathscr{C},E)_{P,I}$ denotes the $(\infty,2)$-category of correspondences in $\mathscr{C}$.
\end{thm} 
\begin{rem}
    By \zcref{lem: ind-inf-schematic maps are schematic}, the six-functor formalism of \zcref{intro: potent D-modules is a 6FF} is equivalent to the formalism with ind-schematic ind-proper maps instead of ind-inf-schematic ind-proper maps.
\end{rem}
We also extend Preygel's equivalence \cite[Theorem 1.3.5]{Preygel} to an equivalence of six-functor formalisms.
\begin{thm}[\zcref{prop: potent D-modules are crystals}]\label{intro: potent D-modules are crystals}
	For separated derived schemes almost of finite type, the objectwise equivalence $\pi_X: \mathsf{Crys}^r_{\mathbb{Z} /2}(X)\xrightarrow{\simeq} \Pot(X)$ of \cite[Theorem 1.3.5]{Preygel} induces an equivalence of six-functor formalisms
	\begin{align*}
		 \mathsf{Crys}^r_{\mathbb{Z} /2}\xrightarrow{\simeq} \Pot_{\mathrm{sch}},
	\end{align*} 
	where $\Pot_{\mathrm{sch}}$ denotes the restriction of the six-functor formalism of \zcref{intro: potent D-modules is a 6FF} to derived schemes almost of finite type. 
\end{thm} 
The proof of \zcref{intro: potent D-modules is a 6FF} proceeds by first constructing a three-functor formalism and then extending it to a six-functor formalism using (a generalization of) \cite{Mann}.
\subsection{The 3-functor formalism}
In this section, we explain first how to construct the underlying three-functor formalism for potent $\mathcal{D}$-modules. Once this is in hand, the proof that it extends to a six-functor formalism is largely formal.
\subsubsection{Functoriality of $\tT$ invariants}
\label{sssec: intro tTate functoriality}
By work of Gaitsgory-Rozenblyum, there is a well-defined six-functor formalism for ind-coherent sheaves on laft prestacks (\cite[Chapter 3 Theorem 5.4.3]{GR17II}). Composing it with the free loop stack functor
\begin{align*}
	\mathcal{L}: \mathsf{PreStk}_{\mathrm{laft}}\to \mathsf{PreStk}_{\mathrm{laft}},
\end{align*} 
we show (\zcref{cor:indcoh-loops-is-a-6FF}) that the composite $\mathsf{IndCoh}(\mathcal{L}-)$ defines a six-functor formalism on all laft prestacks. With this in hand, one could hope that the following hold, which would allow us to conclude the construction.
\begin{enumerate}[label=\roman*)]
	\item $(-)^{\tT}$ defines an $(\infty,2)$-functor.
	\item The underlying $(\infty,1)$-functor of $(-)^{\tT}$ is lax symmetric monoidal.
	\item $(-)^{\tT}$ can be composed with the six-functor formalism $\mathsf{IndCoh}(\mathcal{L}-)$.
\end{enumerate}
Although (i) is true (\zcref{cor: tTate is a 2-functor}), both (ii) and (iii) fail as stated. Failure of (ii) is technical: The definition of $\tT$ invariants requires as input a coherent category $\mathscr{C}$, which is a stable $\infty$-category equipped with a coherent $t$-structure in the sense of \cite[Definition 4.2.7]{Preygel}. More precisely, we show in \zcref{cor: tTate is a 2-functor} that $\tT$ invariants defines an $(\infty,2)$-functor
\begin{align*}
	\bFun(BS^1,\bCoh_t)\to \bDGCat_{\ku}
\end{align*} 
from coherent $k$-linear presentable $\infty$-categories with $S^1$-action\footnote{Our definition of $\bCoh_t$ differs from that in \cite{Preygel} in several ways. See \zcref{ssec: glossary of categories} for the precise definition.} to $\ku$-linear presentable $\infty$-categories. However, the category $\bCoh_t$ is not closed under the Lurie tensor product! We will explain how we repair this after discussing (iii).

The failure of (iii) is even more serious. In order to compose $\tT$ invariants with the functor $\mathsf{IndCoh}(\mathcal{L}-)$, we need the standard $t$-structure on $\mathsf{IndCoh}(\mathcal{L}X)$ to be coherent and $!$-pullback and (when it exists) $* $-pushforward must be left $t$-exact up to a shift. This is not known in the generality of laft prestacks. Now, in \cite[Proposition 4.4.3]{Preygel}, Preygel shows that for any Artin stack $X$ of finite type, $\mathsf{IndCoh}(X)$ has coherent $t$-structure; however, the applications of the six-functor formalism for potent $\mathcal{D}$-modules to the stacks appearing in \cite{ES-Kashiwara} require also a definition of $\Pot(\wih{X})$ when $\wih{X}$ is a formal completion of an algebraic stack $X$ of finite type. Thus, in \zcref{ssec: ind-finite type stacks} below, we define a category of laft stacks which includes algebraic stacks (almost) of finite type as well as their formal completions, denoted $\mathsf{indGStk}_{\mathrm{ft}}$. It is a full subcategory of the category of (ind-tamely presented) ind-geometric stacks of \cite{CautisWilliams,CautisWilliams2}. The category $\mathsf{indGStk}_{\mathrm{ft}}$ is closed under finite limits (\zcref{prop: ind-ft stacks are closed under limits})---so $\mathcal{L}X\in \mathsf{indGStk}_{\mathrm{ft}}$ if $X$ is---and has the property that the standard $t$-structure on $\mathsf{IndCoh}(X)$ for $X\in \mathsf{indGStk}_{\mathrm{ft}}$ is coherent (\zcref{lem: IndCoh is regular for ind-finite type stacks}), thus solving the functoriality issue on objects. For 1-morphisms, Preygel proves that $\tT$ invariants are functorial with respect to functors that are continuous, exact, and left $t$-exact up to a shift \cite[Section 4.6]{Preygel}. The first two properties of $*$-pushforward and $!$-pullback are guaranteed by results of Gaitsgory-Rozenblyum; however, left $t$-exactness up to a shift for $!$-pullback is shown by Gaitsgory only for schemes in \cite[Lemma 7.1.7]{GaitsgoryIndCoh}. The \textit{raison d'\^etre} of the category $\mathsf{indGStk}_{\mathrm{ft}}$ is the following. Any object $X$ has a presentation as a filtered colimit along (almost) finitely presented closed immersions, $X\simeq\colim_\alpha X_\alpha$, where the $X_\alpha$ are algebraic stacks (almost) of finite type (\zcref{lem: ind-finite type stacks have presentations by algebraic stacks}). This property implies that the restriction of the functor $\mathsf{IndCoh}(\mathcal{L}-)$ to the category $(\mathsf{indGStk}_{\mathrm{ft})^{\mathrm{op}}}$ has essential image contained in $\bFun(BS^1,\bCoh_t)$. We use \zcref{lem: ind-finite type stacks have presentations by algebraic stacks} to reduce left $ t$-exactness up to a shift of $!$-pullback to the case for schemes by smooth descent arguments (see \zcref{ssec: 2-functor of potent D-modules}). Functoriality on 2-morphisms is a non-issue, because the 2-morphisms in $\bCoh_t$ are the same as those in the $(\infty,2)$-category of presentable $k$-linear categories. 

Thus, with (iii) repaired, we have a well-defined $(\infty,2)$-functor
\begin{align}\label{eq: intro 3FF potent D-modules}
	\Pot\coloneq \mathsf{IndCoh}(\mathcal{L}-)^{\tT}: (\mathsf{indGStk}_{\mathrm{ft}})^{\mathrm{op}}\to \bDGCat_{\ku}
\end{align} 
to $\ku$-linear stable presentable $\infty$-categories.

\subsubsection{The lax symmetric monoidal structure}
Before explaining how we circumvent (ii), we give a more detailed description of $\tT$ invariants. Using \cite[Proposition 4.4.5]{Preygel}, we define the functor $ (-)^{\tT}$ as the composite
\begin{align}\label{eq: composite tTate invariants}
	\bFun(BS^1,\bCoh_t)\xrightarrow{(-)^{S^1}}\bCoh_{t,k[u]}\xrightarrow{\mathscr{R}_{k[u]}}\bReg_{t,k[u]}\xrightarrow{\oblv_t}\bDGCat_{k[u]}\xrightarrow{-\otimes _{k[u]}\ku}\bDGCat_{\ku}.
\end{align} 
We will now give an overview of the various categories which appear in this definition; full details can be found in \zcref{ssec: glossary of categories}. For a commutative (dg) $k$-algebra $A$, $\bCoh_{t,A}$ is the $(\infty,2)$-category of stable presentable $A $-linear $\infty$-categories with coherent $t$-structure; $\bReg_{t,A}\subset \bCoh_{t,A}$ is the full $(\infty,2)$-subcategory whose objects have regular $t$-structures (in the sense of \cite[Definition 4.3.1]{Preygel}); $\bDGCat_A$ is the $(\infty,2)$-category of stable presentable $A$-linear $\infty$-categories; and $\bFun$ denotes the $(\infty,2)$-category of functors. The functors, in order of appearance in the definition, are $S^1$-invariants, $k[u]$-linear regularization (a $k [u]$-linear enhancement of the regularization operation constructed in \cite[Proposition 4.3.2]{Preygel}), the forgetful functor (forgetting the $t$-structure), and extension of scalars. 

To circumvent the failure of (ii), we use Lurie's theory of $\infty$-operads of \cite{HA}. Namely, we equip the underlying $(\infty,1)$-categories of $\bCoh_{t,A}$ and $\bReg_{t,A}$ with the structures of $\infty$-operads (\zcref{con: Coh_t and Reg_t operads} and \zcref{prop: Coh_t and Reg_t are operads}). We then construct the functor $\mathscr{R}_{k[u]}$ as a map of $\infty$-operads (\zcref{lem: P -> Reg_t is a map of operads} and \zcref{lem: P is equivalent to Coh_t}). The forgetful functor $\oblv_t$ is a map of $\infty$-operads by construction (\zcref{con: Coh_t and Reg_t operads}), and extension of scalars is canonically lax symmetric monoidal, hence also a map of $\infty$-operads. Thus, we construct $\tT$ invariants not as a lax symmetric monoidal functor but as a map of $\infty$-operads. We then show that $\mathsf{IndCoh}(\mathcal{L}-)$ is also a map of $\infty$-operads (\zcref{prop: IndCoh(L-)^S^1 as a map of operads}), which allows us to construct $\Pot$ as a map of $\infty$-operads (\zcref{prop: D^pot is lax symmetric monoidal}). Since the source and target of $\Pot$ are symmetric monoidal categories, this is the desired lax symmetric monoidal structure on the functor \zcref[noname]{eq: intro 3FF potent D-modules}. This gives the desired three-functor formalism.
 
\subsection{The six-functor formalism}
To construct a six-functor formalism from the three-functor formalism \zcref[noname]{eq: intro 3FF potent D-modules}, we use \cite[Theorem 5.14]{CLL_six-functor-formalisms}. In particular, we need to choose certain distinguished classes of maps $I,P\subset E$ in $\mathsf{indGStk}_{\text{ft}}$ such that if $f\in I$, $\Pot(f)$ admits a left adjoint and if $f\in P$, then $\Pot(f) $ admits a right adjoint. Since our starting point is the Gaitsgory-Rozenblyum six-functor formalism for $\mathsf{IndCoh}$, we choose to take $I=E$ to be the class of ind-inf-schematic ind-proper maps and $P$ to be the class of equivalences. Indeed, the Gaitsgory-Rozenblyum six-functor formalism for ind-coherent sheaves gives that $f^!$ has a left adjoint, $f_*$, if $f$ is an ind-inf-schematic ind-proper map of laft prestacks (\cite[Chapter 3 Theorem 5.4.3]{GR17II}). Functoriality of the three-functor formalism then formally gives left and right adjoints after passing through $\tT$ invariants, by applying the functor of $\tT$ invariants to the unit and counit morphisms of the adjunctions. For the free loop stack functor $\mathcal{L}$, we check, via a geometric argument, that it preserves the class $I$ (\zcref{prop:indcoh-loops-6FF-setup-for-potent-d-modules}), noting that it trivially preserves the class of equivalences. 

Now, the application of \cite[Theorem 5.14]{CLL_six-functor-formalisms} requires four distinct inputs.
\begin{enumerate}[label=\arabic*)]
	\item The classes $I,P,E$ and the category $\mathsf{indGStk_{\mathrm{ft}}}$ satisfy the conditions of \cite[Convention 3.1]{CLL_six-functor-formalisms}; $I,P$ are closed under tensor products; and the symmetric monoidal structure on $\mathsf{indGStk}_{\mathrm{ft}}$ preserves pullbacks along morphisms in $I\times I$ and $P\times P$.
	\item Construct $\Pot$ as an $(I,P)$-biadjointable\footnote{By \cite[Example 2.7]{CLL_six-functor-formalisms}, it suffices to check left $I$-adjointability. This is defined in Definitions 2.3 and 2.4 of \textit{op. cit.}} $(\infty,2)$-functor
		\begin{align*}
			(\mathsf{indGStk}_{\text{ft}})^{\text{op}}\to \bDGCat_{\ku}.
		\end{align*} 
	\item Show that the lax symmetric monoidal structure on $\Pot$ is compatible with adjointability (see \cite[Definition 5.13(2)]{CLL_six-functor-formalisms} for the precise formulation).
\end{enumerate}
\subsubsection{Checks on $\mathsf{indGStk}_{\mathrm{ft}}$ and the classes $I$ and $P$}
This is straightforward (see the beginnings of the proofs of \zcref{thm:potent-D-modules-is-a-6FF} and \zcref{prop: indcoh is a 6-functor formalism}).

\subsubsection{$(I,P)$-biadjointability}
Crucially, this uses the already constructed six-functor formalism for $\mathsf{IndCoh}$, whose $(I,P)$-biadjointability follows from results of Gaitsgory-Rozenblyum. We spell this out---in the \cite{CLL_six-functor-formalisms} formalism---in \zcref{prop: indcoh is a 6-functor formalism}. We then show that $\mathsf{IndCoh}\circ \mathcal{L}$ also defines a six-functor formalism in \zcref{cor:indcoh-loops-is-a-6FF}, which reduces to checking that $\mathcal{L}$ preserves the class $I$. This is done by factoring $\mathcal{L}f$ for a map $f\in I$ as the composition of maps which are base changes of maps in $I$. Functoriality of $\tT$ invariants (\zcref{cor: tTate is a 2-functor}) applied to the relevant mates then gives left $I$-adjointability. Right $P$-adjointability follows formally for the class of equivalences (cf. \cite[Example 2.7]{CLL_six-functor-formalisms}). 

\subsubsection{Compatibility of the lax symmetric monoidal structure with biadjointability}
This is checked in \zcref{thm:potent-D-modules-is-a-6FF}. We give a short overview of the proof here. This essentially follows by applying the functor of $\tT$ invariants to the already-known compatibility of \zcref{cor:indcoh-loops-is-a-6FF}. However, since $(-)^{\tT}$ only has underlying $(\infty,1)$-functor that is a map of $\infty$-operads, it does not necessarily preserve mates of \textit{binary} squares of the form required in \cite[Definition 5.13(2)]{CLL_six-functor-formalisms}. It does however preserve unary such squares. Thus, fixing one of the entries of the exterior product functors, we reduce to a unary square. These unary functors are coherent if we plug in a coherent object (since coherent objects are bounded). Since coherent objects generate $\Pot$, we conclude.

\subsection{Comparison with (2-periodic) $\mathcal{D}$-modules}
In this section, we discuss the proof of the second main result of this paper, \zcref{intro: potent D-modules are crystals}. In \zcref{sssec: crystal 3FF}, we construct a three-functor formalism on separated derived schemes almost of finite type for 2-periodicized crystals by extending scalars for the $k$-linear category $\mathsf{Crys}^r(X)$ to $\ku$. In \zcref{sssec: crystal 3FF 2}, we construct an a priori different three-functor formalism given by equipping $\mathsf{Crys}^r(X)$ with the trivial $S^1$-action and taking $\tT$ invariants (which makes sense by \cite[Proposition 4.4.4]{Preygel}). We show (\zcref{lem: periodic crystals are tate invariants of crystals}) that these three-functor formalisms are equivalent. The latter three-functor formalism is more suitable for comparison with potent $\mathcal{D}$-modules. We show that it defines a six-functor formalism (\zcref{lem: crystal 6FF}).

Then, in \zcref{ssec: potent D-modules=crystals}, we construct a natural transformation of six-functor formalisms $\mathsf{Crys}^r_{\mathbb{Z}/2}\Rightarrow\Pot_{\text{sch}}$, whose value on a separated derived scheme almost of finite type is the equivalence of \cite[Theorem 1.3.5]{Preygel}, thus completing the proof of \zcref{intro: potent D-modules are crystals}. 

\subsection{Acknowledgements}
JE would first like to thank his advisor, David Nadler, for his continued support and encouragement during this project. KS thanks Rok Gregoric for illuminating discussions on the topics of homotopy theory and categorical group actions. The authors together would like to thank Pavel Safronov, Sanath Devalapurkar, Swapnil Garg, Yuji Okitani, Jeremy Taylor, Jens Eberhardt, Spencer Tamagni, Germ\'an Stefanich, Sam Gunningham, Sam Raskin, Harrison Chen, David Ben-Zvi, and Justin Hilburn for discussions related to the subject of this work. The authors thank David Nadler for discussions that led to the creation of this paper.
\subsubsection{AI Statement}
AI was used for literature search throughout this project, and for locating errors as the project neared completion. It was also used to fix some issues with LaTeX packages, and to help with typesetting.

\section{Conventions and notation}
\label{sec:conventions and notation}

\subsubsection{}
Fix an algebraically closed field $k$ of characteristic zero, once and for all.

\subsubsection{}
Given a functor $F: \scrC \to \scrD$ between $\infty$-categories,
$\im(F)$ denotes the essential image of $F$.

\subsubsection{Hom-objects and mapping spaces}
The reader will encounter $\infty$-categories in
this paper equipped with a variety of different structures,
each of which we use at different times. Here we
establish notation for the various hom-objects and
mapping spaces associated with these structures.

Let $\scrC$ be an $\infty$-category, and let $a,b \in \scrC$
be two objects. We adopt the following notation.
\begin{enumerate}
	\item \label{hom notation 1}
	$\sfMap_{\scrC}(a,b) \in \Grpd_{\infty}$ denotes the mapping
	space from $a$ to $b$. We drop the subscript ``$\scrC$" 
	when the category is clear from context.

	\item \label{hom notation 2}
	If $(\scrC, \otimes)$ is closed symmetric monoidal, then $\sHom_{\scrC}(a,b) \in \scrC$ denotes
	the internal hom-object from $a$ to $b$, i.e. the object determined by,
	\[\sfMap_{\scrC}(c, \sHom_{\scrC}(a,b)) \simeq \Map_{\scrC}(c \otimes a, b),\]
	for $c \in \scrC$. As above, we drop the subscript ``$\scrC$" when the
	category is clear from context.

	\item \label{hom notation 3} 
	If $\scrC$ is a module category over a
	rigid symmetric monoidal category $\scrR$, 
	then $\Hom_{\scrR}(a,b) \in \scrR$
	is the object determined by
	\[\sfMap_{\scrR}(V, \Hom_{\scrR}(a,b)) \simeq \sfMap_{\scrC}(V \otimes a, b),\]
	for $V \in \scrR$, $a, b \in \scrC$.

	When $\scrR$ is $R\mod$ for a ring spectrum $R$, we alternatively write
	$\Hom_R(a,b)$. For example, if $\scrC$ is an object of $\dgcat_k$, 
	then $\Hom_k(a,b)$ is the usual $k$-module spectrum of maps between $a$ and $b$.

	\item
	The space of $n$-morphisms in $\scrC$
	is given by $\sfMap_{\Cat_{\infty}}([n], \scrC)$
	where $[n]$ denotes the $n$-simplex category.
	We denote by $\Ob(\scrC)$, $\Mor(\scrC)$, and
	$\Mor_n(\scrC)$, respectively, the spaces $0$-morphisms (i.e. objects),
	$1$-morphisms, and $n$-morphisms, respectively.

	%\item \label{hom notation 4}
	%Suppose $\scrC$ has a symmetric monoidal structure,
	%and let $\scrR$ be rigid symmetric monoidal $\infty$-category as above.
	%By \cite[Corollary 3.4.1.7]{HA}, a symmetric monoidal functor
	%$\scrR \to \scrC$ is equivalent to an algebra structure on
	%$\scrC$ as an object of $\scrR\mod$, which we denote by
	%$\otimes_{\scrR}$. We let $\Hom^{\scrR}(a,b)$ denote 
	%the internal hom-object from $a$ to $b$ with respect
	%to the $\otimes_S$ product structure on $\scrC$. Note that this
	%notation does not conflict with the notation of (\labelcref{hom notation 3}).
\end{enumerate}

\begin{warning}
Given two objects $X,Y \in \Stk_k$, we will denote the internal mapping object
of $\Stk_k$ between $X$ and $Y$, i.e. the mapping stack, by
$\Map(X, Y)$.
\end{warning}

We will mainly be interested in categories of sheaves in this paper,
for which we have the following alternative notation.

Let $f: X \to S$ be a derived stack over $S$. Pullback along $f$ gives 
a natural symmetric monoidal functor $f^*:\QCoh(S) \to \QCoh(X)$,
making $\QCoh(X)$ into an algebra over $\QCoh(S)$. Composing
$f^*$ with the canonical symmetric monoidal embedding $\Upsilon: \QCoh(X) \to \IndCoh(X)$
obtains a symmetric monoidal functor $\QCoh(S) \to \IndCoh(X)$,
making $\IndCoh(X)$ into an algebra over $\QCoh(S)$, as well.

\begin{enumerate}
	
	\item We use $\sHom_X(-,-)$ to denote the hom-object
	(\labelcref{hom notation 2}) for $(\scrC, \otimes)$ taken to be either
	$(\QCoh(X), \otimes)$ or $(\IndCoh(X), \overset{!}{\otimes})$.
	When there is potential for confusion, we write $\sHom_{\QCoh(X)}(-,-)$
	or $\sHom_{\IndCoh(X)}(-,-)$.
	
	\item We use $\Hom_S(-,-)$ to denote the hom-object of 
	(\labelcref{hom notation 3}) above for $\scrR = \QCoh(S)$
	and $\scrC$ taken to be either $\QCoh(X)$ or $\IndCoh(X)$.
\end{enumerate}

\brem
There is a useful relationship between $\sHom_{\QCoh(X)}(-,-)$
and $\Hom_S(-,-)$:
\[f_*\sHom_{\QCoh(X)}(-,-) \simeq \Hom_S(-,-).\]
When $S$ is smooth, and $f$ is schematic, then
we similarly have
\[f_*\sHom_{\IndCoh(X)}(-,-) \simeq \Hom_S(-,-).\]

Note that if $X$ is taken to be a stack over itself, i.e.
$S = X$ and $f = \id_X$, then $\Hom_S(-,-) = \sHom_X(-,-)$.
\erem

\brem
In the case when the base stack $S$ is an affine scheme $\Spec R$, 
we also write $\Hom_R(-,-)$ in place of $\Hom_S(-,-)$, by slight abuse of
notation. Indeed note that an action of $\QCoh(\Spec R)$ is an action of $R\mod$,
and that $\Gamma(\Hom_{\Spec R}(-,-)) \simeq \Hom_R(-,-)$.
\erem

\subsubsection{Inclusions}
In order to minimize verbiage, we will use
$i_{?}$ as standing, all-purpose notation 
to denote the inclusion of some space (``$?$")
into another, whenever said inclusion is clear
from context.

\subsubsection{Tate invariants}
\label{sssec:Tate invariants}
Let $\scrC$ be a $k$-linear $\infty$-category equipped with an $S^1$-action.
The invariant category $\scrC^{S^1}$ is linear over $k[u]$, where $u$ has cohomological
degree $2$. We set 
\[\scrC^{tS^1} := \scrC^{S^1} \otimes_{k[u]} \ku,\] 
and call the resulting category the \emph{Tate invariants} of $\scrC$.

When $\scrC$ is a large category, the Tate invariants might not be the correct notion
for a given purpose, so we often use the following construction of Preygel instead.
Assume that $\scrC$ is a $k$-linear stable $\infty$-category with a coherent $t$-structure 
(\cite[Definition 4.2.7]{Preygel}). Let $\Coh(\scrC) \subset \scrC$ be the full subcategory
spanned by bounded almost perfect objects of $\scrC$. Then Preygel defines the \emph{$t$Tate invariants}
of $\scrC$ to be the ind-completion of the Tate invariants of the induced $S^1$-action on $\Coh(\scrC)$:
\[\scrC^{\tT} := \Ind(\Coh(\scrC)^{tS^1}).\]

\subsubsection{Stacks}
\label{sssec:stacks conventions}
We use the notions of geometric and ind-geometric stacks of \cite{CautisWilliams,CautisWilliams2}. In particular, if $X$ is a coherent ind-geometric stack with $S^1$-action, it makes sense to write
\begin{align*}
    \mathsf{IndCoh}(X)^{\tT}
\end{align*}
by \cite[Proposition 5.30]{CautisWilliams} and the results of \cite[Section 3.1]{ES-Kashiwara}. By Remark 5.20 of \textit{op. cit.}, for such stacks the categories of ind-coherent sheaves defined in \cite{CautisWilliams} and in \cite{GaitsgoryIndCoh,GR17,GR17II} are the same.

By \textit{closed immersion}, we mean a closed immersion in the sense of \cite[Section 3.7]{CautisWilliams} that is also schematic.

\subsubsection{Symmetric monoidal categories}
If $(\scrC, \otimes)$ is a symmetric monoidal category, we
denote the tensor unit by $\bbone_{\scrC}$.

\subsubsection{$(\infty,2) $-categories}
Since all $(\infty,1)$-categories considered in this paper will be presentable, we adopt the convention that $(\infty,2)$-categories will be denoted by bold font, while $(\infty,1)$-categories by ordinary font. For example, $\DGCat_k$ will be the $(\infty,1)$-category of stable presentable $k$-linear $(\infty,1)$-categories, while $\bDGCat_k$ will be the $(\infty,2)$-category of stable presentable $k$-linear $(\infty,1)$-categories. Given $(\infty,2)$-categories $\mathscr{C,D}$, the $(\infty,2)$-category of $(\infty,2)$-functors $F:\mathscr{C}\to\mathscr{D}$ will be denoted by $\bFun(\mathscr{C},\mathscr{D})$.
\subsection{Glossary of categories}
\label{ssec: glossary of categories}
${}$ \\
\begin{outline}
	\1 $\Cat_{\infty}$ denotes the $\infty$-category of $\infty$-categories.   
        \2 $\Cat_{\infty}^t$ denotes the category of stable $\infty$-categories equipped with $t$-structures, whose $1$-morphisms are $t$-exact functors,
    and whose higher morphisms are as in $\Cat_{\infty}$.
        \2 $\Cat_{\infty}^{\sm}$ denotes the full subcategory of $\Cat_{\infty}$ spanned by small, idempotent-complete
        $\infty$-categories.
        \2 $\Cat_{\infty, \st}$ denotes the subcategory of $\Cat_{\infty}$ spanned by stable $\infty$-categories,
        whose $1$-morphisms are exact functors, and whose higher morphisms are as in $\Cat_{\infty}$.
            \3 $\Cat_{\infty, \st}^{\sm}$ denotes the subcategory of $\Cat_{\infty}$ given by $\Cat_{\infty, \st} \cap \Cat_{\infty}^{\st}$.
    \1 $\bCat_\infty$ denotes the $(\infty,2)$-category of $\infty$-categories.
    \2 $\bCat_\infty^t$ denotes the $(\infty,2)$-category of $\infty$-categories equipped with $t$-structures, whose $1$-morphisms are $t$-exact functors, and whose higher morphisms are as in $\bCat_\infty$.
            \2 $\bCat_{\infty, \st}$ denotes the subcategory of $\bCat_{\infty}$ spanned by stable $\infty$-categories,
        whose $1$-morphisms are exact functors, and whose higher morphisms are as in $\bCat_{\infty}$.
            \3 $\bCat_{\infty, \st}^{\sm}$ denotes the subcategory of $\bCat_{\infty}$ given by $\bCat_{\infty, \st} \cap \bCat_{\infty}^{\st}$.
    \1 $\Pr^L_{\st}$ denotes the $\infty$-subcategory of $\Cat_{\infty}$ spanned by
       presentable, stable $\infty$-categories, whose $1$-morphisms are exact colimit-preserving functors, 
       and whose higher morphisms are as in $\Cat_{\infty}$.    
    \2 $(\Pr^L_{\st})^t$ denotes the $\infty$-category of objects in $\Pr^L_{\st}$
            equipped with $t$-structures, whose $1$-morphims are $t$-exact functors,
            and whose higher morphisms are as in $\Pr^L_{\st}$.
    \1 $\bPr^L$ denotes the $(\infty,2)$-subcategory of $\bCat_\infty$ spanned by presentable, $(\infty,1)$-categories, whose 1-morphisms are colimit preserving functors, and whose higher morphisms are as in $\bCat_\infty$.
    \2 $\bPr^L_{\st}$ denotes the $\infty$-subcategory of $\bCat_{\infty}$ spanned by
       presentable, stable $\infty$-categories, whose $1$-morphisms are exact colimit-preserving functors, and whose higher morphisms are as in $\bCat_\infty$.
    \2 $(\bPr^L_{\st})^t$ denotes the $(\infty,2)$-category of objects in $\bPr^L_{\st}$
            equipped with $t$-structures, whose $1$-morphims are $t$-exact functors,
            and whose higher morphisms are as in $\bPr^L_{\st}$.
            \1 $\mathsf{Vect}_k$ denotes the (derived) $\infty$-category of dg vector spaces over $k$.
            \2 $\mathsf{CAlg}_k$ denotes the (derived) $\infty$-category of commutative dg algebras over $k$.
            \2 For $A\in\mathsf{CAlg}_k$, $A\mod$ denotes the (derived) $\infty$-category of $A$-modules.
	%\2 $\Cat_{\infty}^{\mathrm{st},t}$ denotes the $\infty$-category of stable $\infty$-categories with $t$-structures and $t$-exact functors.
	\1 $\DGCat_k$ denotes $\Mod_{\Vect_k}(\Pr^L_{\st})$.
        \2 $\DGCat _k^t$ denotes the $\infty$-category of objects in $\DGCat_k$ equipped with
        $t$-structures, whose $1$-morphisms are $t$-exact functors, and whose higher morphisms are 
        as in $\DGCat_k$. Equivalently, $\DGCat_k^t \simeq \Mod_{(\Vect_k,t)}((\Pr^L_{\st})^t)$,
        where $(\Vect_k, t)$ denotes $\Vect_k$ equipped with its standard $t$-structure.
        There is an obvious forgetful functor $\oblv: \DGCat_k^t \to \DGCat_k$ which
        forgets the $t$-structure.
	\2 $\DGCat_k^\otimes $ denotes the $\infty$-operad whose underlying $\infty$-category is $\DGCat_k$ (cf. \cite[Notation 4.8.1.2]{HA}).
    \2 For $A\in \mathsf{CAlg}_k\coloneq \mathsf{CAlg}(\Vect_k)$, $\DGCat_A\coloneq \Mod_{A\mod}(\DGCat_k)$
    \1 $\bDGCat_k$ denotes the $(\infty,2)$-category of $\mathsf{Vect}_k$-module categories in $\bPr^L_{\text{st}}$.
    \2 We define the $(\infty,2)$-categories $\bDGCat_k^t$ and $\bDGCat_A$ by including the non-invertible 2-morphisms in the analogous $(\infty,1)$-categories defined above.
%	\1 $\DGCat_k$ denotes the $(\infty,1)$-category of modules in $\Cat_{\infty}$ over
%	$\Perf_k$ which are stable, and for which the action of $\Perf_k$ is exact in
%	both variables.
	\1 $\Ch_k$ denote the dg-category of strict chain complexes over the field $k$. There is
	a localization map $\Ch_k \to \Vect_k$ inverting quasi-isomorphisms.
	\1 $\QCoh(X)$ denotes the $k$-linear (stable cocomplete) 
	$\infty$-category of quasicoherent complexes on a derived prestack $X$. 
	It is equipped with a natural t-structure, whose heart $\QCoh(X)^{\heart}$ 
	is equivalent to the (ordinary) category of quasicoherent sheaves on $\pi_0(X)$.
	\2 $\Perf(X) \subset \QCoh(X)$ is the full subcategory spanned
	by perfect complexes. If $X$ is perfect (e.g. a quasi-compact, quasi-separated
	derived scheme) in the sense of \cite{BZFN}, then $\QCoh(X) = \Ind(\Perf(X))$.
    \2 Let $X$ be a reasonable ind-geometric stack in the sense of \cite[Definitions 4.1-4.2]{CautisWilliams}. Then we define $\mathsf{Coh}(X)$ as in Section 5.1 of \textit{op. cit.} When $X$ is furthermore coherent and locally almost of finite type, this coincides with the definition of \cite[Section 10]{GaitsgoryIndCoh} and \cite[Chapter 5]{GR17} (see \cite[Remark 5.20]{CautisWilliams}).
	%\1 Given a derived stack $X$, $\Sing(X)$ denotes the the quotient dg-category
	%$\Coh(X)/\Perf(X)$.
	%\2 $\Sing^{\infty}(X)$ denotes the quotient of large dg-categories $\IndCoh(X)/\QCoh(X)$. 
	%It is readily seen that $\Ind(\Sing(X)) \simeq \Sing^{\infty}(X)$.
	\1 For laft prestacks $X$, $\IndCoh(X)$ is defined in detail in \cite{GaitsgoryIndCoh}, but for $X$ as in \zcref{sssec:stacks conventions}, it is simply the ind-completion of $\Coh(X)$ defined above.
	%\1 For $A$ a graded $k$-algebra, $A\mod_{\gr}$ denotes the 
	%$\Rep(\G_m)$-module category of graded modules over $A$.
	%\2 $\Perf_{\gr} A$ denotes the full subcategory of $A\mod_{\gr}$
	%spanned by compact objects.
	\1 $\mathsf{Spc}$ denotes the $\infty$-category of $\infty$-groupoids, whose
	objects are alternatively called \emph{spaces} or \emph{animae}.
	\1 Given an $\infty$-category $\scrC$ with a t-structure, $\scrC^+$ will denote $\cup_{n \geq 0} \scrC^{\geq -n}$,
	and $\scrC^-$ will denote $\cup_{n \geq 0} \scrC^{\leq n}$.
    \1 We let $\bCoh_t$ denote the $(\infty,2)$-category whose objects are $k$-linear stable presentable $\infty$-categories
	equipped with coherent $t$-structure, whose $1$-morphisms are $k$-linear, continuous, and
	left $t$-exact up to a shift functors, and whose higher morphisms are those in $\bDGCat_k^{t}$.
    \2 We let $\mathsf{Coh}_t$ denote the underlying $(\infty,1)$-category of $\bCoh_t$.
	\1 We let $\bReg_t\subset \bCoh_t$ denote the full $(\infty,2)$-subcategory of $\bCoh_t$ spanned by objects with regular $t$-structures.
    \2 We let $\mathsf{Reg}_t$ denote the underlying $(\infty,1)$-category of $\bReg_t$.

\begin{comment}
    \1 We follow \cite{ChenDhillon}. If $S$ is an affine scheme, we define the category of \textit{almost $!$-perfect} sheaves on $S$ to be the full subcategory $\mathsf{APerf}^!(S)\subset\mathsf{IndCoh}(S)$ of eventually coconnective sheaves with coherent cohomology. We then extend to stacks by descent as usual.
    \1 We follow \cite{ChenDhillon,BZChenHelmNadler}. If $Z\hookrightarrow X$ is a closed substack of a stack $X$, then we define the category of \textit{continuous ind-coherent sheaves} to be the full subcategory
\begin{align*}
	\wih{\mathsf{Coh}}(\wih{X_Z})\subset \mathsf{IndCoh}(\wih{X_Z})
\end{align*} 
of sheaves $\calF$ that are $t$-bounded and almost $!$-perfect, such that for any closed substack $i:Z'\hookrightarrow X$ set-theoretically supported on $Z$, $i^!\calF$ has coherent cohomology.
\end{comment}
\end{outline}

%--------------------------------------------------------------------------------------------------%
%--------------------------------------POTENT D-MODULES-------------------------------------------%
%--------------------------------------------------------------------------------------------------%

\section{Preliminaries}
\label{sec:preliminaries}
\subsection{Operads}
\label{ssec:operads}
In order to construct the lax symmetric monoidal structure on $\Pot$, required for the application of \cite[Theorem 5.14]{CLL_six-functor-formalisms}, we need to consider operadic versions of the intermediate categories arising in the definition of $\Pot$. Since by \cite[Proposition 4.5.5]{Preygel}, $\Pot$ is equivalent to the composition
\begin{align*}
	(\ku\otimes _{k[u]}-)\circ\mathscr{R}\circ\mathsf{IndCoh}(\mathcal{L}-)^{S^1},
\end{align*} 
and the functor $\mathscr{R}$ is a functor from $\bCoh_t\to \bReg_t$, one would like to construct a lax symmetric monoidal structure on $\mathscr{R}$. However, the categories $\bCoh_t$ and $\bReg_t$ are not closed under tensor products; hence, at best, they only define $\infty$-operads rather than symmetric monoidal categories. In this section, we execute the construction of these $\infty$-operads and prove some of their basic properties, allowing us to construct $\mathscr{R}$ as a map of $\infty$-operads.
\begin{defn}\label{defn: jointly left t-exact}
		Let $\mathscr{C}_i,\mathscr{D}\in\Mod_A(\DGCat_k)$ for $A\in\mathsf{CAlg}_k$, and let $F: \mathscr{C}_1\times \cdots \times \mathscr{C}_n\to \mathscr{D}$ be a separately continuous, exact, $A$-multilinear functor. We say $F$ is \textit{jointly left $t$-exact up to a shift} if, for some $N\in \mathbb{Z} $,
		\begin{align*}
			F(\mathscr{C}_1^{\ge 0},\ldots,\mathscr{C}_n^{\ge 0})\subset \mathscr{D}^{\ge N}.
		\end{align*} 
		Denote the category of such functors by $\Fun_A^{\text{coh}^n}(\prod \mathscr{C}_i,\mathscr{D})$.
	\end{defn} 
	We will need the following multilinear version of \cite[Lemma 4.6.1]{Preygel}.
	\begin{lem}\label{lem: multilinear regularization functor}
		Suppose $\mathscr{C}_i\in \bReg_{t,A}$ for $i\in \{1,\ldots,n\} $ and $\mathscr{D}\in\bCoh_{t,A}$. Then postcomposition with the functor $\varepsilon_{\mathscr{D}}: \mathscr{R}(\mathscr{D})\to \mathscr{D}$ of \cite[Proposition 4.3.2]{Preygel}, which is canonically $A$-linear by \zcref{lem: A-linear regularization} below, induces an equivalence 
\begin{align*}
	(\varepsilon_{\mathscr{D}})_*:\Fun_A^{\mathrm{coh}^n}\left(\prod \mathscr{C}_i,\mathscr{R}(\mathscr{D})\right)\to \Fun_A^{\mathrm{coh}^n}\left(\prod \mathscr{C}_i,\mathscr{D}\right).\footnotemark
\end{align*} 
\footnotetext{The functor $(\varepsilon_{\mathscr{D}})_*$ indeed preserves jointly left $t$-bounded functors, because $\varepsilon_{\mathscr{D}}: \mathscr{R}(\mathscr{D})^+\to \mathscr{D}^+$ is a $t$-exact equivalence.}
	\end{lem} 
	\begin{proof}
		Below, we use the convention that on small categories $A$-linearity means $\mathsf{Perf}(A)$-linearity. If $n=0$, then the statement to be shown is that we have a $A$-linear equivalence
		\begin{align*}
			\varepsilon_{\mathscr{D}}:\mathscr{R}(\mathscr{D})^+\xrightarrow{\simeq}\mathscr{D}^+.
		\end{align*} 
		Indeed, this restriction is $A$-linear by the same argument as in \zcref{lem: A-linear regularization} below, and it is an equivalence by \cite[Proposition 4.3.2]{Preygel}. Now, let $n>0$ and $F\in\Fun_A^{\text{coh}^n}(\prod \mathscr{C}_i,\mathscr{D})$. Restricting $F$ to tuples of coherent objects, it admits a lift through the $k$-linear equivalence
		\begin{align*}
			\varepsilon_{\mathscr{D}}:\mathscr{R}(\mathscr{D})^+\xrightarrow{\simeq}\mathscr{D}^+.
		\end{align*} 
		Thus, we can ind-extend the resulting functor in each variable to obtain a lift $\widetilde {F}$ of $F$. Since every object of $\mathscr{C}_i^{\ge 0}$ is a filtered colimit of objects of $\mathsf{Coh}(\mathscr{C}_i)^{\ge 0}$, $\widetilde{F}$ satisfies the same bound as $F$; hence, $(\varepsilon_{\mathscr{D}})_*$ is essentially surjective. For fully faithfulness, we make the following observation.
		\begin{align*}
			\Fun_A^{L,L}(\mathscr{C}_i\times \mathscr{C}_j,\mathscr{E})\simeq&\Fun_A^{L}(\mathscr{C}_i,\Fun_A^{L}(\mathscr{C}_j,\mathscr{E}))\\
			\simeq&\Fun_A^{\text{ex}}(\mathsf{Coh}(\mathscr{C}_i),\Fun_A^L(\mathscr{C}_j,\mathscr{E}))\\
			\simeq&\Fun_A^{\text{ex}}(\mathsf{Coh}(\mathscr{C}_i),\Fun_A^{\text{ex}}(\mathsf{Coh}(\mathscr{C}_j),\mathscr{E}))\\
			\simeq&\Fun_A^{\text{ex,ex}}(\mathsf{Coh}(\mathscr{C}_i)\times \mathsf{Coh}(\mathscr{C}_j),\mathscr{E}),
		\end{align*} 
where in the second and third lines we use the universal property of $\mathsf{Ind}(-)$. Note that $\mathsf{Ind}(-)$ extends the $\mathsf{Perf}(A)$-linearity to $\mathsf{Mod}_A$-linearity. Thus, we have a natural equivalence
\begin{align}\label{eq: multilinear ind universal prop}
	\Fun_A^{L^n}\left(\prod \mathscr{C}_i,\mathscr{E}\right)\simeq\Fun_A^{\text{ex}^n}\left(\prod \mathsf{Coh}(\mathscr{C}_i),\mathscr{E}\right),
\end{align} 
where $\Fun_A^{L^n}$ (resp. $\Fun_A^{\text{ex}^n}$) denotes separately continuous (resp. exact) multilinear functors. We next check that this equivalence restricts to jointly left $t$-exact up to a shift functors. Clearly, restriction preserves this condition, so it remains to check that if $F$ is a jointly left $t$-exact up to a shift functor in the right-hand category, then it has an extension which is also jointly left $t$-exact up to a shift. Thus, suppose $F\in \Fun_A^{\text{ex}^n}\left(\prod \mathsf{Coh}(\mathscr{C}_i),\mathscr{E}\right)$ such that
\begin{align*}
	F(\mathsf{Coh}(\mathscr{C}_1)^{\ge 0},\ldots,\mathsf{Coh}(\mathscr{C}_n)^{\ge 0})\subset \mathscr{E}^{\ge N}
\end{align*} 
for some $N\in \mathbb{Z} $. Since $\mathscr{C}_i^{\ge 0}\simeq \mathsf{Ind}(\mathsf{Coh}(\mathscr{C}_i)^{\ge 0})$, the functor $\widetilde{F}:\prod \mathscr{C}_i\to \mathscr{E}$ corresponding to $F$ under the equivalence \zcref[noname]{eq: multilinear ind universal prop} evaluates on tuples of nonnegative objects as a filtered colimit of values on tuples of nonegative coherent objects. Since $\mathscr{E}^{\ge N}$ is closed under filtered colimits, we conclude that $\widetilde{F}$ satisfies the same bound as $F$. We now apply this equivalence for $\mathscr{E}= \mathscr{R}(\mathscr{D})$:
\begin{align*}
	\Fun_A^{\text{coh}^n}\left(\prod\mathscr{C}_i,\mathscr{R}(\mathscr{D})\right)\simeq&\Fun_A^{\text{ex}^n,\text{lb}}\left(\prod \mathsf{Coh}(\mathscr{C}_i),\mathscr{R}(\mathscr{D})^+)\right)\\
	\xrightarrow{\simeq}&\Fun_A^{\text{ex}^n,\text{lb}}\left(\prod \mathsf{Coh}(\mathscr{C}_i),\mathscr{D}^+\right)\\
	\simeq&\Fun_A^{\text{coh}^n}\left( \prod \mathscr{C}_i,\mathscr{D} \right) ,
\end{align*} 
where the equivalence in the second line is given by $(\varepsilon_{\mathscr{D}})_*$.
	\end{proof}

 \begin{lem}\label{lem: A-linear regularization}
	Let $\mathscr{C}$ be a stable presentable category with coherent $t$-structure, and suppose that $\mathscr{C}$ is $A$-linear, for $A\in \mathsf{CAlg}_k$. Then both the category $\mathscr{R}(\mathscr{C})$ and the functor $\varepsilon_{\mathscr{C}}: \mathscr{R}(\mathscr{C})\to \mathscr{C}$ of \cite[Proposition 4.3.2]{Preygel} are canonically $A$-linear.
\end{lem} 	
\begin{proof}
	For $\mathscr{R}(\mathscr{C})$, there is an induced $\mathsf{Perf}(A)$-action on $\mathsf{Coh}(\mathscr{C})$; hence, there is a canonical $A\mod$ action on $\mathscr{R}(\mathscr{C})\coloneq \mathsf{Ind}(\mathsf{Coh}(\mathscr{C}))$ given by ind-extension. For $A$-linearity of $\varepsilon_{\mathscr{C}}$, the canonical $\mathsf{Perf}(A)$-action on $\mathsf{Coh}(\mathscr{C})$ induces the $A$-linear structure on $\varepsilon_{\mathscr{C}}$. Indeed, let $M\in A\mod$ and $c\in \mathscr{R}(\mathscr{C})$. Write $M\simeq\colim_n M_n$ for $M_n\in \mathsf{Perf}(A)$ and $c\simeq\colim_m c_m$ for $c_m\in \mathsf{Coh}(\mathscr{C})$. Since $\varepsilon_{\mathscr{C}}$ is continuous by \cite[Proposition 4.3.2]{Preygel}, the image of $M\otimes c$ under $\varepsilon_{\mathscr{C}}$ is given by
	\begin{align*}
		\varepsilon_{\mathscr{C}}(M\otimes c)\simeq& \colim_{n,m}\varepsilon_{\mathscr{C}}(M_n\otimes c_m)\\
		\simeq& \colim_{n,m}M_n\otimes c_m \\
		\simeq&  \colim_n M_n\otimes (\colim_m \varepsilon_{\mathscr{C}}(c_m))\\
		\simeq& M\otimes \varepsilon_{\mathscr{C}}(c),
	\end{align*} 
	where we note that the $A$-action on $\mathscr{R}(\mathscr{C})$ commutes with colimits in $\mathsf{Coh}(\mathscr{C})$ by construction.
\end{proof}
	
\begin{comment}
\begin{prop}\label{prop: regularization is lax monoidal}
	The restriction of the regularization functor $\mathscr{R}:\bCoh_t\to \DGCat^{t}$ to the image of $\mathsf{IndCoh}(\mathcal{L}-)^{S^1}:\mathsf{indGStk}_{\mathrm{ft}}\to \bCoh_t$ is lax symmetric monoidal.
\end{prop} 
\begin{proof}
	By \zcref{lem: multilinear regularization functor}, there is a symmetric monoidal structure on regularization. Combining this with \zcref{lem: exterior product in Coh_t}, we conclude.
\end{proof} 
	
\begin{prop}\label{prop: tTate is lax symmetric monoidal}
	The restriction of the functor
	\begin{align*}
		(-)^{\tT}:\bCoh_t\to \bDGCat_{\ku}
	\end{align*} 
	to the image of $\mathsf{IndCoh}(\mathcal{L}-):\mathsf{indGStk}_{\mathrm{ft}}\to \bCoh_t$ is lax symmetric monoidal.
\end{prop} 
\begin{proof}
	By \zcref{lem: S^1-invariants commute with tensors} and \zcref{prop: regularization is lax monoidal}, it suffices to check that extension of scalars is lax symmetric monoidal. Indeed, it is symmetric monoidal.
\end{proof} 
\end{comment}  
\begin{con}\label{con: Coh_t and Reg_t operads}
	Let $ \Coh_{t,A}^\otimes $ be the following category. We choose a simplicial model of the total category of the $\infty$-operad $\DGCat_A^\otimes \coloneq \mathsf{Mod}_A(\Pr^L_{\mathrm{st}})^{\otimes } $ over $\mathsf{Fin}_*$ (see \cite[Theorem 4.5.2.1]{HA} for the existence of this operad). We replace each tuple $(\mathscr{C}_1,\ldots,\mathscr{C}_n)$ by all possible tuples $((\mathscr{C}_1,t_1),\ldots,(\mathscr{C}_n,t_n))$ (including the empty tuple), where the $t _i$ are coherent $t$-structures on the $\mathscr{C}_i$. Between two such objects, we retain only those mapping space components whose component $A$-multilinear functors are jointly left $t$-exact up to a shift, and then we take the simplicial nerve. The simplicial set we have obtained still satisfies the weak Kan complex property: If inner horns have edges that are jointly left $t$-exact up to a shift, they are filled by outer edges in $\DGCat_A^\otimes $ which must be jointly left $t$-exact up to a shift. This category is equipped with the obvious forgetful functor
	\begin{align*}
		\oblv_t^\otimes :\Coh_{t,A}^\otimes \to \DGCat_A^\otimes ,
	\end{align*} 
	given by forgetting the $t$-structure. In the same way, we can construct the category $\Reg_{t,A}^\otimes $.
\end{con}  
\begin{rem}
	The mapping space in $\Coh_{t,A}^\otimes $ over a map $\alpha:\braket{n}\to \braket{m}$ in $N(\mathsf{Fin}_*)$ is 
	\begin{align*}
		\mathsf{Map}^\alpha_{\Coh_{t,A}^\otimes }((\mathscr{C}_i),(\mathscr{D}_j))\coloneq \prod_{j=1}^m\Fun_A^{\text{coh}^{|\alpha^{-1}(j)|}}\left( \prod_{i\in\alpha ^{-1}(j)}\mathscr{C}_i,\mathscr{D}_j \right) ^\simeq.
	\end{align*}
\end{rem} 
\begin{prop}\label{prop: Coh_t and Reg_t are operads}
	The composition 
	\begin{align*}
		\pi:\Coh_{t,A}^\otimes \xrightarrow{\oblv_t^\otimes }\DGCat_A^\otimes \to N(\mathsf{Fin}_*)
	\end{align*} 
	exhibits $\Coh_{t,A}^\otimes $ as an $\infty$-operad. Similarly, the composition
	\begin{align*}
		\Reg_{t,A}^\otimes \hook \mathsf{Coh}_{t,A}^\otimes\xrightarrow{\oblv_t^\otimes }\DGCat_A^\otimes \to N(\mathsf{Fin}_*)
	\end{align*} 
	exhibits $\Reg_{t,A}^\otimes $ as an $\infty$-operad.
\end{prop} 
\begin{proof}
	We will prove the claim for $\Coh_{t,A}^\otimes $; the proof for $\Reg_{t,A}^\otimes $ is the same. We show that the axioms (1-3) of \cite[Definition 2.1.1.10]{HA} hold for $\pi$. Throughout this proof, we omit the $t$-structures for ease of notation.
	\begin{enumerate}[label=\arabic*)]
		\item Let $\alpha:\braket{m}\to \braket{n}$ be an inert morphism in $N(\mathsf{Fin}_*)$. Let $\alpha^{-1}(j)=\{i_j\} $, and write
			\begin{align*}
				\mathscr{C}'=(\mathscr{C}_{i_1},\ldots,\mathscr{C}_{i_n})\in (\Coh_{t,A}^\otimes )_{\braket{n}}
			\end{align*} 
			Now, we define an arrow $F: \mathscr{C}\to \mathscr{C}'$ over $\alpha$ whose component arrows are identity functors. Explicitly, we let $\mathscr{C}\coloneq (\mathscr{C}_1,\ldots,\mathscr{C}_m)$, and we set the component functor of $F$ for output $j\in \{1,\ldots,n\} $ to be given by the identity functor
			\begin{align*}
				F_j\coloneq \id_{\mathscr{C}_{i_j}}: \prod_{i\in\alpha ^{-1}(j)}\mathscr{C}_i=\mathscr{C}_{i_j}\to \mathscr{C}_{i_j}.
			\end{align*} 
			We claim that $F$ is $\pi$-coCartesian. Note that $ \pi$ is an inner fibration by \cite[Proposition 2.3.1.5]{HTT}. Thus, by the dual of Proposition 2.4.4.3 of \textit{op. cit.}, it suffices to show that precomposition with $F$ induces an equivalence
			\begin{align*}
				F^*: \mathsf{Map}^\beta(\mathscr{C}',\mathscr{D})\xrightarrow{\simeq}\mathsf{Map}^{\beta\alpha}(\mathscr{C},\mathscr{D})
			\end{align*} 
			for every $\mathscr{D}=(\mathscr{D}_1,\ldots,\mathscr{D}_r)$ and every $\beta:\braket n\to \braket r$, not necessarily inert. Now, the source and target of $F^*$ are respectively
			\begin{align*}
				\prod_{k=1}^r\Fun_A^{\text{coh}^{|\beta^{-1}(k)|}}\left( \prod_{j\in\beta ^{-1}(k)} \mathscr{C}_{i_j} ,\mathscr{D}_k\right)^\simeq 
			\end{align*} 
			and
			\begin{align*}
				\prod_{k=1}^r\Fun_A^{\text{coh}^{|(\beta\alpha)^{-1}(k)|}}\left( \prod_{i\in(\beta\alpha)^{-1}(k)}\mathscr{C}_i,\mathscr{D}_k \right) ^\simeq.
			\end{align*} 
			For each $k>0$, the map
			\begin{align*}
				\beta^{-1}(k)\to &(\beta\alpha)^{-1}(k)\\
				j\mapsto &i_j
			\end{align*} 
			is a bijection. Thus, the functor
			\begin{align*}
				r_k:\prod_{i\in(\beta\alpha)^{-1}(k)}\mathscr{C}_i\xrightarrow{\simeq}&\prod_{j\in\beta ^{-1}(k)}\mathscr{C}_{i_j}\\
				(x_i)_i\mapsto &(x_{i_j})_j
			\end{align*} 
			is an equivalence. Moreover, given a functor
			\begin{align*}
				G_k:\prod_{j\in\beta^{-1}(k)}\mathscr{C}_{i_j}\to \mathscr{D}_k,
			\end{align*} 
			we have
			\begin{align*}
				(G\circ F)_k:\prod_{i\in(\beta\alpha)^{-1}(k)}\mathscr{C}_i\to &\mathscr{D}_k\\
				(x_i)_i\mapsto &G_k\left((F_j(x_{i_j}))_{j\in\beta ^{-1}(k)}\right)=G_k\left( (x_{i_j})_{j\in\beta ^{-1}(k)} \right) .
			\end{align*} 
			Thus, we have shown that
			\begin{align*}
				(G\circ F)_k\simeq G_k\circ r_k,
			\end{align*} 
			so precomposition with $F$ induces a reindexing equivalence, as desired.
		\item Now, let $\alpha:\braket{m}\to \braket{n}$ be an arbitrary morphism in $N(\mathsf{Fin}_*)$. Choose projections
			\begin{align*}
				p_j: \mathscr{C}'\to \mathscr{D}_j
			\end{align*} 
			over $\rho_j:\braket{n}\to \braket{1}$. By the already proved item (1), these projections are coCartesian. We would like to show that the canonical map
			\begin{align*}
				\mathsf{Map}^\alpha(\mathscr{C},\mathscr{C}')\to \prod_{j=1}^n \mathsf{Map}^{\rho_j\alpha}(\mathscr{C},\mathscr{D}_j)
			\end{align*} 
			is an equivalence. Indeed, the source is simply
			\begin{align*}
				\prod_{j=1}^n\Fun_A^{\text{coh}^{|\alpha^{-1}(j)|}}\left( \prod_{i\in\alpha^{-1}(j)}\mathscr{C}_i,\mathscr{D}_j \right) ^\simeq,
			\end{align*} 
			and the $j^{\text{th}}$ factor of the target is the $j^{\text{th}}$ factor of this product. Postcomposition with $p_j$ projects precisely onto that factor, giving the desired equivalence.
		\item Given objects $\mathscr{C}_1,\ldots,\mathscr{C}_n\in\bCoh_{t,A}$, the tuple $\mathscr{C}\coloneq (\mathscr{C}_1,\ldots,\mathscr{C}_n)$ is an object over $\braket{n} $, and the projections to the $\mathscr{C}_i$ are coCartesian by (1). When $n=0$ we can use the empty tuple.
	\end{enumerate}
\end{proof}
\begin{con}\label{con: S^1-equivariant Coh_t operad}
	We record for future use the following construction, due to Lurie (see \cite[Remark 2.2.5.4]{HA}). Let $\mathscr{O}^\otimes $ be an $\infty$-operad, and let $K\in \mathsf{Spc}$. Then we define an $\infty$-operad
	\begin{align*}
		\Fun(K,\mathscr{O})^\otimes \coloneq \Fun(K,\mathscr{O}^\otimes )\times _{\Fun(K,N(\mathsf{Fin}_*)}N(\mathsf{Fin}_*).
	\end{align*} 
	In particular, there is an $\infty$-operad $\Fun(BS^1,\Coh_t)^\otimes $. Its $n$-ary operations are precisely $S^1$-equivariant $n$-ary operations in $\Coh_t ^\otimes $. 
\end{con} 
\begin{con}\label{con: construction of P}
	We also construct the $\infty$-operad $\mathscr{P}^\otimes $ as the following arrow operad. Let
	\begin{align*}
		\mathscr{Q}^\otimes \coloneq \Fun(\Delta^1,\Coh_{t,A}^\otimes )\times _{\Fun(\Delta^1,N(\mathsf{Fin}_*))}N(\mathsf{Fin}_*),
	\end{align*} 
	be the arrow $\infty$-operad of \cite[Example 4.3(d)]{LinskensNardinPol} with $\mathscr{I}\coloneq \Delta^1$. Let $\mathscr{P}\subset \mathscr{Q} $ be the full subcategory of the underlying category of $\mathscr{Q}^\otimes $ spanned by objects $\varepsilon_{\mathscr{C}}$, where for $\mathscr{C}\in\bCoh_{t,A}$,
	\begin{align*}
		\varepsilon_{\mathscr{C}}: \mathscr{R}(\mathscr{C})\to \mathscr{C}
	\end{align*} 
	is the canonical functor. By the discussion in the beginning of \cite[Section 2.2.1]{HA}, the category $\mathscr{P}^\otimes $ spanned by the objects $P=\bigoplus _i\varepsilon_{\mathscr{C}_i}$ for $\mathscr{C}_i\in \bCoh_{t,A}$ is an $\infty$-operad.
\end{con} 
\begin{rem}
	The space of operations in $\mathscr{P}^\otimes $ between tuples $(\varepsilon_{\mathscr{C}_1},\ldots,\varepsilon_{\mathscr{C}_n})$ and $\varepsilon_{\mathscr{D}}$ is
	\begin{align*}
		\Fun^{\text{coh}^n}\left( \prod \mathscr{R}(\mathscr{C}_i),\mathscr{R}(\mathscr{D}) \right)^\simeq \underset{\Fun^{\text{coh}^n}\left( \prod \mathscr{R}(\mathscr{C}_i),\mathscr{D} \right)^\simeq }{\times }\Fun^{\text{coh}^n}\left(\prod\mathscr{C}_i,\mathscr{D}\right)^\simeq.
	\end{align*} 
\end{rem} 
\begin{lem}\label{lem: P is equivalent to Coh_t}
	The functor given by the composite
	\begin{align*}
		p:\mathscr{P}^\otimes \to \Fun(\Delta^1,\Coh_{t,A}^\otimes )\xrightarrow{\ev_1} \Coh_{t,A}^\otimes ,
	\end{align*} 
	where the first functor is projection and the second the target functor, is an equivalence of $\infty$-operads.
\end{lem} 
\begin{proof}
	We first show that $p$ is a map over $N(\mathsf{Fin}_*)$: There is a diagram where the inner squares commute,
	\[
	\begin{tikzcd}
		\mathscr{P}^\otimes \arrow[r]\arrow[d]&\Fun(\Delta^1,\Coh_{t,A}^\otimes )\arrow[r]\arrow[d]&\Coh_{t,A}^\otimes \arrow[d]\\
		N(\mathsf{Fin}_*)\arrow[r]&\Fun(\Delta^1,N(\mathsf{Fin}_*))\arrow[r]&N(\mathsf{Fin}_*)
	\end{tikzcd}.
	\] 
	Thus, it is sufficient to show that $\mathscr{P}^\otimes \to \Coh_{t,A}^\otimes $ is an equivalence of categories; preservation of inert morphisms then follows formally, so that it is also an equivalence of $\infty$-operads. Essential surjectivity follows from \cite[Proposition 4.3.2]{Preygel}. For fully faithfulness, let $\alpha:\braket{m}\to \braket{n}$ be an arrow in $N(\mathsf{Fin}_*)$. Recall that an object of $\mathscr{P}^\otimes_{\braket{m}} $ is a tuple 
	\begin{align*}
		\varepsilon_{\mathscr{C}}\coloneq (\varepsilon_{\mathscr{C}_1},\ldots,\varepsilon_{\mathscr{C}_m})
	\end{align*} 
	where the $\mathscr{C}_i\in\bCoh_{t,A}$ and $\varepsilon_{\mathscr{C}_i}: \mathscr{R}(\mathscr{C}_i)\to \mathscr{C}_i$ are the maps of \textit{loc. cit.} Then, an arrow lying over $\alpha$ in $\mathscr{P}^\otimes $ consists of arrows
	\begin{align*}
		F:& \mathscr{C}\to \mathscr{D}\coloneq (\mathscr{D}_1,\ldots,\mathscr{D}_n)\\
		\widetilde{F}:&\mathscr{RC}\to \mathscr{RD},
	\end{align*} 
	where $\mathscr{RC}$ (resp. $\mathscr{RD}$) denotes $(\mathscr{R}(\mathscr{C}_1),\ldots,\mathscr{R}(\mathscr{C}_m))$ (resp. $(\mathscr{R}(\mathscr{D}_1),\ldots,\mathscr{R}(\mathscr{D}_n))$), such that
	\begin{align*}
		\varepsilon_{\mathscr{D}}\circ\widetilde{F}\simeq F\circ \varepsilon_{\mathscr{C}}.
	\end{align*} 
	Thus, we have a Cartesian square
	\[
	\begin{tikzcd}
		\mathsf{Map}^\alpha_{\mathscr{P}^\otimes }(\varepsilon_{\mathscr{C}},\varepsilon_{\mathscr{D}})\arrow[r]\arrow[d,"p_*"]&\mathsf{Map}_{\Coh_{t,A}^\otimes }^\alpha(\mathscr{RC},\mathscr{RD})\arrow[d,"(\varepsilon_{\mathscr{D}})_*"]\\
		\mathsf{Map}_{\Coh_{t,A}^\otimes }^\alpha(\mathscr{C},\mathscr{D})\arrow[r,"(\varepsilon_{\mathscr{C}})^*"]&\mathsf{Map}^\alpha_{\Coh_{t,A}^\otimes }(\mathscr{RC},\mathscr{D})
	\end{tikzcd}.
	\] 
	Now, \zcref{lem: multilinear regularization functor} implies that the right vertical arrow $(\varepsilon_{\mathscr{D}})_*$ is an equivalence. Since the square is Cartesian, this implies that the left vertical arrow $p_*$ is an equivalence, too. Taking the union over all $\alpha$ shows that $p$ is fully faithful on the total categories.
\end{proof}	

\begin{lem}\label{lem: P -> Reg_t is a map of operads}
	The functor given by the projection
	\begin{align*}
		\mathscr{P}^\otimes \to \Fun(\Delta^1,\Coh_{t,A}^\otimes )\xrightarrow{\ev_0}\Coh_{t,A}^\otimes  
	\end{align*} 
	factors through $\Reg_{t,A}^\otimes $. Moreover, the functor $q: \mathscr{P}^\otimes \to \Reg_{t,A}^\otimes $ thus obtained is a map of $\infty$-operads.
\end{lem}
\begin{proof}
	The desired factorization follows from the definition of $\mathscr{P}^\otimes$. As in the proof of \zcref{lem: P is equivalent to Coh_t}, the map $q$ is a map over $N(\mathsf{Fin}_*)$. We now check preservation of inert morphisms. Let $\alpha:\braket{m}\to \braket{n}$ be an inert morphism, and let $P\coloneq (\varepsilon_{\mathscr{C}_1},\ldots,\varepsilon_{\mathscr{C}_m})$ be an object of $\mathscr{P}^\otimes _{\braket{m}}$. Writing $\alpha^{-1}(j)=\{i_j\} $, there is a canonical arrow over $\alpha$,
	\begin{align*}
		e_\alpha: P\to P_\alpha\coloneq (\varepsilon_{\mathscr{C}_{i_1}},\ldots,\varepsilon_{\mathscr{C}_{i_n}}),
	\end{align*} 
whose $j^{\text{th}}$ component is the square
\[
\begin{tikzcd}
	\mathscr{R}(\mathscr{C}_{i_j})\arrow[r,"\id"]\arrow[d,"\varepsilon_{\mathscr{C}_{i_j}}"]&\mathscr{R}(\mathscr{C}_{i_j})\arrow[d,"\varepsilon_{\mathscr{C}_{i_j}}"]\\
	\mathscr{C}_{i_j}\arrow[r,"\id"]&\mathscr{C}_{i_j}
\end{tikzcd}.
\] 
Moreover, $e_\alpha$ is inert, since its image under $p: \mathscr{P}^\otimes \to \Coh_{t,A}^\otimes $ is the canonical inert projection $\Coh_{t,A}^\otimes $ and $p$ is an equivalence of $\infty$-operads by \zcref{lem: P is equivalent to Coh_t}. It remains to show that an arbitrary inert arrow
\begin{align*}
	f:P\to P'
\end{align*} 
over $\alpha$ has $q(f)$ inert. Since both $e_\alpha$ and $f$ are coCartesian lifts of $\alpha$ with the same source, the dual version of \cite[Remark 2.4.1.9]{HTT} implies that there is an equivalence
\begin{align*}
	\xi:P_\alpha \xrightarrow{\simeq}P'
\end{align*} 
over $\id_{\braket{n}}$ such that $f\simeq \xi\circ e_\alpha$. Applying the functor  $q$, we obtain
\begin{align*}
	q(f)\simeq q(\xi)\circ q(e_\alpha).
\end{align*} 
Since $q(\xi)$ is an equivalence and $q(e_\alpha)$ is inert by the previous, we conclude that $q(f)$ is inert.
\end{proof} 
\begin{rem}
	The reason that the target map $\mathscr{P}^\otimes \to \Coh_{t,A}^\otimes $ is an equivalence is essentially because of \cite[Lemma 4.6.1]{Preygel}, of which \zcref{lem: multilinear regularization functor} is a multilinear generalization. On the other hand, the source map $\mathscr{P}^\otimes \to \Reg_{t,A}^\otimes $ is not an equivalence, as, if $\mathscr{C},\mathscr{D}\in\bCoh_{t,A}$ and $\widetilde{F}: \mathscr{R}(\mathscr{C})\to \mathscr{R}(\mathscr{D})$, then there is in general no way to obtain a functor $F: \mathscr{C}\to \mathscr{D}$ which makes the following diagram commute.
	\[
	\begin{tikzcd}
		\mathscr{R}(\mathscr{C})\arrow[r,"\widetilde{F}"]\arrow[d,"\varepsilon_{\mathscr{C}}"]&\mathscr{R}(\mathscr{D})\arrow[d,"\varepsilon_{\mathscr{D}}"]\\
		\mathscr{C}\arrow[r,"F"]&\mathscr{D}
	\end{tikzcd}.
	\] 
\end{rem} 
\begin{rem}
	The composite of maps of $\infty$-operads
	\begin{align*}
	\mathscr{R}_A^\otimes \coloneq qp^{-1}:\Coh_{t,A}^\otimes \to \Reg_{t,A}^\otimes 
	\end{align*} 
	of \zcref{lem: P is equivalent to Coh_t} and \zcref{lem: P -> Reg_t is a map of operads} is an operadic and $A$-linear version of the regularization functor of \cite[Proposition 4.3.2]{Preygel}.
\end{rem}

\subsection{Six-Functor Formalisms}
In this section, we review the theory of three-functor and six-functor formalisms. Since we need the more recent, $(\infty,2)$-functorial version of \cite{CLL_six-functor-formalisms}, we relate it to the usual story of \cite{Mann}. We also formulate the six-functor formalism $\mathsf{IndCoh}$ of \cite{GR17II} in the \cite{CLL_six-functor-formalisms} formalism.

We begin with the following definition of 3-functor formalisms, following \cite{Mann}.
\begin{defn}(\cite[Definition 2.4]{SixFunctors})\label{def: scholze 6FF}
	Let $\scrC$ be an $\infty$-category admitting finite limits, and equip $\mathscr{C}$ with the Cartesian monoidal structure. Let $E$ be a class of
	morphisms in $\scrC$ stable under pullback and composition and containing all isomorphisms.
	A 3-functor formalism $D$ is a lax symmetric monoidal functor
	\[D: \Corr(\scrC,E) \to \Cat_{\infty}.\]
	A 3-functor formalism is furthermore called a 6-functor formalism if
	the functors $- \otimes A$, $f^*$, and $f_!$ all have right adjoints.
\end{defn}
While Mann's definition is formulated with $\Cat_{\infty}$ as the target
category, the theory works just as well if we change the target of the
lax symmetric monoidal functor to a different symmetric monoidal $\infty$-category,
such as $\DGCat_k$ (\cite{CLL_six-functor-formalisms}).
\subsubsection{Ind-coherent sheaves}
Our main example of interest in this paper is the 6-functor formalism $\mathsf{IndCoh}^{\text{GR}}$ constructed in \cite[Theorem 5.4.3]{GR17II} (cf. also \cite[Theorem 2.4.2]{BeraldoDiffOp} and \cite[\S8.7]{SixFunctors}). We can restrict $\mathsf{IndCoh}^{\text{GR}}$ to the full subcategory $\mathsf{PreStk}_{\text{laft}}$, to obtain a functor
\begin{align*}
	\mathsf{IndCoh}^{\text{GR,res}}: \mathsf{PreStk}_{\text{laft}}^{\text{op}}\to \bDGCat_k.
\end{align*} 
\begin{prop}\label{prop: indcoh is a 6-functor formalism}
	The functor $\mathsf{IndCoh}^{\mathrm{GR,res}}$ induces a lax symmetric monoidal $(\infty,2)$-functor
	\begin{align*}
	\mathsf{IndCoh}: \bSpan(\mathsf{PreStk}_{\mathrm{laft}},E)_{P,I}\to \bDGCat_k,
	\end{align*} 
where $E=I$ is the class of ind-inf-schematic ind-proper maps and $P\subset E$ is the class of equivalences.
\end{prop} 
\begin{proof}
	We seek to apply \cite[Theorem 5.14]{CLL_six-functor-formalisms}. We first show that $ I$ (resp. $P$) is left cancellable, which means that if $g\circ f,g\in I$ (resp. $\in P$), then $f\in I$ (resp. $\in P$) as well. For $P$, this follows by two-out-of-three for equivalences. For $I$, suppose that $f:X\to Y$, $g:Y\to Z$ such that if $f,g\circ f\in I$. Factor $f$ as
	\begin{align*}
		X\to X\times _ZY\to Y,
	\end{align*} 
	where the second map is the projection. The first map is a base change of $\Delta_g$ and the second is a base change of $gf$. Hence, it suffices to check that $\Delta_g \in I$, which follows from \zcref{lem: relative diagonals of maps in I are in I} below. Now, since $P$ is the class of equivalences, it follows (cf. \cite[Example 2.7]{CLL_six-functor-formalisms}) that it suffices to check that $\mathsf{IndCoh}^{\text{GR,res}}$ is left $I$-adjointable and that it has a lax symmetric monoidal structure satisfying condition (2a) of Definition 5.13 of \textit{op. cit.} (condition (2b) is trivially satisfied for any lax symmetric monoidal structure). The lax symmetric monoidal structure in question is supplied by \cite[Chapter 3 Section 6.3.1]{GR17II}.\footnote{A right-lax symmetric monoidal structure in the sense of \textit{loc. cit.} is a lax symmetric monoidal structure in the sense of \cite[Definition 5.12]{CLL_six-functor-formalisms}.} Similarly, left $I$-adjointability is supplied by \cite[Chapter 3 Section 5.4.4]{GR17II}. For \cite[Definition 5.13(2a)]{CLL_six-functor-formalisms}, let $Y\in \mathsf{PreStk}_{\text{laft}}$ and $i:X\to X'$ be a map in $I$. Consider the commutative square
	\[
	\begin{tikzcd}
		\mathsf{IndCoh}^{\text{GR,res}}(X')\otimes \mathsf{IndCoh}^{\text{GR,res}}(Y)\arrow[r,"i^!\otimes \id_Y^!"]\arrow[d,"E_{X',Y}"]&\mathsf{IndCoh}^{\text{GR,res}}(X)\otimes \mathsf{IndCoh}^{\text{GR,res}}(Y)\arrow[d,"E_{X,Y}"]\\
		\mathsf{IndCoh}^{\text{GR,res}}(X'\times Y)\arrow[r,"(i\times \id_Y)^!"]&\mathsf{IndCoh}^{\text{GR,res}}(X\times Y)
	\end{tikzcd},
	\] 
	where $E$ denotes the exterior product functor arising from the lax symmetric monoidal structure. Both horizontal arrows admit left adjoints, given respectively by $i_*\otimes \id_{Y*}$ and $(i\times \id_Y)_*$. Denote the naturality isomorphisms by
	\begin{align*}
		\alpha:E_{X,Y}\circ (i^!\otimes \id_Y^!)\xrightarrow{\simeq}&(i\times \id_Y)^!\circ E_{X',Y}\\
		\beta:E_{X',Y}\circ(i_*\otimes \id_{Y*})\xrightarrow{\simeq}&(i\times \id_Y)_*\circ E_{X,Y}.
	\end{align*} 
We claim that the canonical left mate of $\alpha$ is given by $\beta^{-1}$. Indeed, this mate is by definition
\begin{align*}
	&(i\times \id_Y)_*\circ E_{X,Y}\xrightarrow{(i\times \id_Y)_*E_{X,Y}(\eta_{i}\otimes \eta_{\id_Y})}(i\times \id_Y)_*\circ E_{X,Y}\circ (i^!\otimes \id_Y^!)(i_*\otimes \id_{Y*})\xrightarrow{(i\times \id_Y)_*\alpha (i_*\otimes \id_{Y*})}\\
	&(i\times \id_Y)_*(i\times \id_Y)^!E_{X',Y}(i_*\otimes \id_{Y*})\xrightarrow{\varepsilon_{i\times \id_Y}E_{X',Y}(i_*\otimes \id_{Y*})}E_{X',Y}(i_*\otimes \id_{Y*}),
\end{align*} 
where for an ind-inf-schematic ind-proper map of laft prestacks $f$, $\eta_{f}$ (resp. $\varepsilon_{f}$) is the unit (resp. counit) of $f_*\dashv f^!$. Denote this mate by $m$; we will show that $\beta\circ m\simeq\id$. Now, under the adjunction $(i\times \id_Y)_*\dashv (i\times \id_Y)^!$, the map $\beta\circ m$ has adjunct
\begin{align}\label{eq: adjunct of beta composed with mate}
	\begin{split}
	&E_{X,Y}\xrightarrow{E_{X,Y}(\eta_i\otimes \eta_{\id_Y})}E_{X,Y}(i^!\otimes \id_Y^!)(i_*\otimes \id_{Y*})\xrightarrow{\alpha(i_*\otimes \id_{Y*})}\\
	&(i\times \id_Y)^!E_{X',Y}(i_*\otimes \id_{Y*})\xrightarrow{(i\times \id_Y)^!\beta}(i\times \id_Y)^!(i\times \id_Y)_*E_{X,Y}.
	\end{split}
\end{align} 
It suffices, therefore, to identify this adjunct with the adjunct to the identity, which is given by
\begin{align}\label{eq: adjunct to identity}
	E_{X,Y}\xrightarrow{\eta_{i\times \id_Y}E_{X,Y}}(i\times \id_Y)^!(i\times \id_Y)_*E_{X,Y}.
\end{align} 
Consider the two correspondences given by $\ell\coloneq (X\xleftarrow{\id_X}X\xrightarrow{i}X')$ and $r\coloneq (X'\xleftarrow{i}X\xrightarrow{\id_X}X$. Their composition is
\begin{align*}
	r\ell=(X\xrightarrow{p_1}X\times _{X'}X\xrightarrow{p_2}X),
\end{align*} 
where $p_i:X\times _{X'}X\to X$ are the canonical projections. The unit of the adjunction $\ell\dashv r$ is the $2$-morphism
\begin{align*}
	\delta:\id_X\Rightarrow r\ell
\end{align*} 
represented by the following diagram
\begin{equation}\label{cd: naturality of diagonals}
\begin{tikzcd}
	X\arrow[d,equal]&X\arrow[r,"\id_X"]\arrow[l,swap,"\id_X"]\arrow[d,"\Delta_i"]&X\arrow[d,equal]\\
	X&X\times _{X'}X\arrow[l,swap,"p_1"]\arrow[r,"p_2"]&X
\end{tikzcd},
\end{equation}
where $\Delta_i$ is the relative diagonal of $i$. Now, applying the functor $\mathsf{IndCoh}^{\text{GR}}(-)\otimes \mathsf{IndCoh}(Y)$ to $\delta$ gives the natural transformation $\eta_i\otimes \eta_{\id_Y}$. On the other hand, applying the functor $\mathsf{IndCoh}^{\text{GR}}(-\times Y)$ to $\delta$ gives the natural transformation $\eta_{i\times \id_Y}$, where in the diagram \zcref[noname]{cd: naturality of diagonals}, taking the product with $Y$ sends $\Delta_i$ to $\Delta_{i\times \id_Y}$. Further, the naturality equivalences for $r$ and $\ell$ are respectively given by $\alpha$ and $\beta$. We can now use that the exterior product $E_{-,Y}$ is a natural transformation of $(\infty,2)$-functors
\begin{align*}
	\mathsf{IndCoh}^{\text{GR}}(-)\otimes \mathsf{IndCoh}(Y)\to \mathsf{IndCoh}^{\text{GR}}(-\times Y).
\end{align*} 
By naturality, we thus obtain a commutative diagram
\[
	\begin{tikzcd}[column sep = huge]
	E_{X,Y}\arrow[d,equal]\arrow[r,"E_{X,Y}(\eta_{i}\otimes \eta_{\id_Y})"]&E_{X,Y}(i^!\otimes \id_Y^!)(i_*\otimes \id_{Y*})\arrow[d,"((i\times \id_Y)^!\beta)\circ(\alpha(i_*\otimes \id_{Y*}))"]\\
	E_{X,Y}\arrow[r,"\eta_{i\times \id_Y}E_{X,Y}"]&(i\times \id_Y)^!(i\times \id_Y)_*E_{X,Y}
\end{tikzcd},
\] 
which shows that \zcref[noname]{eq: adjunct of beta composed with mate} is equivalent to \zcref[noname]{eq: adjunct to identity}, as desired.
\end{proof} 
\begin{lem}\label{lem: relative diagonals of maps in I are in I}
	Let $f:X\to Y$ be an ind-inf-schematic ind-proper. Then the relative diagonal
	\begin{align*}
		\Delta_f:X\to X\times _YX
	\end{align*} 
	is ind-inf-schematic ind-proper.
\end{lem} 
\begin{proof}
First, assume that $f$ is ind-inf-schematic, and let $S\to X\times _YX$ be an $S$-point. We then claim that there is an isomorphism
\begin{align*}
S\underset{X\times _YX}{\times }X\simeq S\underset{X\times _YS}{\times }S.
\end{align*} 
Since $X\times _YS$ is an ind-inf-scheme and ind-inf-schemes are closed under finite limits, the above isomorphism allows us to conclude. To see the desired isomorphism, observe that there is a Cartesian square
\[
\begin{tikzcd}
	S\arrow[r]\arrow[d]&X\arrow[d]\\
	S\times _YX\arrow[r]&X\times _YX
\end{tikzcd},
\] 
where we use that $S\times _YX\simeq S\times _X(X\times _YX)$. Using this square, we see that we then have isomorphisms
\begin{align*}
	S\underset{S\times _YX}{\times }S\simeq S\underset{S\times _YX}{\times }(S\times _YX \underset{X\times _YX}{\times }X)\simeq S\underset{X\times _YX}{\times }X
\end{align*} 
as desired.

Now, suppose $f$ is ind-inf-schematic ind-proper. It suffices to check that $\Delta_f$ is ind-proper by the argument above. By definition, ind-properness can be checked at the level of the underlying reduced prestacks. Moreover, since the functor $(-)_{\mathrm{red}}$ is a right adjoint, we can assume all prestacks are reduced and classical and all limits are taken in reduced, classical prestacks.\footnote{Recall that $(-)_{\mathrm{red}}$ is by definition the right adjoint to the composition of inclusions $\mathsf{PreStk}_{\mathrm{cl,red}}\hookrightarrow \mathsf{PreStk}_{\mathrm{cl}}\hookrightarrow \mathsf{PreStk}$. Thus, if $X\times _ZY$ is a fiber product in $\mathsf{PreStk}$, then $(X\times _ZY)_{\text{red}}\simeq X_{\text{red}}\times^{\text{red}} _{Z_{\text{red}}}Y_{\text{red}}\coloneq (X_{\text{red}}\times _{Z_{\text{red}}}Y_{\text{red}})_{\text{red}}$.}
If $S\to X\times _YX$ is an $S$-point, we need to check that the projection
\begin{align*}
	S\underset{X\times _YX}{\times }X\to S
\end{align*} 
is ind-proper. Further, recall that since $f$ is ind-inf-schematic, the presentation
\begin{align*}
	S\underset{X\times _YX}{\times }X\simeq S\underset{X\times _YS}{\times }S
\end{align*} 
provides two maps $S\rightrightarrows X\times _YS$. Now, $X\times _YS\to S$ is a base change of $f$; hence, it is ind-proper. We can thus write
\begin{align*}
	X\times _YS\eqcolon P\simeq\colim_i P_i
\end{align*} 
with $P_i\to S$ proper for all $i$. Since $S$ is affine, there is some $i$ for which the two maps $S\rightrightarrows P$ factor as
\begin{align*}
	S\to P_i\to P.
\end{align*} 
Since $P_i\to P$ is a closed immersion in reduced classical prestacks, it is a monomorphism. Thus, we see that
\begin{align*}
	S\underset{X\times _YS}{\times }S\simeq S\times _{P_i}S,
\end{align*} 
which is now a fiber product of classical schemes. Now, observe that there is a Cartesian diagram
\[
\begin{tikzcd}
	S\arrow[r]\arrow[d]&P_i\arrow[d]\\
	P_i\arrow[r]&P_i\times _SP_i
\end{tikzcd},
\] 
since we have the isomorphism
\begin{align*}
	P_i\times _SP_i\times _{P_i}S\simeq P_i.
\end{align*} 
Thus, we see that the map $S\to P_i$ is a base change of the relative diagonal of $P_i\to S$. Since $P_i\to S$ is proper, it is in particular separated; hence, $P_i\to P_i\times _SP_i$ is a closed immersion, so $S\to P_i$ is a closed immersion as well. In particular, the composition
\begin{align*}
	S\times _PS\xrightarrow{\simeq}S\times _{P_i}S\to S
\end{align*} 
is a closed immersion; hence, proper as desired.
\end{proof}
We now relate the Gaitsgory-Rozenblyum 6-functor formalism of ind-coherent sheaves to \zcref{def: scholze 6FF}. Let $\oblv_k:\bDGCat_k\to \bCat_{\infty}$ be the $(\infty,2)$-functor which forgets the $k$-linear structure. Composing $\mathsf{IndCoh}$ with this functor we obtain an $(\infty,2)$-functor
\begin{align*}
	\oblv_k\circ\mathsf{IndCoh}: \bSpan(\mathsf{PreStk}_{\mathrm{laft}},E)_{I,P}\to \bCat_{\infty}.
\end{align*}
\begin{cor}
	The underlying $(\infty,1)$-functor of $\oblv_k\circ\mathsf{IndCoh}$ is lax symmetric monoidal, where the $\Cat_{\infty}$ is equipped with the Cartesian monoidal structure.
\end{cor} 
\begin{proof}
	The $(\infty,2)$-functor $\mathsf{IndCoh}$ is lax symmetric monoidal by \zcref{prop: indcoh is a 6-functor formalism}; in particular, its underlying $(\infty,1)$-functor is also lax symmetric monoidal. We claim that the underlying $(\infty,1)$-functor of $\oblv_k$ is canonically lax symmetric monoidal, which follows from \cite[Propositions 4.8.1.4 and 4.8.1.15]{HA} (see also Variant 4.8.1.11 of \textit{op. cit.}).
\end{proof} 
\subsubsection{Ind-coherent sheaves on loop spaces}
Let $\calL-\coloneq \Map(S^1,-)$ be the free loop space functor, and denote by
\begin{align*}
	\mathsf{IndCoh}(\calL-): (\mathsf{PreStk}_{\laft})^{\op}\to& \bFun(BS^1,\bDGCat_k)\\
	X\mapsto &\mathsf{IndCoh}(\calL X),
\end{align*} 
where the $S^1$-action on $\mathsf{IndCoh}(\calL X)$ is induced by the loop rotation $S^1$-action on $\calL X$. Denote by $\oblv $ the forgetful functor $\bFun(BS^1,\bDGCat_k)\to \bDGCat_k$.

\begin{lem}\label{lem:forgetful-indcoh-loops-is-a-6FF}
	The composition $F\coloneq \oblv\circ\mathsf{IndCoh}(\calL-): (\mathsf{PreStk}_{\laft})^{\op} \to \bDGCat_k$ satisfies the conditions of \cite[Theorem 5.17]{CLL_six-functor-formalisms}; i.e., it is lax symmetric monoidal $(I,P)$-biadjointable with respect to the classes $I$ and $P$ given respectively by ind-inf-schematic ind-proper maps and equivalences.
\end{lem} 
\begin{proof}
	Note that $F$ is lax symmetric monoidal because $\calL=\Map(S^1,-)$ is strictly symmetric monoidal, as it preserves limits. To check that $F $ is lax symmetric monoidal $(I,P)$-biadjointable, it is sufficient to check the following.
\begin{enumerate}[label=\arabic*)]
	\item $\calL$ preserves pullbacks.
	\item $\calL$ preserves the classes $I$ and $P$.
\end{enumerate}	
Indeed, $\mathsf{IndCoh}$ defines a 6-functor formalism on $\mathsf{PreStk}_{\mathrm{laft}}$ by \zcref{prop: indcoh is a 6-functor formalism}; hence, its composition with $\calL$ will satisfy the conditions of the cited theorem if $\calL$ satisfies (1) and (2). Since $\calL$ is a limit, (1) is immediate. We will now prove (2). Let $f:X\to Y$ be an ind-inf-schematic ind-proper map of laft prestacks. There is a chain of isomorphisms
\begin{align*}
	X\times _{X\times X}(X\times _YX)\simeq X\times _{X\times X}(X\times X\times _{Y\times Y}Y)\simeq X\times _{Y\times Y}Y\simeq X\times _Y\calL Y.
\end{align*} 
Thus, $\calL f$ factors as the composition
\begin{align*}
	\calL X\xrightarrow{g} X\times _Y\calL Y\xrightarrow{p} \calL Y,
\end{align*} 
where $p$ is the projection and $g$ is the composition
\begin{align*}
	\calL X\xrightarrow{\simeq} X\times _{X\times X}X\to X\times _{X\times X}(X\times _YX)\xrightarrow{\simeq} X\times _Y\calL Y.
\end{align*} 
In the above, $X\to X\times _YX$ is the relative diagonal $\Delta_f$ of $f$. In particular, $p$ is a base change of $f$ and $g$ is a base change of $\Delta_f$. Thus, it is sufficient to show that ind-inf-schematic ind-proper maps are preserved by base change, composition, and relative diagonals. Preservation under base change and composition is by construction, and preservation under relative diagonals is \zcref{lem: relative diagonals of maps in I are in I}.
\end{proof}

\begin{cor}\label{cor:indcoh-loops-is-a-6FF}
	The $(\infty,2)$-functor $\mathsf{IndCoh}(\calL-):(\mathsf{PreStk}_{\mathrm{laft}})^{\mathrm{op}}\to \bFun(BS^1,\bDGCat_k)$ satisfies the conditions of \cite[Theorem 5.14]{CLL_six-functor-formalisms}; i.e., it is lax symmetric monoidal $(I,P)$-biadjointable for the classes of maps $I$ and $P$ respectively given by ind-inf-schematic ind-proper maps and equivalences.
\end{cor} 
\begin{proof}
	It is lax symmetric monoidal, as the functor $\oblv:S^1\bmod\to \bDGCat_k$ is strictly symmetric monoidal and the composition $F\coloneq \oblv\circ \mathsf{IndCoh}(\calL-)$ is lax symmetric monoidal by \zcref{lem:forgetful-indcoh-loops-is-a-6FF}. This lax symmetric monoidal structure is moreover $S^1$-equivariant: If $X,Y\in \mathsf{PreStk}_{\text{laft}}$, then the equivalence
	\begin{align*}
		\mathcal{L}(X\times Y)\simeq \mathcal{L}X\times \mathcal{L}Y
	\end{align*} 
	is canonically $S^1$-equivariant for the diagonal loop rotation action. Applying the Gaitsgory-Rozenblyum 6-functor formalism, we thus obtain an exterior product functor
	\begin{align*}
		\mathsf{IndCoh}(\mathcal{L}X)\otimes \mathsf{IndCoh}(\mathcal{L}Y)\to \mathsf{IndCoh}(\mathcal{L}X\times \mathcal{L}Y)\simeq \mathsf{IndCoh}(\mathcal{L}(X\times Y))
	\end{align*} 
	which is $S^1$-equivariant for the diagonal $S^1$-action, and this equivariance extends to arbitrary finite products, giving the desired lax symmetric monoidal structure.

	We now check condition (2) of Definition 5.13 of \emph{op. cit.}. To do this, we need to check that the exterior product functor admits an $S^1$-equivariant lift, in which case the condition follows from strict monoidality of $\oblv$ and \zcref{lem:forgetful-indcoh-loops-is-a-6FF}. The relevant mates in condition (2a) of \textit{loc. cit.} are $S^1$-equivariant by the lax symmetric monoidality constructed above, and their images in $\DGCat_k$ under the forgetful functor are isomorphisms by \zcref{lem:forgetful-indcoh-loops-is-a-6FF}. Since this forgetful functor detects equivalences, we conclude that \cite[Definition 5.13(2)]{CLL_six-functor-formalisms} is satisfied.

	We would now like to check that $\mathsf{IndCoh}(\calL -)$ satisfies condition (1) of \cite[Definition 5.13]{CLL_six-functor-formalisms}. Thus, suppose given a Cartesian diagram in $\mathsf{PreStk}_{\mathrm{laft}}$
	\[
	\begin{tikzcd}
		X'\arrow[r]\arrow[d]&X\arrow[d]\\
		Y'\arrow[r]&Y
	\end{tikzcd}
	\] 
	for which we would like to check that $\mathsf{IndCoh}(\calL-)$ is $(I,P)$-biadjointable. Observe that, by definition, we are merely checking if certain natural transformations in $S^1\bmod$ are isomorphisms. Since an $S^1$-equivariant map is an isomorphism if and only if its underlying non-equivariant map is, we see that $\mathsf{IndCoh}(\calL-)$ is $(I,P)$-biadjointable if and only if $F$ is. 
\end{proof}
 
\section{Potent $D$-modules}

Following \cite{BZPotent}, we make the following defintion.

\begin{defn}
	Let $X$ be a derived stack such that $\mathsf{IndCoh}(\mathcal{L} X)\in\bFun(BS^1,\bCoh_t)$. The category of \emph{potent $D$-modules}
	on $X$ is defined to be
	\[\Pot(X) := \IndCoh(\calL X)^{\tT} \, (= \Ind(\Coh(\calL X)^{tS^1})).\]
\end{defn}

\begin{rem}
	The name in the above definition takes inspiration from the equivalence
	under Koszul duality between cyclic sheaves on the unipotent loop
	space of $X$ and ordinary $D$-modules on $X$. Potent $D$-modules
	correspond to cyclic sheaves on the full---or \emph{potent}---loop space rather
	than just the \emph{uni}potent loop space of $X$.
\end{rem}

\subsection{$\tT$ functoriality}
\label{ssec: tTate functoriality}

\begin{comment}
	
Suppose $F: \scrC \to \scrD \in \bCoh_t$. Then there exists a functor $\scrR(F): \scrR(\scrC) \to \scrR(\scrD)$
defined as follows. There are canonical functors, $j_{\scrC}: \scrC_{<\infty} \to \scrR(\scrC)$
and $j_{\scrD}: \scrD_{<\infty} \to \scrR(\scrD)$, obtained as the restriction to eventually connective
subcategories of the right adjoints to the natural functors, $\scrR(\scrC) \to \scrC$ and $\scrR(\scrD) \to \scrD$, respectively.
Let $F':\Coh(\scrC) \to \scrR(\scrD)$ denote the functor obtained by composition from
the following sequence of functors,
\[\Coh(\scrC) \xhookrightarrow{i_{\scrC}} \scrC_{<\infty} \xrightarrow{F_{<\infty}} \scrD_{<\infty} \xrightarrow{j_{\scrD}} \scrR(\scrD).\]
Moreover, the functor $F'$ is exact.

\end{comment} 

\begin{lem}\label{lem: regularization is a 2-functor}
    Regularization defines an $(\infty,2)$-functor $\scrR(-): \bCoh_{t,A} \to \bReg_{t,A}$.
\end{lem}
\begin{proof}
	Let $\mathscr{D}\in \bCoh_{t,A}$, and denote by $\varepsilon_{\mathscr{D}}:\mathscr{R}(\mathscr{D})\to \mathscr{D}$ the canonical functor of \cite[Proposition 4.3.2]{Preygel}, which is shown to be $A$-linear in \zcref{lem: multilinear regularization functor}. It is sufficient to check the following universal property: For any $\mathscr{C}\in \bReg_{t,A}$, postcomposition with $\varepsilon_{\mathscr{D}}$ induces an equivalence of $\infty$-categories
	\begin{align*}
		\left( \varepsilon_{\mathscr{D}} \right) _*:\Fun_A^{\text{coh}}(\mathscr{C},\mathscr{R}(\mathscr{D}))\xrightarrow{\simeq}\Fun_A^{\text{coh}}(\mathscr{C},\mathscr{D}),
	\end{align*} 
	which is moreover natural in $\mathscr{C}\in \bReg_{t,A}$. Indeed, once we have this universal property naturally in $\mathscr{C}\in \bReg_{t,A}$, this will show that $\mathscr{R}$ is an $(\infty,2)$-categorical right adjoint to the inclusion $\bReg_{t,A}\hookrightarrow \bCoh_{t,A}$. Now, for any coherent functor $F: \mathscr{C}\to \mathscr{D}$, \cite[Lemma 4.6.1]{Preygel} supplies an equivalence
	\begin{align*}
		\varepsilon_{\mathscr{D}} \mathscr{R}(F)\simeq F \varepsilon_{\mathscr{C}}.
	\end{align*} 
	Thus, the functor $(\varepsilon_{\mathscr{D}})_*$ is essentially surjective, since $\varepsilon_{\mathscr{C}}$ is an equivalence for $\mathscr{C}\in\bReg_{t,A}$. For fully faithfulness, we can again use \cite[Lemma 4.6.1]{Preygel} and the universal property of $\mathsf{Ind}(-)$ to obtain the natural equivalence
	\begin{align*}
		\Fun_A^{\text{coh}}(\mathscr{C},\mathscr{R}(\mathscr{D}))\simeq&\Fun_A^{\text{ex,lb}}(\mathsf{Coh}(\mathscr{C}),\mathscr{R}(\mathscr{D})),
	\end{align*} 
	where $\Fun_A^{\text{ex,lb}}$ denotes exact and left $t$-exact up to a shift functors. Similarly, we have the natural equivalence
	\begin{align*}
		\Fun_A^{\text{coh}}(\mathscr{C},\mathscr{D})\simeq\Fun_A^{\text{ex,lb}}(\mathsf{Coh}(\mathscr{C}),\mathscr{D}),
	\end{align*} 
	again by the universal property of $\mathsf{Ind}(-)$. Now, the category on the right-hand-side admits only functors which are in particular left $t$-bounded; hence, it is in fact equivalent to $\Fun_A^{\text{ex,lb}}(\mathsf{Coh}(\mathscr{C}),\mathscr{D}^+)$, since $\mathsf{Coh}(\mathscr{C})$ contains only $t$-bounded objects of $\mathscr{C}$. Thus, by the above, we have the chain of equivalences
	\begin{align*}
		\Fun_A^{\text{coh}}(\mathscr{C},\mathscr{R}(\mathscr{D}))\simeq&\Fun_A^{\text{ex,lb}}(\mathsf{Coh}(\mathscr{C}),\mathscr{R}(\mathscr{D})^+)\\
		\xrightarrow{\simeq}&\Fun_A^{\text{ex,lb}}(\mathsf{Coh}(\mathscr{C},\mathscr{D}^+)\\
		\simeq&\Fun_A^{\text{coh}}(\mathscr{C},\mathscr{D}),
	\end{align*} 
	where we use that $\varepsilon_{\mathscr{D}}$ is an equivalence when restricted to the bounded below subcategories in the second line. Since postcomposition with $\varepsilon_{\mathscr{D}}$ is naturally compatible with precomposition with $\mathscr{C}'\to \mathscr{C}$, these equivalences exhibit $\mathscr{R}$ as the $(\infty,2)$-categorical right adjoint to the inclusion $\bReg_{t,A}\hookrightarrow \bCoh_{t,A}$.
	\end{proof}
	 
\begin{lem}\label{lem: S^1-invariants is a 2-functor}
	Taking $S^1$-invariants defines an $(\infty,2)$-functor
	\begin{align*}
		(-)^{S^1}:\bFun(BS^1,\bDGCat_k)\to \bDGCat_{k[u]}.
	\end{align*} 
\end{lem} 
\begin{proof}
Let $\mathscr{C}\in\bFun(BS^1,\bDGCat_k)$. Equipping $\mathsf{Vect}$ with the trivial $S^1$-action, $\mathsf{Vect}^{S^1}$ acts on $\mathscr{C}^{S^1}$. Since $\mathsf{Perf}(k[u])\simeq\mathsf{Perf}(k)^{S^1}\subset \mathsf{Vect}^{S^1}$, we can ind-extend this inclusion to obtain a symmetric monoidal functor
	\begin{align*}
		k[u]\mod\to \mathsf{Vect}^{S^1}.
	\end{align*} 
	Thus, we can restrict the $\mathsf{Vect}^{S^1}$-action along this functor to obtain the desired $k[u]$-linear structure on $\mathscr{C}^{S^1}$. Now, evaluation at $k\in\mathsf{Vect}$ identifies $\Hom_{\bFun(BS^1,\bDGCat_k)}(\mathsf{Vect},\mathscr{C})$ with $\mathscr{C}^{S^1}$. A 1-morphism $\mathscr{C}\to\mathscr{D}$ in $\bFun(BS^1,\bDGCat_k)$ acts on this mapping category by postcomposition, and a 2-morphism $F\Rightarrow G$ acts by whiskering. These actions are functorial; thus, taking $S^1$-invariants is $(\infty,2)$-functorial in $\mathsf{Vect}^{S^1}$-module categories, and restriction of scalars as above gives the desired functoriality in $\bDGCat_{k[u]}$.
\end{proof}

\begin{cor}
\label{cor: tTate is a 2-functor}
    Regularized Tate invariants defines an $(\infty,2)$-functor,
    \[(-)^{\tT}: \bFun(BS^1, \bCoh_t) \to \bDGCat_{\ku}.\]
    \end{cor}

\begin{proof}
	To show the first part of the proposition, it is sufficient by \zcref{lem: regularization is a 2-functor} to demonstrate that each
	of the operations $(-)^{S^1}: \bFun(BS^1, \bCoh_t) \to \bCoh_{t,k[u]}$ and $- \otimes_{k[u]} \ku$ is an $(\infty, 2)$-functor. For the former, we will show that the $(\infty,2)$-functor of \zcref{lem: S^1-invariants is a 2-functor} restricts to an $(\infty,2)$-functor
	\begin{align*}
		\bFun(BS^1,\bCoh_t)\to \bCoh_{t,k[u]}.
	\end{align*} 
	To see this, note that $S^1$-invariants of categories with coherent $t$-structure still have coherent $t$-structures by \cite[Lemma 4.5.3]{Preygel}. Now, given a 1-morphism $F: \mathscr{C}\to \mathscr{D}$ in $\bFun(BS^1,\bCoh_t)$, $F^{S^1}$ is continuous and $k[u]$-linear by \zcref{lem: S^1-invariants is a 2-functor}. Since the forgetful functor $\mathscr{C}^{S^1}\to \mathscr{C}$ is $t$-exact (see \cite[Section 3.1.3]{ES-Kashiwara}), we see that $F^{S^1}$ is left $t$-exact up to a shift if $F$ is. Finally, since the mapping categories of $\bCoh_{t,k[u]}$ are full subcategories of the mapping categories of $\bDGCat_{k[u]}$, we obtain the desired restriction on 2-morphisms. Thus, we can form the composite
	\begin{align*}
		\bFun(BS^1,\bCoh_t)\xrightarrow{(-)^{S^1}}\bCoh_{t,k[u]}\xrightarrow{\mathscr{R}_{k[u]}}\bReg_{t,k[u]}\xrightarrow{\oblv_t}\bDGCat_{k[u]}\xrightarrow{-\otimes _{k[u]}\ku}\bDGCat_{\ku},
	\end{align*} 
	where $\oblv_t$ forgets the $t$-structure, $\mathscr{R}_{k[u]}$ is the $k[u]$-linear version of the regularization $(\infty,2)$-functor constructed in \zcref{lem: regularization is a 2-functor}, and $-\otimes _{k[u]}\ku$ denotes the $(\infty,2)$-functor of extension of scalars. This composite computes $\tT$ invariants defined in \zcref{sssec:Tate invariants} by \cite[Proposition 4.5.5]{Preygel}.
\end{proof}
\begin{cor}\label{cor: tTate preserves adjunctions}
Suppose that $L: \scrC \leftrightarrows \scrD: R$ is an adjunction in $\bFun(BS^1,\bCoh_t)$. Then the induced functors 
	$L^{\tT}: \scrC^{\tT} \leftrightarrows \scrD^{\tT}: R^{\tT}$
	are adjoint to one another.
\end{cor}
\begin{proof}
	This follows immediately by applying the $(\infty,2)$-functor of $\tT$ invariants of \zcref{cor: tTate is a 2-functor} to the unit and counit of the adjunction $L\dashv R$, since $L,R$ are 1-morphisms in $\Fun(BS^1,\bCoh_t)$ by assumption.
\end{proof}

\subsection{Ind-finite type stacks}
\label{ssec: ind-finite type stacks}
One obstruction to constructing a six-functor formalism for potent $\mathcal{D}$-modules is that the domain of definition of $\mathsf{IndCoh}(\calL-)$ is correspondences of all laft prestacks. However, the functor $(-)^{\tT}$ can only be applied to objects of $\bCoh_t$. In particular, all $!$-pullbacks must be coherent functors for the composition to be well-defined. Thus, we need to restrict our domain of definition to a suitable subcategory which will include all stacks of interest in \cite{ES-Kashiwara}. 

Recall the category $\mathsf{indGStk}$ of \cite{CautisWilliams}, which is the category of ind-geometric stacks. These are convergent stacks $X$ which admit an expression $X\simeq \colim X_\alpha$ as a filtered colimit of truncated geometric stacks along closed immersions.\footnote{See Definition 3.1 of \textit{op. cit.} for the definition of geometric stack. Observe also that because the affine atlas is faithfully flat, rather than smooth, such stacks are not necessarily Artin.} Following \textit{op. cit.}, we call such a stack \textit{reasonable} if it admits a presentation with finitely presented structure maps and \textit{coherent} if it admits a reasonable presentation whose terms are coherent geometric stacks (see Definition 3.3 of \textit{op. cit.}). Finally, we call such a stack \textit{ind-tamely presented} if it is reasonable and its structure map to $\Spec k$ satisfies the conditions of \cite[Proposition 5.5]{CautisWilliams2}. Following \textit{loc. cit.}, we denote the category of ind-tamely presented $\mathsf{indGStk}^{\mathrm{tm}}$. The following example shows that filtered colimits of Artin stacks almost of finite type and with affine diagonal are ind-tamely presented; in particular, formal completions of such stacks are ind-tamely presented.
\begin{exmp}\label{ex: aft artin stacks are tamely presented}
	Every Artin stack almost of finite type and with affine diagonal is tamely presented over $\Spec k$. Indeed, such a stack $X$ is necessarily 1-Artin by the affine diagonal hypothesis; hence, it admits a smooth atlas $U\to X$ with $U$ affine and almost of finite type. Since the diagonal of $X$ is affine, so is $U\to X$. Thus, $U\to X$ is an affine and almost finitely presented flat cover. Now, given an affine scheme $S$, $U\times S\to S$ is strictly tamely presented, because the composite $U\to X\to \pt$ is affine almost of finite presentation, hence strictly tamely presented (cf. \cite[Definitions 3.2 and 4.6]{CautisWilliams2}). This suffices to conclude that $X$ is tamely presented by \cite[Definition 4.6]{CautisWilliams2}, since then $U\times S\to X\times S$ is a strictly tamely presented flat cover.
\end{exmp} 
 
\begin{lem}\label{lem: ind-finite type stacks have presentations by algebraic stacks}
	Let $X$ be a laft ind-geometric stack. Then $X$ is ind-tamely presented if and only if $X$ admits a reasonable presentation by truncated algebraic stacks almost of finite type. Moreover, when these conditions hold, every stage of every reasonable presentation of $X$ is an algebraic stack almost of finite type. 
\end{lem} 
\begin{proof}
	Suppose first that $X$ is ind-tamely presented, and let $X\simeq\colim_\alpha X_\alpha$ be a reasonable presentation of $X$ as an ind-geometric stack. By \cite[Proposition 5.5]{CautisWilliams2}, each $X_\alpha$ is tamely presented. Since the inclusion $X_\alpha\to X$ is almost finitely presented, $X_\alpha$ is laft. Moreover, each $X_\alpha$ is truncated by definition of an ind-geometric presentation. We will now show that it has a smooth atlas. Choose a strictly tamely presented flat cover $\Spec B\to X_\alpha$ with $B$ strictly tamely presented over $k$. If $ X_\alpha$ is $n$-truncated, then so is $B$; thus, there is a presentation $B=\colim_i B_i$ by finitely $n$-presented $k$-algebras $B_i$ with $B$ flat over each $B_i$ and each $B_i$ almost finitely presented over $k$, since $k$ is Noetherian. Since $X_\alpha$ is laft and $ n$-truncated and the $B_i$ are $n$-truncated as well, we have
	\begin{align*}
		X_\alpha(B)\simeq\colim_i X_\alpha(B_i).
	\end{align*} 
	In particular, the point $\Spec B\to X_\alpha$ in $X_\alpha(B)$ came from some finite stage; hence, there is a factorization
	\begin{align*}
		\Spec B\xrightarrow{f_i} \Spec B_i\xrightarrow{g_i} X_\alpha
	\end{align*} 
	for some $i$, where $f_i$ is flat and $g_i$ is affine (because $X_\alpha$ has affine diagonal) and almost finitely presented. Write the base change $\Spec C\coloneq \Spec B \times _{X_\alpha}\Spec B_i$, and let $p:\Spec C\to \Spec B$ be the projection. Note that $C$ is $ n$-truncated by flatness over $B_i$. Since $p$ is a base change of $g_i$, $C\otimes_B H^0(B)$ is almost finitely presented over $H^0(B)$. Moreover, $M \coloneq \tau^{\ge -n-1}(C\otimes _B H^0(B))$ is finitely $(n+1)$-presented as a $(C\otimes_B H^0(B))$-module, since it is the truncation of the free rank 1 $(C\otimes_B H^0(B))$-module. Thus, \cite[Proposition 6.1.4.5]{SAG} applied to $M$ as a module over $H^0(B)$ with tor amplitude bound 0 implies that the locus where $M$ is flat over $H^0(B)$ is open in $|\Spec (C\otimes_B H^0(B))|=|\Spec C|$. We claim that this is precisely the flat locus of $p$. Indeed, if $C_x$ is flat over $B$, then $(C\otimes _B H^0(B))_x$ is discrete and flat over $H^0(B)$, so $M_x\simeq (C\otimes _BH^0(B))_x$ is flat. Conversely, if $M_x$ is flat, then \cite[Propositions 8.5]{PortaYu} implies that $C_x$ is flat to order $n+1$ over $B$. Thus, \cite[Propositions 8.6]{PortaYu} gives that $C_x$ is flat over $B$, since both $B$ and $C_x$ are $n$-truncated. We conclude that the flat locus of $p$ is open in $\Spec C$. Now, let $F_i\subset |\Spec B_i|$ be the flat locus of $g_i$. Its inverse image under the projection $\Spec C\to \Spec B_i$ is the flat locus of $p$, since pointwise flatness is preserved and detected by flat base change. Moreover, this projection is fpqc; hence, \cite[Tag 02JY]{Stacks} applied to its classical truncation shows that $F _i$ is open. Give $F_i$ its induced derived open subscheme structure to obtain an open subscheme of $\Spec B_i$. We claim that the image of $f_i$ lies in $F_i$. Indeed, for each $x\in\Spec B$, the flat map $f_i$ induces a faithfully flat map on local rings at $x$ and $f_i(x)$. Thus, flatness of $g_if_i$ at $x$ implies flatness of $g_i$ at $f_i(x)$, since this can be checked flat locally on $X_\alpha$. Since $g_if_i$ is surjective, $F_i\to X_\alpha$ is surjective as well. It is flat and locally almost of finite presentation, so an affine cover of $F_i$ gives an fppf atlas of $X_\alpha$. By \cite[Th\'eor\`eme 2.1]{Toen-flat-is-smooth}, $X_\alpha$ thus admits a smooth atlas. Since it is quasi-compact and laft, it is an algebraic stack almost of finite type.

	Now, suppose that $X$ admits a reasonable presentation by truncated algebraic stacks almost of finite type. These stacks necessarily have affine diagonal by definition of geometric stacks; hence, they are tamely presented by \zcref{ex: aft artin stacks are tamely presented}. \cite[Proposition 5.5]{CautisWilliams2} then implies that $X$ is ind-tamely presented as desired.
\end{proof} 
This lemma motivates the following definition.
\begin{defn}\label{def: ind-ft stacks}
	Let
\begin{align*}
	\mathsf{indGStk}_{\text{ft}}\subset \mathsf{indGStk}^{\text{tm}}
\end{align*} 
be the full subcategory consisting of laft stacks $X$ for which $X_{\text{red}}$ is a finite type algebraic stack.
\end{defn}  
\begin{rem}
	The condition on the reduction in \zcref{def: ind-ft stacks} is used in \zcref{lem: bounded below iff bounded below on reduced stack} below, which is needed to prove that potent $\mathcal{D}$-modules is a well-defined $(\infty,2)$-functor on this subcategory.
\end{rem} 

\begin{prop}\label{prop: ind-ft stacks are coherent}
	Let $X\in \mathsf{indGStk}_{\mathrm{ft}}$. Then $X$ is a coherent ind-geometric stack.
\end{prop} 
\begin{proof}
	Since all stacks in the reasonable presentation $X\simeq X_\alpha$ are almost of finite type, it suffices to check that the $X_\alpha$ are coherent geometric stacks. Indeed, the smooth affine atlas $U\to X_\alpha$ can be chosen to be almost of finite type, hence Noetherian; this suffices by \cite[Definition 3.3]{CautisWilliams}.
\end{proof} 
\begin{prop}\label{prop: ind-ft stacks are closed under limits}
	The category $\mathsf{indGStk}_{\mathrm{ft}}$ is closed under finite limits, 
\end{prop} 
\begin{proof}
	We follow the strategy of \cite[Proposition 4.21]{CautisWilliams}. By \cite[Proposition 4.4.2.6]{HTT}, it suffices to check closure under fiber products. By \cite[Proposition 4.21]{CautisWilliams}, a fiber product in $\mathsf{indGStk}_{\text{ft}}$ is ind-geometric, and it is laft by finite limit closure of $\mathsf{PreStk}_{\mathrm{laft}}$. Hence, it suffices to check that such a fiber product $X'\coloneq X\times _YY'$ is ind-tamely presented and has $X'_{\text{red}}$ an algebraic stack of finite type. Suppose first that $X$ and $Y'$ are truncated algebraic stacks almost of finite type, and let $Y\simeq\colim Y_\alpha$ be a reasonable presentation by truncated algebraic stacks almost of finite type, which exists by \zcref{lem: ind-finite type stacks have presentations by algebraic stacks}. By \cite[Proposition 4.7]{CautisWilliams} and \cite[Example 7.3.4.7]{HTT}, we have $X'\simeq\colim_{\beta \ge \alpha}X'_\beta$, where $X'_\beta\coloneq X\times _{Y_\beta} Y'$, for some $\alpha$. Each $X'_\beta$ is an algebraic stack almost of finite type and the transition maps $X'_\beta\to X'_{\beta'}$ are closed immersions by \cite[Proposition 3.23]{CautisWilliams}. Thus, we can write $X'\simeq\colim_{\beta\ge \alpha,n\ge 0}\tau_{\le n}X'_\beta$, and each term is truncated and almost of finite type, with the transition maps almost finitely presented closed immersions; hence, $X'\in \mathsf{indGStk}_{\mathrm{ft}}$ by \zcref{lem: ind-finite type stacks have presentations by algebraic stacks}.

	Now, suppose given reasonable presentations $X\simeq \colim X_\alpha$ and $Y'\simeq\colim Y'_\beta$, which are necessarily by truncated algebraic stacks almost of finite type by \zcref{lem: ind-finite type stacks have presentations by algebraic stacks}. Then $X'\simeq\colim_{\alpha,\beta}X'_{\alpha\beta}$, where $X'_{\alpha\beta}\coloneq X_\alpha \times _Y Y'_{\beta}$. Each $X'_{\alpha\beta}$ admits a reasonable presentation by algebraic stacks almost of finite type by the previous argument. The transition maps $X'_{\alpha\beta}\to X'_{\alpha'\beta'}$ are composites of base changes of the original almost finitely presented closed immersions; hence, they are almost ind-finitely presented ind-closed immersions. By \zcref{lem: ind-finite type stacks have presentations by algebraic stacks}, each $X'_{\alpha\beta}$ is ind-tamely presented. Since the transition maps are almost ind-finitely presented ind-closed immersions, \cite[Proposition 5.10]{CautisWilliams2} implies that $X'\simeq\colim_{\alpha,\beta}X'_{\alpha\beta}$ is ind-tamely presented. 

	It remains to check that $X'_{\text{red}}$ is an algebraic stack of finite type. To see this, we claim that there is an isomorphism
	\begin{align*}
		X'_{\text{red}}\simeq\left( X_{\text{red}}\times _{Y_{\text{red}}}^{\text{cl}}Y'_{\text{red}} \right) _{\text{red}},
	\end{align*} 
	where $\times ^{\text{cl}}$ denotes the fiber product in classical stacks. Since $(-)_{\text{red}}$ is right adjoint to the inclusion $\mathsf{Stk}_{\text{cl,red}}\hookrightarrow \mathsf{Stk}$, it takes limits in the target to limits in the source. We then conclude by observing that the fiber product in reduced classical stacks is the reduction of the fiber product in classical stacks.
\end{proof}
This category is still quite interesting, as it includes all examples of interest in \cite{ES-Kashiwara}, as well as their formal completions.
\begin{prop}\label{prop: formal completions are ind-ft}
    Let $Z\hookrightarrow X$ be an almost finitely presented inclusion of a closed substack of an algebraic stack $X$ almost of finite type and with affine diagonal, and let $\wih{X_Z}$ be the formal completion. Then $Z,X,\wih{X_Z}\in \mathsf{indGStk}_{\mathrm{ft}}$.
\end{prop}
\begin{proof}
	By \zcref{ex: aft artin stacks are tamely presented}, $X\in \mathsf{indGStk}_{\text{ft}}$. By pulling back the smooth affine atlas of $X$ to $Z$, we find that $Z$ is also equipped with an smooth affine atlas. Since $Z\hookrightarrow X$ is almost finitely presented, $Z$ is almost of finite type; hence, it is also an object of $\mathsf{indGStk}_{\text{ft}}$. Finally, $(\wih{X_Z})_{\text{red}}\simeq Z_{\text{red}}$, and the arguments of \cite[Lemmas 3.34 and 3.35]{ES-Kashiwara} imply that it has a reasonable presentation by algebraic stacks almost of finite type, so $\wih{X_Z}\in \mathsf{indGStk}_{\text{ft}}$ by \zcref{lem: ind-finite type stacks have presentations by algebraic stacks}.
\end{proof}
\begin{cor}
	In the setting of \zcref{prop: formal completions are ind-ft}, $\mathcal{L}Z,\mathcal{L}X,\mathcal{L}(\wih{X_Z})\in \mathsf{indGStk}_{\mathrm{ft}}$.
\end{cor} 
\begin{proof}
    This follows immediately from \zcref{prop: ind-ft stacks are closed under limits}, \zcref{prop: formal completions are ind-ft}, and \cite[Lemma 3.57]{ES-Kashiwara}.
\end{proof}
\begin{cor}\label{cor: ind-finite type stacks have IndCoh=Ind(Coh)}
If $X\in \mathsf{indGStk}_{\mathrm{ft}}$, then $\mathsf{Ind}(\mathsf{Coh}(X))\simeq \mathsf{IndCoh}(X)$, where the latter category is defined by descent.
\end{cor} 
\begin{proof}
	This follows from \cite[Remark 5.20 and Proposition 5.30]{CautisWilliams} and \cite[Lemma 3.35]{ES-Kashiwara}.
\end{proof} 
\begin{lem}\label{lem: CW pushforward = GR pushforward}
	Let $f:X\to Y$ be an ind-inf-schematic ind-proper map in $\mathsf{indGStk}_{\mathrm{ft}}$. Under the canonical equivalences of \zcref{cor: ind-finite type stacks have IndCoh=Ind(Coh)}, there is a natural equivalence
	\begin{align*}
		f_*^{\mathrm{GR}}\simeq f_*^{\mathrm{CW}},
	\end{align*} 
	where $f_*^{\mathrm{GR}}$ denotes the Gaitsgory-Rozenblyum $* $-pushforward (\cite[Chapter 3 Theorem 5.4.3]{GR17II}) and $f_*^{\mathrm{CW}}$ denotes the Cautis-Williams $*$-pushforward (see \cite[Definition 5.16]{CautisWilliams} and the paragraph following it). Moreover, the left adjoint $f_*$ to $f^!$ is left $t$-exact.
\end{lem} 
\begin{proof}
	Choose reasonable presentations $X\simeq\colim_\alpha X_\alpha$ and $Y\simeq\colim_\beta Y_\beta$ with inclusions $i_\alpha:X_\alpha\to X$ and $j_\beta:Y_\beta\to Y$. By \cite[(5.18)]{CautisWilliams}, we have colimit presentations along the CW pushforwards $i_{\alpha*}^{\mathrm{CW}}$ and $j_{\beta*}^{CW}$,
	\begin{align*}
		\mathsf{IndCoh}(X)\simeq&\colim_\beta \mathsf{IndCoh}(X_\alpha)\\
		\mathsf{IndCoh}(Y)\simeq&\colim_\beta \mathsf{IndCoh}(Y_\beta).\footnotemark
	\end{align*} 
	\footnotetext{We do not need to write $\mathsf{IndCoh}^{\mathrm{GR}}$ and $\mathsf{IndCoh}^{\mathrm{CW}}$, because for $Z\in \mathsf{indGStk}_{\mathrm{ft}}$, these categories are equivalent by \cite[Remark 5.20]{CautisWilliams} and \zcref{prop: ind-ft stacks are coherent}.}
	Now, for every truncated affine scheme $S$ almost of finite type, \cite[Proposition 4.7]{CautisWilliams} gives an equivalence
	\begin{align*}
		X(S)\simeq&\colim_\alpha X_\alpha(S)\\
		Y(S)\simeq&\colim_\beta Y_\beta(S).
	\end{align*} 
	Since such test affines are precisely those which define $\mathsf{PreStk}_{\mathrm{laft}}$, the presentations for $X$ and $Y$ are colimits in that category. Thus, \cite[Lemma 10.1.3]{GaitsgoryIndCoh} gives
	\begin{align*}
		\mathsf{IndCoh}(X)\simeq&\lim_{\alpha}\mathsf{IndCoh}(X_\alpha)\\
		\mathsf{IndCoh}(Y)\simeq&\lim_\beta \mathsf{IndCoh}(Y_\beta),
	\end{align*} 
	where the transition functors are respectively $i_{\alpha\beta}^!$ and $j_{\alpha\beta}^!$ for $i_{\alpha\beta}:X_{\alpha}\to X_\beta$ and similarly for $j_{\alpha\beta}$. Passing to left adjoints, we obtain presentations
	\begin{align*}
		\mathsf{IndCoh}(X)\simeq&\colim_\alpha\!{}^{\mathrm{GR}} \mathsf{IndCoh}(X_\alpha)\\
		\mathsf{IndCoh}(Y)\simeq&\colim_\beta\!{}^{\mathrm{GR}} \mathsf{IndCoh}(Y_\beta),
	\end{align*} 
	where we use the notation $\colim^{\mathrm{GR}}$ to denote the fact that the transition functors here are $i_{\alpha\beta*}^{\mathrm{GR}}$ and $j_{\alpha\beta*}^{\mathrm{GR}}$, respectively. Now, fix $\alpha$. By \cite[Proposition 4.7]{CautisWilliams}, choose $f_{\alpha\beta}:X_\alpha\to Y_\beta$ with $j_\beta f_{\alpha\beta}\simeq fi_\alpha$. We claim that $f_{\alpha\beta}$ is proper and schematic. Indeed, given an almost finite type affine scheme $S$ and an $S$-point $S\to Y_\beta$, $(X_\alpha \times _{Y_\beta}S)_{\mathrm{red}}$ is a closed substack of $(X\times _YS)_{\mathrm{red}}$. By \zcref{lem: ind-inf-schematic maps are schematic}, $X\times _YS$ is an ind-scheme; moreover, by the argument of \textit{loc. cit.}, its reduction is a scheme. Thus, $(X_\alpha \times _{Y_\beta}S)_{\mathrm{red}}$ is a quasi-compact classical, reduced scheme, proper over $S_{\mathrm{red}}$ by the definition of ind-properness. Again applying the argument of \zcref{lem: ind-inf-schematic maps are schematic}, we find that $X_\alpha \times _{Y_\beta}S$ is a derived scheme, while properness is detected on the classical reduction. Thus, $f_{\alpha\beta}$ is schematic and proper. Now, for $\mathcal{E}\in \mathsf{Coh}(X_\alpha)$, both the Gaitsgory-Rozenblyum and Cautis-Williams $*$-pushforwards along $f_{\alpha\beta}$ are the usual derived direct image for coherent sheaves. Indeed, for the former this follows from \cite[(3.17)]{DrinfeldGaitsgory}, and for the latter this is the compatibility following \cite[(5.14)]{CautisWilliams}. Hence, functoriality of (Gaitsgory-Rozenblyum) $\mathsf{IndCoh}$ and the construction preceding \cite[Definition 6.4]{CautisWilliams} give natural equivalences
	\begin{align*}
		f_*^{\mathrm{GR}}(i_{\alpha*}\mathcal{E})\simeq j_{\beta*}(f_{\alpha\beta*}\mathcal{E})\simeq f_*^{\mathrm{CW}}(i_{\alpha*}\mathcal{E}).
	\end{align*} 
	Since these equivalences are compatible with passage to higher stages $X_\gamma$ by functoriality of direct image, we conclude that $f_*^{\mathrm{GR}}$ and $f_*^{\mathrm{CW}}$ agree on $\mathsf{Coh}(X)\simeq\lim_\alpha \mathsf{Coh}(X_\alpha)$. Continuity extends this equivalence to $\mathsf{IndCoh}(X)$. Left $t$-exactness follows from the paragraph following \cite[Definition 6.4]{CautisWilliams}.
\end{proof} 
\begin{rem}\label{rem: CW pushforward = GR pushforward}
	Using \zcref{lem: CW pushforward = GR pushforward}, if $f$ is an ind-inf-schematic ind-proper map in $\mathsf{indGStk}_{\mathrm{ft}}$, we will simply write $f_*$ for the left adjoint to $f^!$.
\end{rem}
\begin{lem}\label{lem: ind-inf-schematic maps are schematic}
	Let $f:X\to Y$ be an ind-inf-schematic map in $\mathsf{indGStk}_{\mathrm{ft}}$. Then $f$ is ind-schematic.
\end{lem} 
\begin{proof}
	We claim it suffices to test ind-schematicity on almost finite type affines. Indeed, let $S\to Y$ be an $S$-point for $S$-arbitrary, and fix $n$. Then $\tau^{\ge -n}S\to Y$ factors through $S_n$ for some $n$-truncated almost finite type affine $S_n$, since $Y$ is laft. If $X\times _Y S_n$ is an ind-scheme, then base change to $\tau^{\ge -n}S$ gives an ind-scheme on $n$-truncated tests. Since the prestack $X\times _YS$ is convergent, we conclude that such tests were sufficient.
    
    Now, let $S\to Y$ be an $S$-point of $Y$, where $S$ is an almost finite type affine scheme, and let $Z=X\times _YS$. By \zcref{prop: ind-ft stacks are closed under limits}, $Z\in \mathsf{indGStk}_{\text{ft}}$, so $Z_{\text{red}}$ is an algebraic stack of finite type. Since it is also an ind-inf-scheme, it is thus a reduced, classical ind-scheme. If we write it as a filtered colimit of reduced, classical schemes along closed immersions, $Z_{\text{red}}\simeq\colim_i W_i$, then a reduced, classical smooth affine atlas of $Z_{\text{red}}$ factors through some $W_i$. Thus, the monomorphism $W_i\to Z_{\text{red}}$ is an isomorphism, so $Z _{\text{red}}$ is in fact a scheme. Now, by \zcref{lem: ind-finite type stacks have presentations by algebraic stacks}, we can write $Z=\colim_\alpha Z_\alpha$ with $Z_\alpha$ truncated algebraic stacks almost of finite type and transition maps closed immersions. Each $(Z_\alpha)_{\text{red}}$ is a closed subscheme of $Z_{\text{red}}$; hence, $Z_\alpha^{\text{cl}}$ is a scheme by \cite[Tag 0BPW]{Stacks}. By \cite[Definition 8.5.7(iii)]{KhanLectures}, $Z_\alpha$ is Deligne-Mumford, and by \cite[Corollary 1.6.7.4]{SAG}, $Z_\alpha$ is a derived scheme. Let $Z'=\colim_\alpha Z_\alpha$, where the colimit is taken in prestacks. By \cite[Proposition 4.7]{CautisWilliams}, $Z'\to Z$ is an equivalence on truncated affine tests. Since $Z$ is convergent, $Z$ is equivalent to the image of $Z'$ in convergent prestacks (under the right adjoint to the inclusion $\mathsf{PreStk}_{\mathrm{conv}}\hookrightarrow \mathsf{PreStk}$). Thus, \cite[Chapter 2 Corollary 1.4.4]{GR17II} implies that $Z$ is an ind-scheme. 
\end{proof}

\subsection{The $(\infty,2)$-functor of potent $\mathcal{D}$-modules}
\label{ssec: 2-functor of potent D-modules}
\begin{comment}
Thus, we can restrict the six-functor formalism $\mathsf{IndCoh}(\calL-)$ to this class. Restricting $\mathsf{IndCoh}$ to this subcategory gives another 6-functor formalism which is the one we will use in the sequel. In particular, we record for future use that the classes $I$ and $P$ we will use are ind-schematic ind-proper maps and equivalences, respectively.
\end{comment} 
In this section, we construct the $(\infty,2)$-functor $\Pot: \mathsf{indGStk}_{\text{ft}}^{\text{op}}\to \bDGCat_{\ku}$ and show that it is $(I,P)$-biadjointable. That is, we will show the following, which, together with \zcref{cor: tTate is a 2-functor}, supplies the construction.
\begin{prop}\label{prop:indcoh-loops-6FF-setup-for-potent-d-modules}
	Restricting the six-functor formalism $\mathsf{IndCoh}(\calL -)$ of \zcref{cor:indcoh-loops-is-a-6FF} to $\mathsf{indGStk}_{\mathrm{ft}}$, we obtain an $(I,P)$-biadjointable $(\infty,2)$-functor
	\begin{align*}
		\mathsf{IndCoh}(\calL -): (\mathsf{indGStk}_{\mathrm{ft}})^{\mathrm{op}}\to \bFun(BS^1,\bCoh_t),
	\end{align*} 
	for the classes of maps $I$ and $P$ given respectively by ind-inf-schematic ind-proper maps and equivalences.
	\end{prop} 
\begin{proof}
		By \zcref{lem: IndCoh is regular for ind-finite type stacks}, which applies to $\mathcal{L}X$ by \zcref{prop: ind-ft stacks are closed under limits} and \zcref{prop: ind-ft stacks are coherent}, $\mathsf{IndCoh}(\mathcal{L}X)$ has regular $t$-structure; hence, it is an object of $\bFun(BS^1,\bCoh_t)$.

	Now, by \zcref{cor:indcoh-loops-is-a-6FF}, it suffices to show that $f^!$ is coherent for any $f$ and $f_* $ is coherent for ind-inf-schematic and ind-proper $f$.

	We will show that for $ f:X \to Y$ with $X,Y\in \mathsf{indGStk_{\text{ft}}}$, $f^!:\mathsf{IndCoh}(Y)\to \mathsf{IndCoh}(X)$ is coherent. By \zcref{prop: ind-ft stacks are closed under limits}, this allows us to conclude that $(\mathcal{L}f)^!$ is coherent because $\mathcal{L}$ is a finite limit. Note that $f^!$ is continuous and exact by the construction of the six-functor formalism for $\mathsf{IndCoh}$ (\cite[Chapter 3 Theorem 5.4.3]{GR17II}). Thus, it remains to check that it is left $t$-exact up to a shift. By \zcref{lem: bounded below iff bounded below on reduced stack}, we can assume $X$ is reduced and classical.\footnote{Since $f^!$ is continuous, there is in fact a uniform lower bound. Indeed, if there were no uniform bound, choose $\mathcal{F}_n\in \mathsf{IndCoh}(Y)^{\ge 0}$ such that $f^!\mathcal{F}_n\notin \mathsf{IndCoh}(X)^{\ge -n}$ and $\mathcal{F}=\oplus \mathcal{F}_n$ is nonnegative, so $f^!\mathcal{F}_n$ is bounded below. However, every $f^!\mathcal{F}_n$ is a retract of $f^!\mathcal{F}$, contradicting the choices for sufficiently large $n$.} By definition of $\mathsf{indGStk}_{\text{ft}}$, $X$ is thus an algebraic stack of finite type and with affine diagonal. Let $i:Y_{\text{red}}\to Y$ be the inclusion. Then there is a commutative diagram
	\[
	\begin{tikzcd}
		X\arrow[r,"f_{\text{red}}"]\arrow[dr,"f",swap]&Y_{\text{red}}\arrow[d,"i"]\\
							 &Y
	\end{tikzcd};
	\] 
	hence, by \zcref{lem: bounded below iff bounded below on reduced stack}, it is sufficient to show that $f_{\text{red}}$ is left $t$-exact up to a shift. In particular, we can assume that $Y=Y_{\text{red}}$ is reduced and classical as well. Further, by \zcref{lem: atlas detects bounds}, we can assume that $X$ is a classical affine scheme. Now, let $V\to Y$ be a smooth affine atlas, and let $U\to V\times _YX$ be a smooth affine atlas as well. Then we have the following commutative diagram, where both vertical arrows are smooth surjections. 
	\[
	\begin{tikzcd}
		U\arrow[r,"g"]\arrow[d,"q"]&V\arrow[d,"p"]\\
		X\arrow[r,"f"]&Y
	\end{tikzcd}.
	\] 
Now, since $q^!f^!\simeq g^!p^!$ and $p,q$ are smooth affine atlases for $X$ and $Y$, respectively, we conclude by \zcref{lem: atlas detects bounds} that it suffices to check that $g^!$ is left $t$-exact up to a shift. Thus, we can reduce to the case that $X$ and $Y$ are classical affine schemes. In this case, since $f:X\to Y$ is a finite type morphism of schemes, we conclude by \cite[Lemma 7.1.7]{GaitsgoryIndCoh}.

If $f\in I$, then $\mathcal{L}f$ in $I$ by \zcref{lem:forgetful-indcoh-loops-is-a-6FF}. In addition, $\mathcal{L}X,\mathcal{L}Y\in \mathsf{indGStk}_{\mathrm{ft}}$ by \zcref{prop: ind-ft stacks are closed under limits}. \zcref{lem: CW pushforward = GR pushforward} (see also \zcref{rem: CW pushforward = GR pushforward}) then shows that the pushforward $(\mathcal{L}f)_*$ is left $t$-exact, as required.
\end{proof} 
 
\begin{lem}\label{lem: IndCoh is regular for ind-finite type stacks}
	Let $X$ be a coherent ind-geometric stack. Then the standard $t$-structure on $\mathsf{IndCoh}(X)$ is regular. In particular, it is coherent.
\end{lem} 
\begin{proof}
The first two conditions of \cite[Definition 4.2.7]{Preygel} follow from \cite[Proposition 5.21]{CautisWilliams}. We will now show that $\mathsf{IndCoh}(X)^\heart$ is a locally coherent abelian category. By \cite[Proposition 5.30]{CautisWilliams}, we have an equivalence
	\begin{align*}
		\mathsf{IndCoh}(X)\simeq \mathsf{Ind}(\mathsf{Coh}(X).
	\end{align*} 
	Thus, we have an equivalence
	\begin{align*}
		\mathsf{IndCoh}(X)^\heart\simeq \mathsf{Ind}_1(\mathsf{Coh}(X)^\heart),
	\end{align*} 
	where $\mathsf{Ind}_1$ is the 1-categorical ind-completion. Since $\mathsf{Coh}( X)^\heart$ is a small abelian category, we conclude that $\mathsf{IndCoh}( X)^\heart$ is locally coherent; hence, the $t$-structure on $\mathsf{IndCoh}( X)$ satisfies the conditions of \cite[Definition 4.2.7]{Preygel} and hence is coherent. It is thus regular by Proposition 4.3.2 of \textit{op. cit.}
\end{proof}
\begin{lem}\label{lem: bounded below iff bounded below on reduced stack}
	Let $X\in \mathsf{indGStk}_{\text{ft}}$, and denote by $i:X_{\mathrm{red}}\to X$ the inclusion. Then $\mathcal{F}\in \mathsf{IndCoh}(X)^{+} $ if and only if $i^!\mathcal{F}\in \mathsf{IndCoh}(X_{\mathrm{red}})^{+}$. 
\end{lem} 
\begin{proof}
Assume that $\mathcal{F}\in \mathsf{IndCoh}(X)^+$. Suppose first that $X$ is an almost finite type algebraic stack, and let $U\to X$ be a smooth affine atlas. Then $U_{\text{red}}$ is a smooth affine atlas of $X_{\text{red}}$, and we have a Cartesian diagram
\begin{equation}\label{cd: atlas and reduction}
\begin{tikzcd}
	U_{\text{red}}\arrow[r,"j"]\arrow[d,"q"]&U\arrow[d,"p"]\\
	X_{\text{red}}\arrow[r,"i"]&X
\end{tikzcd}.
\end{equation} 
Commutativity of the above diagram implies that 
\begin{align*}
	q^!i^!\mathcal{F}\simeq j^!p^!\mathcal{F}.
\end{align*} 
By \cite[Lemma 7.1.7]{GaitsgoryIndCoh}, $j^!$ is left $t$-exact up to a finite shift. We conclude that $i^!\mathcal{F}\in \mathsf{IndCoh}(X_{\text{red}})^+$ by \zcref{lem: atlas detects bounds} below.

Now, suppose that $i^!\mathcal{F}\in \mathsf{IndCoh}(X_{\text{red}})^+$ and $X$ is an almost finite type algebraic stack. We claim that for the atlas $U$, $j^!\mathcal{G}\in \mathsf{IndCoh}(U_{\text{red}})^+$ if and only if $\mathcal{G}\in \mathsf{IndCoh}(U)^+$. This will allows us to conclude by base change along the diagram \zcref[noname]{cd: atlas and reduction} and \zcref{lem: atlas detects bounds}. One direction is supplied by the above; hence, suppose $j^!\mathcal{G}\in \mathsf{IndCoh}(U_{\text{red}})^{\ge n}$. By \cite[Proposition 5.21]{CautisWilliams}, it suffices to check that
\begin{align*}
	\Hom(\mathcal{E},\mathcal{G})\in \mathsf{Vect}^{\ge 1}
\end{align*} 
for every $\mathcal{E}\in \mathsf{Coh}(U)^{\le n-1}$. Suppose first that $\mathcal{E}\simeq j_*M[-m]$ for $M\in \mathsf{Coh}(U_{\text{red}})^\heart$ and $m\le n-1$. In this case,
\begin{align*}
	\Hom(j_*M[-m],\mathcal{G})\simeq\Hom(M[-m],j^!\mathcal{G})\in \mathsf{Vect}^{\ge 1}
\end{align*} 
by adjunction. Now, take an arbitrary $\mathcal{E}\in \mathsf{Coh}(U)^{\le n-1}$, and let 
\begin{align*}
	\mathcal{I}\coloneq \ker(H^0(\mathcal{O}_U)\to j_*\mathcal{O}_{U_{\text{red}}}).
\end{align*} 
Choose an $N$ such that $\mathcal{I}^N\simeq 0$.\footnote{Such an $N$ exists because $\mathcal{I}$ is the nilradical of $H^0(\mathcal{O}_U)$ and $U$ is Noetherian.} Now, consider a cohomology sheaf $H^m(\mathcal{E})\in \mathsf{Coh}(U)^\heart$. There is a filtration
\begin{align*}
	0=\mathcal{I}^NH^m(\mathcal{E})\subset \mathcal{I}^{N-1}H^m(\mathcal{E})\subset \cdots\subset H^m(\mathcal{E}).
\end{align*} 
Each subquotient is annihilated by $\mathcal{I}$; hence, it is isomorphic to $j_*M _m$ with $M_m\in \mathsf{Coh}(U_{\text{red}})^\heart$. Now, shift by $-m$ to obtain that, after applying the functor $\Hom(-,\mathcal{G})$, each subquotient vanishes.\footnote{Because $\mathcal{E}\in \mathsf{Coh}(U)^{\le n-1}$, $H^m(\mathcal{E})\simeq 0$ when $m\ge n$.} Induction on the length of the filtration gives the required vanishing for $\Hom(H^m(\mathcal{E})[-m],\mathcal{G})$. Now, use the Postnikov filtration for $\mathcal{E}$: Since $\mathcal{E}$ is bounded, there is a filtration
\begin{align*}
	0\simeq\tau^{\le a-1}\mathcal{E}\to \tau^{\le a}\mathcal{E}\to \cdots\to \tau^{\le b}\mathcal{E}\simeq \mathcal{E}.
\end{align*} 
For each $a\le m\le b$, there is a canonical distinguished triangle
\begin{align*}
	\tau^{\le m-1}\mathcal{E}\to \tau^{\le m}\mathcal{E}\to H^m(\mathcal{E})[-m].
\end{align*} 
Applying the functor $\Hom(-,\mathcal{G})$ to these triangles and using induction on $m$ gives that $\Hom(\mathcal{E},\mathcal{G})\in \mathsf{Vect}^{\ge 1}$ as desired.

Now, suppose that $X\in \mathsf{indGStk}_{\text{ft}}$ and $\mathcal{F}\in \mathsf{IndCoh}(X)$ such that $i^!\mathcal{F}\in \mathsf{IndCoh}(X_{\mathrm{red}})^{\ge n}$, and let $X\simeq\colim X_\alpha$ be a reasonable presentation with $j_\alpha:X_\alpha\to X$ the inclusions. By \cite[Proposition 4.7]{CautisWilliams}, the map $X_{\text{red}}\to X$ factors through some finite stage $X_\alpha$. Since $X_{\mathrm{red}}$ is an algebraic stack of finite type by definition of $\mathsf{indGStk}_{\mathrm{ft}}$, we can choose a smooth atlas $U\to X_{\mathrm{red}}$ with $U$ affine. By \cite[Proposition 4.7]{CautisWilliams}, the composite $U\to X$ factors through some $X_\alpha$. Pulling back the classical closed immersion $X_\alpha^{\text{cl}}\to X^{\text{cl}}$ along $X_{\text{red}}\to X^{\text{cl}}$, the resulting closed substack of $X_{\text{red}}$ contains the atlas $U$; hence, it is equivalent to $X_{\text{red}}$. This same argument applies for any $\beta\geq \alpha$. Now, fix $\beta \ge \alpha$, and let $\mathcal{E}\in \mathsf{Coh}(X_\beta)^{\le n-1}$. Let $\mathcal{I}_\beta$ be the nilpotent\footnote{This can be checked on a Noetherian smooth affine atlas.} ideal of $X_{\mathrm{red}}\simeq (X_{\beta})_{\mathrm{red}}$ in $X_\beta^{\mathrm{cl}}$. Thus, as above, its powers equip each cohomology sheaf $H^m(\mathcal{E})$ with a finite filtration with successive quotients of the form $i_{\beta*}M$ for $M\in \mathsf{Coh}(X_{\mathrm{red}})^\heart$. For $m\le n-1$, we have, by adjunction,
\begin{align*}
	\Hom(i_{\beta*}M[-m],j^!_\beta \mathcal{F})\simeq\Hom(M[-m],i^!\mathcal{F})\in \mathsf{Vect}^{\ge 1},
\end{align*} 
since $i^!\mathcal{F}\in \mathsf{IndCoh}(X_{\mathrm{red}})^{\ge n}$. As in the argument above, we apply this first to the filtration of each $H^m(\mathcal{E})$ and then to the Postnikov filtration of $\mathcal{E}$ to obtain that $\Hom(\mathcal{E},j^!_\beta \mathcal{F})\in \mathsf{Vect}^{\ge 1}$. By \cite[Proposition 5.30]{CautisWilliams}, such $\mathcal{E}$ generate $\mathsf{IndCoh}(X_\beta)^{\le n-1}$ under filtered colimits; hence, $j_\beta^!\mathcal{F}\in \mathsf{IndCoh}(X_\beta)^{\ge n}$. Passing to a cofinal subsystem implies that we can assume this is true for all $\beta$. Hence, \cite[Proposition 5.21]{CautisWilliams} implies then that $\mathcal{F}\in \mathsf{IndCoh}(X)^{\geq n}$. Conversely, if $\mathcal{F}\in \mathsf{IndCoh}(X)^{\geq n}$, then \textit{loc. cit.} applied to the truncated geometric substack $X_{\text{red}}\to X$ allows us to conclude that $i^!\mathcal{F}\in \mathsf{IndCoh}(X_{\text{red}})^{\geq n}$.
\end{proof}
\begin{lem}\label{lem: atlas detects bounds}
	Let $X$ be an algebraic stack almost of finite type and with affine diagonal, and let $p:U\to X$ be a smooth affine atlas. Then $\mathcal{F}\in \mathsf{IndCoh}(X)^+$ if and only if $p^!\mathcal{F}\in \mathsf{IndCoh}(U)^+$.
\end{lem} 
\begin{proof}
	Let $U_d\subset U$ in the locus where $p$ has relative dimension $d$. These loci are open and closed because the relative dimension of a smooth morphism is locally constant. Since $U$ is affine, it is quasi-compact, so only finitely many $U_d$ are nonempty. By \cite[Equation (3.12)]{DrinfeldGaitsgory} and the remarks surrounding it, on each $U_d$ we have
	\begin{align*}
		(p^!\mathcal{F})|_{U_d}\simeq K_{U_d /X}\otimes (p^*\mathcal{F})|_{U_d}[d],
	\end{align*} 
	where $K_{U_d/X}$ is the relative canonical bundle of $U_d\to X$. Tensoring with the relative canonical bundle is $t$-exact. Thus, since only finitely many shifts $d$ occur, $p^!\mathcal{F}$ is bounded below if and only if $p^*\mathcal{F}$ is bounded below.
\end{proof}
	
\subsection{The lax symmetric monoidal structure}
\label{ssec: lax symmetric monoidal structure}
In this section, we construct the lax symmetric monoidal structure on the $(\infty,2)$-functor of \zcref{prop:indcoh-loops-6FF-setup-for-potent-d-modules}.

In order to construct the six-functor formalism $\Pot$ in question, one would like to now show that the $(\infty,2)$-functor $(-)^{\tT}$ of \zcref{cor: tTate is a 2-functor} is lax symmetric monoidal. However, the tensor product of categories in $\bCoh_t$ does not necessarily inherit a coherent $t$-structure \cite[Example 21]{tensor-product-of-coherent-is-not-coherent}, so this is not a well-posed question. We will instead use the formalism of $\infty$-operads of \zcref{ssec:operads} and construct $(-)^{\tT}$ as a map of $\infty$-operads. Its composition with $\mathsf{IndCoh}(\mathcal{L}-)$ will then also be a map of $\infty$-operads between symmetric monoidal categories, thus endowing $\Pot$ with a lax symmetric monoidal structure.

\subsubsection{A lax symmetric monoidal structure on $\mathsf{IndCoh}(\mathcal{L}-)^{S^1}$}

\begin{prop}\label{prop: IndCoh(L-)^S^1 as a map of operads}
	The functor
	\begin{align*}
		\mathsf{IndCoh}(\mathcal{L}-): \left((\mathsf{indGStk}_{\mathrm{ft}})^{\mathrm{op}}\right)^\otimes \to \Fun(BS^1,\Coh_{t})^\otimes  ,
	\end{align*} 
	is a map of $\infty$-operads, where $\Fun(BS^1,\Coh_t)^\otimes $ is the $\infty$-operad of \zcref{con: S^1-equivariant Coh_t operad}.
\end{prop} 
\begin{proof}
	By \zcref{cor:indcoh-loops-is-a-6FF}, we have a map of $\infty$-operads
	\begin{align*}
		\mathsf{IndCoh}(\mathcal{L}-):\left((\mathsf{indGStk}_{\mathrm{ft}})^{\mathrm{op}}\right)^\otimes \to \Fun(BS^1,\DGCat_k)^\otimes .
	\end{align*} 
We will show that this map lifts through postcomposition with the forgetful functor
	\begin{align*}
		\oblv_t^\otimes :\Coh_{t}^\otimes \to \DGCat_k^\otimes 
	\end{align*} 
	of \zcref{con: Coh_t and Reg_t operads}, which is a map of $\infty$-operads by \zcref{prop: Coh_t and Reg_t are operads}. Indeed, this follows from \zcref{lem: exterior product in Coh_t} applied repeatedly. For clarity, we will spell this out explicitly. Let $\alpha:\braket{n}\to \braket m$ be a map in $N(\mathsf{Fin}_*)$, and let
	\begin{align*}
		(X_1,\ldots,X_n)\to (Y_1,\ldots,Y_m)
	\end{align*} 
	be a map in $(\mathsf{indGStk}_{\text{ft}}^{\text{op}})^\otimes $ lying over $\alpha$. For each $j\in \{1,\ldots,m\} $, such a morphism is specified by maps of stacks
	\begin{align*}
		f_j:Y_j\to \prod_{i\in\alpha^{-1}(j)}X_i,
	\end{align*} 
	and the $j^{\text{th}}$ component operation of its image under $A^\otimes $ is
	\begin{align*}
		\prod_{i\in\alpha^{-1}(j)}\mathsf{IndCoh}(\mathcal{L}X_i)\xrightarrow{E_{\alpha^{-1}(j)}}\mathsf{IndCoh}\left(\prod_{i\in\alpha^{-1}(j)}\mathcal{L}X_i\right)\xrightarrow{(\mathcal{L}f_j)^!}\mathsf{IndCoh}(\mathcal{L}Y_j),
	\end{align*} 
	where $E_{\alpha^{-1}(j)}$ is the exterior product functor. Applying \zcref{lem: exterior product in Coh_t} repeatedly to this functor, we see that this operation is a well-defined operation in $\Fun(BS^1,\Coh_t)^\otimes  $, where we note that $(\mathcal{L}f_j)^{!}$ is a 1-morphism in $\bFun(BS^1,\bCoh_{t})$ by \zcref{prop:indcoh-loops-6FF-setup-for-potent-d-modules}. Note that the same argument applies if $\alpha^{-1}(j)$ is empty, where only \zcref{prop:indcoh-loops-6FF-setup-for-potent-d-modules} is required.
\end{proof} 
\begin{lem}\label{lem: exterior product in Coh_t}
	Let $X,Y\in \mathsf{indGStk}_{\mathrm{ft}}$. Then there is a separately continuous and $S^1$-equivariant exterior product functor
	\begin{align*}
		E_{X,Y}:\mathsf{IndCoh}(\mathcal{L}X)\times \mathsf{IndCoh}(\mathcal{L}Y)\to \mathsf{IndCoh}(\mathcal{L}X\times \mathcal{L}Y),
	\end{align*} 
where the $S^1$-actions on the source and target are the diagonal ones. Moreover, the functor $E_{X,Y}$ is jointly left $t$-exact up to a shift.
\end{lem} 
\begin{proof}
	By \zcref{cor:indcoh-loops-is-a-6FF}, there is a separately continuous and $S^1$-equivariant exterior product functor
	\begin{align*}
		\mathsf{IndCoh}(\mathcal{L}X)\times \mathsf{IndCoh}(\mathcal{L}Y)\to \mathsf{IndCoh}(\mathcal{L}X\times \mathcal{L}Y),
	\end{align*} 
	and we will show it is jointly left $t$-exact up to a shift.
	By \cite[Proposition 5.30]{CautisWilliams}, \zcref{prop: ind-ft stacks are closed under limits}, and \zcref{prop: ind-ft stacks are coherent}, if $\mathcal{F}\in \mathsf{IndCoh}(\mathcal{L}X)^{\ge a}$ and $\mathcal{G}\in \mathsf{IndCoh}(\mathcal{L}Y)^{\ge b}$, then
	\begin{align*}
		\mathcal{F}\simeq\colim_i \mathcal{F}_i,\qquad \mathcal{G}\simeq\colim_j \mathcal{G}_j,
	\end{align*} 
	where $\mathcal{F}_i\in \mathsf{Coh}(\mathcal{L}X)^{\ge a}$ and $\mathcal{G}_j\in \mathsf{Coh}(\mathcal{L}Y)^{\ge b}$. Since exterior product preserves these colimits separately and $\mathsf{IndCoh}(\mathcal{L}X\times \mathcal{L}Y)^{\ge 0}$ is closed under filtered colimits, it suffices to prove the non-equivariant assertion for coherent $\mathcal{F},\mathcal{G}$. Now, choose reasonable presentations $\mathcal{L}X\simeq\colim X_\alpha$ and $\mathcal{L}Y\simeq\colim Y_\beta$ with inclusion $i_\alpha:X_\alpha\to \mathcal{L}X$ and $j_\beta:Y_\beta\to \mathcal{L}Y$. By \cite[Proposition 5.21]{CautisWilliams}, there are sheaves $\mathcal{F}_\alpha\in \mathsf{Coh}(X_\alpha)^{\ge a}$ and $\mathcal{G}_\beta\in \mathsf{Coh}(Y_\beta)^{\ge b}$ such that $\mathcal{F}\simeq i_{\alpha*}\mathcal{F}_\alpha$ and $\mathcal{G}\simeq j_{\beta*}\mathcal{G}_\beta$. Moreover, we have
	\begin{align*}
		\mathcal{F}\boxtimes \mathcal{G}\simeq(i_{\alpha}\times j_\beta)_*(\mathcal{F}_\alpha\boxtimes \mathcal{G}_\beta).
	\end{align*} 
	Since the pushforwards along $i_\alpha$ and $j_\beta$, hence along $i_\alpha \times j_\beta$, are $t$-exact by \textit{loc. cit.}, it suffices to prove the non-equivariant assertion for the stacks $X_\alpha$ and $Y_\beta$. Thus, let $p:U\to X_\alpha$ and $q:V\to Y_\beta$ be faithfully flat affine atlases. For coherent objects, the exterior product is the ordinary exterior product of quasi-coherent sheaves; hence,
	\begin{align*}
		(p\times q)^*(\mathcal{F}_\alpha\boxtimes \mathcal{G}_\beta)\simeq p^*\mathcal{F}_\alpha\boxtimes q^*\mathcal{G}_\beta.
	\end{align*} 
	Since flat quasi-coherent pullback is $t$-exact, the two factors have lower bounds $a$ and $b$, respectively. Now, since $U=\Spec A$ and $V=\Spec B$, hence $U\times V$, are affine, the exterior product above corresponds to an ordinary tensor product $M\otimes N$ as an $(A\otimes B)$-module. Moreover, its cohomology splits as
	\begin{align*}
		H^k(M\otimes N)\simeq\bigoplus _{i+j=k}H^i(M)\otimes H^j(N).
	\end{align*} 
	When $k<a+b$, the cohomology vanishes, so the exterior product has a lower bound $a+b$. Faithfully flat pullback detects this bound.
\end{proof}
\subsubsection{A lax symmetric monoidal structure on $S^1$-invariants}
\begin{lem}\label{lem: S^1-invariants is lax symmetric monoidal}
	Taking $S^1$-invariants defines a lax symmetric monoidal functor
	\begin{align*}
		(-)^{S^1}:\Fun(BS^1,\DGCat_k)\to \DGCat_{k[u]}.
	\end{align*} 
\end{lem} 
\begin{proof}
	The $k$-linear version of $(-)^{S^1}$ is right adjoint to the trivial action functor
	\begin{align*}
		\triv:\DGCat_k\to \Fun(BS^1,\DGCat_k),
	\end{align*} 
	which is symmetric monoidal with the diagonal action on tensor products. Thus, \cite[Corollary 7.3.2.7]{HA} with $\mathscr{O}^\otimes \coloneq N(\mathsf{Fin}_*)$, the commutative $\infty$-operad, induces the lax symmetric monoidal structure on the $k$-linear functor of $S^1$-invariants. Thus, $S^1$-invariants gives a lax symmetric monoidal functor (cf. \cite[Section 3.3.3]{HA})
    \begin{align*}
        \Fun(BS^1,\DGCat_k)=\mathsf{Mod}_{\Vect}(\Fun(BS^1,\DGCat_k))\to\mathsf{Mod}_{\Vect^{S^1}}(\DGCat_k),
    \end{align*}
    and there is a natural restriction functor $\mathsf{Mod}_{\Vect^{S^1}}(\DGCat_k)\to \DGCat_{k[u]}$ induced by restriction along $\DGCat_{k[u]}\to\Vect^{S^1}$. This functor is canonically lax symmetric monoidal, so the composition
    \begin{align*}
        \mathsf{Mod}_{\Vect}(\Fun(BS^1,\DGCat_k))\to\mathsf{Mod}_{\Vect^{S^1}}(\DGCat_k)\to\DGCat_{k[u]}
    \end{align*}
    is the desired lax symmetric monoidal $S^1$-invariants functor.
\begin{comment}
%explanation of part of the induced lax symmetric monoidal structure below.	
	The right adjoint $(-)^{S^1}$ as a map of $\infty$-operads includes the following data. For the two inert maps $\braket{2}\to \braket{1}$ in $N(\mathsf{Fin}_*)$, this supplies an identification $$\left( (-)^{S^1} \right) ^{\otimes }(\mathscr{C},\mathscr{D})\simeq(\mathscr{C}^{S^1},\mathscr{D}^{S^1}),$$ where each maps picks one of the two entries of a pair. Now, for the active morphism,
	\begin{align*}
		\xi:\braket{2}\to& \braket{1}\\
		1,2\mapsto&1,
	\end{align*} 
	the tensor product is encoded by a coCartesian arrow lying over $\xi$
	\begin{align*}
		(\mathscr{C},\mathscr{D})\to \mathscr{C}\otimes \mathscr{D}.
	\end{align*} 
	Applying $\left( (-)^{S^1} \right) ^\otimes $ to this arrow gives a map
	\begin{align*}
		\left( \mathscr{C}^{S^1},\mathscr{D}^{S^1} \right) \to (\mathscr{C}\otimes \mathscr{D})^{S^1}
	\end{align*} 
	also over $\xi$. An arrow over $\xi$ with these endpoints thus gives a map
	\begin{align*}
		\mathscr{C}^{S^1}\otimes \mathscr{D}^{S^1}\to (\mathscr{C}\otimes \mathscr{D})^{S^1}.
	\end{align*} 
	One can obtain the rest of the coherence data by using higher active and inert maps in $N(\mathsf{Fin}_*)$. 
\end{comment} 
\end{proof}
\begin{cor}\label{cor: S^1-invariants is an operad map}
	Taking $S^1$-invariants defines a map of $\infty$-operads
	\begin{align*}
		(-)^{S^1}:\Fun(BS^1,\Coh_t)^\otimes \to \Coh_{t,k[u]}^\otimes ,
	\end{align*} 
	where $\Fun(BS^1,\Coh_t)^\otimes $ is the $\infty$-operad of \zcref{con: S^1-equivariant Coh_t operad}.
\end{cor} 
\begin{proof}
	Let $F:\prod \mathscr{C}_i\to \mathscr{D}$ be an $S^1$-equivariant $n$-ary operation in $\Coh_t^\otimes $, and write $U_{\mathscr{C}}: \mathscr{C}^{S^1}\to \mathscr{C}$ for the forgetful functor. The induced operation satisfies
	\begin{align*}
		U_{\mathscr{D}}F^{S^1}(c_1,\ldots,c_n)\simeq F(U_{\mathscr{C}_1}c_1,\ldots,U_{\mathscr{C}_n}c_n).
	\end{align*} 
	The invariant categories $\mathscr{C}_i^{S^1}$ have coherent $t$-structures and the forgetful functors $U_{\mathscr{C}_i}$ detect their $t$-bounds by \cite[Lemma 4.5.3]{Preygel}. Hence, by \zcref{con: Coh_t and Reg_t operads}, the $\infty$-operad map of \zcref{lem: S^1-invariants is lax symmetric monoidal} restricts as desired.
\end{proof} 
\subsubsection{Potent $\mathcal{D}$-modules are lax symmetric monoidal}
We are now ready for the main result of this section.
\begin{prop}\label{prop: D^pot is lax symmetric monoidal}
	The functor
	\begin{align*}
		\Pot: \mathsf{indGStk}_{\mathrm{ft}}^{\mathrm{op}}\to \DGCat_{\ku}
	\end{align*} 
	extends to a map of $\infty$-operads
	\begin{align*}
		\Pot: \left(\mathsf{indGStk}_{\mathrm{ft}}^{\mathrm{op}}\right)^\otimes \to \DGCat_{\ku}^\otimes .
	\end{align*} 
	In particular, $\Pot$ is lax symmetric monoidal.
\end{prop}
\begin{rem}\label{rem: lax symmetric monoidal structure is exterior product}
	The lax symmetric monoidal structure on $\Pot$ is given by exterior tensor product. More precisely, the proof of \zcref{prop: D^pot is lax symmetric monoidal} below shows that if $X,Y\in \mathsf{indGStk}_{\text{ft}}$, then the induced map
	\begin{align*}
		\Pot(X)\otimes \Pot(Y)\to \Pot(X\times Y)
	\end{align*} 
	is the exterior product functor
	\begin{align*}
		\mathcal{F}\otimes \mathcal{G}\mapsto \mathcal{F}\boxtimes \mathcal{G}.
	\end{align*} 
\end{rem} 
\begin{proof}
	By \zcref{prop: IndCoh(L-)^S^1 as a map of operads} and \zcref{cor: S^1-invariants is an operad map}, we have the following map of $\infty$-operads
	\begin{align*}
		\mathsf{IndCoh}(\mathcal{L}-)^{S^1}:(\mathsf{indGStk}_{\text{ft}}^{\text{op}})^\otimes \to \Coh_{t,k[u]}^\otimes .
	\end{align*} 
	Thus, by \zcref{lem: P is equivalent to Coh_t} and \zcref{lem: P -> Reg_t is a map of operads}, we have the following composition of maps of $\infty$-operads
	\begin{align}\label{eq: k-linear D^pot}
		(\mathsf{indGStk}_{\text{ft}}^{\text{op}})^\otimes \xrightarrow{\mathsf{IndCoh}(\mathcal{L}-)^{S^1}} \Coh_{t,k[u]}^\otimes \xrightarrow{p^{-1}}\mathscr{P}^\otimes \xrightarrow{q} \Reg_{t,k[u]}^\otimes\hook\Coh_{t,k[u]}^\otimes \xrightarrow{\oblv_t^\otimes }  \DGCat_{k[u]}^\otimes ,
	\end{align} 
	where the map $\oblv_t^\otimes :\Coh_{t,k[u]}^\otimes \to \DGCat_{k[u]}^\otimes $ is the forgetful functor of \zcref{con: Coh_t and Reg_t operads}, shown to be a map of $\infty$-operads in \zcref{lem: P -> Reg_t is a map of operads}. Moreover, the image in $\DGCat_{k[u]}^\otimes $ of an object $X\in \mathsf{indGStk}_{\text{ft}}$ is precisely
	\begin{align*}
		\mathscr{R}(\mathsf{IndCoh}(\mathcal{L}X)^{S^1})
	\end{align*} 
	as a $k[u]$-linear category. Now, extending scalars to $\ku$ via the strict symmetric monoidal functor $\DGCat_{k[u]}^\otimes \to \DGCat_{\ku}^\otimes $ produces a map of $\infty$-operads
	\begin{align*}
		\Pot: \left(\mathsf{indGStk}_{\mathrm{ft}}^{\mathrm{op}}\right)^\otimes \to \DGCat_{\ku}^\otimes .
	\end{align*} 
	Since both $\mathsf{indGStk}_{\text{ft}}^{\text{op}}$ and $\DGCat_{\ku}$ are symmetric monoidal categories, this map thus defines a lax symmetric monoidal structure on $\Pot$.
\end{proof}

	\subsection{The main theorem}
	We now combine the results of the previous sections to construct a six-functor formalism on potent $\mathcal{D}$-modules.
\begin{thm}\label{thm:potent-D-modules-is-a-6FF}
	The $(\infty,2)$-functor
	\begin{align*}
		\Pot:\mathsf{indGStk}_{\mathrm{ft}}^{\mathrm{op}}\to \bDGCat_{\ku}
	\end{align*} 
	is lax symmetric monoidal $(I,P)$-biadjointable, in the sense of \cite[Definition 5.13]{CLL_six-functor-formalisms}, for the classes of maps $E=I$ and $P$ given respectively by ind-inf-schematic ind-proper maps and equivalences. Hence, $\Pot$ defines a six functor formalism
    \begin{align*}
        (\Pot)^\otimes:\bSpan(\mathsf{indGStk}_{\mathrm{ft}},E)^\otimes_{P,I}\to\bDGCat_{\ku}^\otimes.
    \end{align*}
\end{thm} 
\begin{proof}
	We remark that the conditions of \cite[Theorem 5.14]{CLL_six-functor-formalisms} on $\mathsf{indGStk}_{\text{ft}}$ and $I,P,E$ are satisfied because $\mathsf{indGStk}_{\text{ft}}$ is a full subcategory of $\mathsf{PreStk}_{\text{laft}}$ closed under finite limits by \zcref{prop: ind-ft stacks are closed under limits}. The conditions on the latter, for the classes of maps $I$ and $P$, were checked in \zcref{prop: indcoh is a 6-functor formalism}. By \zcref{prop: D^pot is lax symmetric monoidal}, $\Pot$ is lax symmetric monoidal. Furthermore, 
	\begin{align*}
		\mathsf{IndCoh}(\mathcal{L}-): (\mathsf{indGStk}_{\text{ft}})^{\mathrm{op}}\to \bFun(BS^1,\bCoh_t)
	\end{align*} 
	is $(I,P)$-biadjointable  by \zcref{prop:indcoh-loops-6FF-setup-for-potent-d-modules}. Thus, using that regularized Tate invariants $(-)^{\tT}$ is an $(\infty,2)$-functor (\zcref{cor: tTate is a 2-functor}) and that $(-)^{\tT}$ preserves adjunctions \zcref{cor: tTate preserves adjunctions}, \cite[Definition 5.13(1)]{CLL_six-functor-formalisms} is satisfied. It remains to check Condition (2) of \textit{loc. cit.} Let $i:X'\to X$ be an ind-inf-schematic ind-proper map in $\mathsf{indGStk}_{\text{ft}}$, and let $Y\in \mathsf{indGStk}_{\text{ft}}$. By \zcref{rem: lax symmetric monoidal structure is exterior product}, it suffices to check that the commutative square
	\begin{equation}\label{eq: left adjointable D^pot}
		\begin{tikzcd}[column sep=large]
		\Pot(X)\otimes \Pot(Y)\arrow[r,"i^{!,\text{pot}}\otimes \id_Y^{!,\text{pot}}"]\arrow[d,"E^{\text{pot}}_{X,Y}"]&\Pot(X')\otimes \Pot(Y)\arrow[d,"E^{\text{pot}}_{X',Y}"]\\
		\Pot(X\times Y)\arrow[r,"(i\times \id_Y)^{!,\text{pot}}"]&\Pot(X'\times Y)
	\end{tikzcd}
	\end{equation}
	is horizontally left adjointable, where $E^{\text{pot}}_{X,Y},E^{\text{pot}}_{X',Y}$ are the exterior product functors. For $\mathcal{F}\in \Pot(X')$ and $\mathcal{G}\in \Pot(Y)$, the left mate of \zcref[noname]{eq: left adjointable D^pot} is the map
	\begin{align*}
		m_{\mathcal{F},\mathcal{G}}:(i\times \id_Y)_*^{\mathrm{pot}}(\mathcal{F}\boxtimes \mathcal{G})\to i_*^{\mathrm{pot}}\mathcal{F}\boxtimes \mathcal{G}
	\end{align*} 
	whose adjunct under $(i\times \id_Y)_*^{\mathrm{pot}}\dashv (i\times \id_Y)^{!,\mathrm{pot}}$ is
	\begin{align*}
		\mathcal{F}\boxtimes \mathcal{G}\xrightarrow{\eta_i(\mathcal{F})\boxtimes \id_{\mathcal{G}}}i^{!,\mathrm{pot}}i_*^{\mathrm{pot}}\mathcal{F}\boxtimes \mathcal{G}\xrightarrow{\simeq}(i\times \id_Y)^{!,\mathrm{pot}}(i_*^{\mathrm{pot}}\mathcal{F}\boxtimes \mathcal{G}).
	\end{align*} 
	Here, $\eta_i$ is the unit of the adjunction $i_*^{\mathrm{pot}}\dashv i^{!,\mathrm{pot}}$, and the second arrow is the naturality equivalence of \zcref[noname]{eq: left adjointable D^pot}. 

	Fix $\mathcal{G}\in \mathsf{Coh}(\mathsf{IndCoh}(\mathcal{L}Y)^{S^1})$, and denote its image in $\Pot(Y)$ by $\mathcal{G}'$. We will prove the required adjointability by showing that the relevant mate is invertible at $(\mathcal{F},\mathcal{G}')$ for every $\mathcal{F}\in \Pot(X')$. Since the $\mathcal{G}'$ generate $\Pot(Y)$ and both sides of the mate preserve colimits in the second variable, this allows us to conclude that the mate is invertible. By \zcref{lem: exterior product in Coh_t}, exterior product with $\mathcal{G}$ gives an $S^1$-equivariant functor
	\begin{align*}
		E_{X,Y}(-,\mathcal{G}): \mathsf{IndCoh}(\mathcal{L}X)\to \mathsf{IndCoh}(\mathcal{L}X\times \mathcal{L}Y).
	\end{align*} 
	By \textit{loc. cit.} and since $\mathcal{G}$ is bounded, this functor is coherent. Its naturality square for $( \mathcal{L}i)^!$ is horizontally left adjointable by \zcref{cor:indcoh-loops-is-a-6FF}. Indeed, its mate is the mate of \textit{loc. cit.} evaluated at $\mathcal{G}$ in the second variable. Now, apply the $(\infty,2)$-functor $\tT$ to this square, and use \cite[Observation 1.6.1]{HaineAnabel}. Since $\mathcal{G}$ is bounded, it defines an $S^1$-equivariant nullary operation. By the nullary case of \zcref{lem: multilinear regularization functor}, regularization sends it to its lift along $\mathscr{R}(\mathsf{IndCoh}(\mathcal{L}Y)^{S^1})^+\simeq (\mathsf{IndCoh}(\mathcal{L}Y)^{S^1})^{+}$. Inverting $u$ sends this lift to $\mathcal{G}'$. Therefore, the operadic construction of \zcref{prop: D^pot is lax symmetric monoidal} identifies exterior product with $\mathcal{G}$, after applying $(-)^{\tT}$, with exterior product with $\mathcal{G}'$. Moreover, this identification is compatible with the pullback naturality squares for $(\mathcal{L}i)^!$ and $i^{!,\mathrm{pot}}$. Since $(-)^{\tT}$ preserves the units and counits of the horizontal adjunctions by \zcref{cor: tTate is a 2-functor}, it follows that 
	\begin{align*}
		m_{\mathcal{F},\mathcal{G}'}\simeq (m'_{\mathcal{G}})^{\tT}_{\mathcal{F}},
	\end{align*} 
	where $m'_{\mathcal{G}}$ denotes the unary functor $\mathcal{F}\mapsto m'_{\mathcal{F},\mathcal{G}}$, with $m'$ the invertible mate of the square corresponding to \zcref[noname]{eq: left adjointable D^pot} for the six-functor formalism of \zcref{cor:indcoh-loops-is-a-6FF}. Since the right-hand side is invertible, we conclude that $m_{\mathcal{F},\mathcal{G}'}$ is invertible as desired.
	\end{proof}

\section{Compatibility with ordinary $\mathcal{D}$-modules}
\subsection{A six-functor formalism for right crystals}
\label{sssec: crystal 3FF}
In this section, we define a six-functor formalism for periodic right crystals by extension of scalars, and we show it is equivalent to the six-functor formalism given by $\tT$ invariants for the trivial $S^1$-action on right crystals. This latter six-functor formalism is the one we will compare to potent $\mathcal{D}$-modules in \zcref{ssec: potent D-modules=crystals} below.
\subsubsection{A three-functor formalism for periodic crystals}
Consider the category of separated derived schemes almost of finite type $\mathsf{Sch}_{\text{aft}}$. For $X\in \mathsf{Sch}_{\text{aft}}$, let
\begin{align*}
	\mathsf{Crys}^r(X)\coloneq &\mathsf{IndCoh}(X_{\mathrm{dR}}),\\
	\mathsf{Crys}^r_{\mathbb{Z} /2}(X)\coloneq& \mathsf{IndCoh}(X_{\text{dR}})\otimes \ku\mod.
\end{align*} 
By \cite[Chapter 4 Theorem 2.1.2]{GR17II}, $\mathsf{Crys}^r$ defines (by restriction) a lax symmetric monoidal $(\infty,2)$-functor 
\begin{align*}
	(\mathsf{Sch}_{\mathrm{aft}})^{\mathrm{op}}\to \bDGCat_{k}.
\end{align*} 
Since extension of scalars is $(\infty,2)$-functorial and lax symmetric monoidal, we obtain a lax symmetric monoidal $(\infty,2)$-functor
\begin{align*}
	\mathsf{Crys}^r_{\mathbb{Z} /2}:(\mathsf{Sch}_{\mathrm{aft}})^{\mathrm{op}}\to \bDGCat_{\ku}.
\end{align*} 
\subsubsection{Another three-functor formalism for periodic crystals}
\label{sssec: crystal 3FF 2}
By \cite[Proposition 4.4.4]{Preygel}, there is an a priori different assignment we could define, which sends $X\in \mathsf{Sch}_{\mathrm{aft}}$ to
\begin{align*}
	\mathsf{Crys}^r(X)^{\tT},
\end{align*} 
where we equip $\mathsf{Crys}^r(X)$ with the trivial $S^1$-action. We will show below that this assignment also defines a three-functor formalism; in fact, these three-functor formalisms are equivalent, which we will prove in \zcref{lem: periodic crystals are tate invariants of crystals} below.
\begin{lem}\label{lem: crystals factor through Coh_t}
	The composition of $(\infty,2)$-functors
	\begin{align*}
		(	\mathsf{Sch}_{\mathrm{aft}})^{\mathrm{op}}\xrightarrow{(-)_{\mathrm{dR}}} (\mathsf{PreStk}_{\mathrm{laft}})^{\mathrm{op}}\xrightarrow{\mathsf{IndCoh}} \bDGCat_k
	\end{align*} 
	factors through the forgetful functor
	\begin{align*}
		\oblv_t:\bCoh_t\to \bDGCat_k.
	\end{align*} 
\end{lem} 
\begin{proof}
	\cite[Proposition 4.4.4]{Preygel} shows that $\mathsf{IndCoh}(X_{\mathrm{dR}})$ is coherent for $X\in \mathsf{Sch}_{\mathrm{aft}}$. We now claim that given a map $f:X\to Y$ in $\mathsf{Sch}_{\mathrm{aft}}$, then $f^{!}_{\mathrm{dR}}: \mathsf{Crys}^r(Y)\to \mathsf{Crys}^r(X)$ is coherent. To see this, note that the right crystal $t$-structure is defined via the pullback along the canonical map $q_X:X\to X_{\mathrm{dR}}$ \cite[Section 4.2.2]{GRCrystals}. More precisely, the definition of \textit{loc. cit.} gives that $\mathcal{F}\in \mathsf{Crys}^r(X)^{\ge n}$ if and only if $q_X^!\mathcal{F}\in \mathsf{IndCoh}(X)^{\ge n}$ and similarly for $q_Y$. Now, we have a commutative diagram
	\[
	\begin{tikzcd}
		\mathsf{IndCoh}(Y_{\mathrm{dR}})\arrow[r,"q_Y^!"]\arrow[d,"f_{\mathrm{dR}}^!"]&\mathsf{IndCoh}(Y)\arrow[d,"f^!"]\\
		\mathsf{IndCoh}(X_{\mathrm{dR}})\arrow[r,"q_X^!"]&\mathsf{IndCoh}(X)
	\end{tikzcd},
	\]
	so we conclude that $f_{\mathrm{dR}}^!$ is left $t$-exact up to a shift, hence coherent, by \cite[Lemma 7.1.7]{GaitsgoryIndCoh} applied to $f^!$. Since 2-morphisms in $\bCoh_t$ are the same as those in $\bDGCat_k$, we obtain the desired factorization.
\end{proof} 
By \zcref{lem: crystals factor through Coh_t}, there is a well-defined $(\infty,2)$-functor
\begin{align*}
	\mathsf{Crys}^r(-)^{\tT}: (\mathsf{Sch}_{\mathrm{aft}})^{\mathrm{op}}\to \bDGCat_{\ku}.
\end{align*} 
\begin{lem}\label{lem: tTate of crystals is lax symmetric monoidal}
	The underlying $(\infty,1)$-functor of the $(\infty,2)$-functor $\mathsf{Crys}^r(-)^{\tT}$ has a canonical lax symmetric monoidal structure.
\end{lem} 
\begin{proof}
	We first show that the exterior product on right crystals is left $t$-exact up to a shift. This will imply, by a similar to argument to that in \zcref{prop:indcoh-loops-6FF-setup-for-potent-d-modules}, that the underlying $(\infty,1)$-functor of the composite 
	\begin{align*}
		(	\mathsf{Sch}_{\mathrm{aft}})^{\mathrm{op}}\xrightarrow{(-)_{\mathrm{dR}}} (\mathsf{PreStk}_{\mathrm{laft}})^{\mathrm{op}}\xrightarrow{\mathsf{IndCoh}} \bCoh_t
	\end{align*} 
	of \zcref{lem: crystals factor through Coh_t} induces a map of $\infty$-operads
	\begin{align*}
		\left((\mathsf{Sch}_{\mathrm{aft}})^{\mathrm{op}}\right)^\otimes \to \Coh_t^\otimes .
	\end{align*} 
	Indeed, naturality of the exterior product for ind-coherent sheaves gives, for $f:X\to Y$ a map in $\mathsf{Sch}_{\mathrm{aft}}$ and $q_X:X\to X_{\mathrm{dR}}$,
	\begin{align*}
		q^!_{X\times Y}(\mathcal{F}\boxtimes \mathcal{G})\simeq q_X^!\mathcal{F}\boxtimes q_Y^!\mathcal{G}.
	\end{align*} 
	Hence, the exterior product on right crystals is left $t$-exact up to a shift if and only if it is so on ind-coherent sheaves on almost finite type schemes, which we proved above in \zcref{lem: exterior product in Coh_t}. Moreover, the same bound holds on circle invariants for the trivial $S^1$-action by forgetting equivariance, as in the proof of \textit{loc. cit.} Thus, we have a well defined map of $\infty$-operads
	\begin{align*}
		\left((\mathsf{Sch}_{\mathrm{aft}})^{\mathrm{op}}\right)^\otimes \to& \Coh_{t,k[u]}^\otimes \\
		X\mapsto &\mathsf{Crys}^r(X)^{S^1}.
	\end{align*} 
	Composing it with the map $\oblv_t^\otimes \circ qp^{-1}$ of \zcref[noname]{eq: k-linear D^pot} and extending scalars along $k[u]\to \ku$, we obtain the desired lax symmetric monoidal structure.
\end{proof} 
\subsubsection{Periodic crystals are $\tT$ invariants of crystals}
We now show that the two a priori different three-functor formalisms defined in \zcref{sssec: crystal 3FF} and \zcref{sssec: crystal 3FF 2} are equivalent.
\begin{lem}\label{lem: periodic crystals are tate invariants of crystals}
    There is a natural equivalence of $(\infty,2)$-functors $(\mathsf{Sch}_{\mathrm{aft}})^{\mathrm{op}}\to \bDGCat_{\ku}$,
    \begin{align*}
        \xi:\mathsf{Crys}^r_{\mathbb{Z}/2}\xRightarrow{\simeq}\mathsf{Crys}^r(-)^{\tT}.
    \end{align*}
    Moreover, this equivalence induces an equivalence of the underlying lax symmetric monoidal $(\infty,1)$-functors.
\end{lem}
\begin{proof}
	Equipping an object with the trivial $S^1$-equivariant structure defines a functor
	\begin{align*}
		t_X: \mathsf{Crys}^r(X)\to \mathsf{Crys}^r(X)^{S^1},
	\end{align*} 
	and these functors commute with $f_{\mathrm{dR}}^!$ for $f:X\to Y$. They also commute with exterior products in all arities, since the diagonal equivariant structure on an exterior product of trivially equivariant objects is also trivial.  Moreover, $t_X$ is $t$-exact, by \cite[Section 3.1.3]{ES-Kashiwara}. Thus, since $\mathsf{Crys}^r(X)$ is regular and $t_X$ is exact, it lifts to a $k$-linear functor
	\begin{align*}
		\tilde{t}_X: \mathsf{Crys}^r(X)\to \mathscr{R}(\mathsf{Crys}^r(X)^{S^1})
	\end{align*} 
	by \zcref{lem: multilinear regularization functor}. Moreover, the same lemma gives a natural lift in $X$; hence, there is a well-defined $k$-linear natural transformation of lax symmetric monoidal functors
	\begin{align*}
		\mathsf{Crys}^r\to \mathscr{R}(\mathsf{Crys}^r(-)^{S^1}).
	\end{align*} 
	Extending $k[u]$-linearly thus gives a natural transformation
	\begin{align*}
		\mathsf{Crys}^r\otimes k[u]\mod\to \mathscr{R}(\mathsf{Crys}^r(-)^{S^1}),
	\end{align*} 
	and inverting $u$ gives the desired natural transformation. It is an equivalence by \cite[Lemma 4.5.4 and Proposition 4.5.5]{Preygel} (cf. also Remark 4.5.6 of \textit{op. cit.}).
	\end{proof}
	
\begin{lem}\label{lem: crystal 6FF}
	The $(\infty,2)$-functor $\mathsf{Crys}^r_{\mathbb{Z} /2}$ induces a lax symmetric monoidal $(\infty,2)$-functor
	\begin{align*}
		\bSpan(\mathsf{Sch}_{\mathrm{aft}},\mathrm{proper})_{\mathrm{equiv},\mathrm{proper}}\to \bDGCat_{\ku},
	\end{align*} 
	which we also denote by $\mathsf{Crys}^r_{\mathbb{Z} /2}$ by abuse of notation, for the classes of maps $E=I$ and $P$ given respectively by proper maps and equivalences.
\end{lem} 
\begin{proof}
	Since the extension of scalars functor $\bDGCat_k\to \bDGCat_{\ku}$ is strictly symmetric monoidal, it suffices to show that the functor
	\begin{align*}
		(\mathsf{Sch}_{\text{aft}})^{\text{op}}\to &\bDGCat_k\\
		X\mapsto &\mathsf{IndCoh}(X_{\text{dR}})\eqcolon \mathsf{Crys}^r(X)
	\end{align*} 
	induces a lax symmetric monoidal $(\infty,2)$-functor
	\begin{align*}
		\mathsf{Crys}^r:\bSpan(\mathsf{Sch}_{\text{aft}},E)_{I,P}\to \bDGCat_k.
	\end{align*} 
	We write $\mathsf{Crys}^r:(\mathsf{Sch}_{\text{aft}})^{\text{op}}\to\bDGCat_k$ as the composition of $\mathsf{IndCoh}$ with $(-)_{\text{dR}}$, and we check that it satisfies the conditions of \cite[Definition 5.13]{CLL_six-functor-formalisms} and Theorem 5.14 of \textit{op. cit.}, which will allow us to conclude.

	We first check that given maps $f:X\to Y$ and $g:Y\to Z$ in $\mathsf{Sch}_{\text{aft}}$ such that $g$ and $g \circ f$ are proper, then $f$ is proper. Since properness is checked on the classical truncation, we can assume that $X,Y,Z$ are classical. We factor $f$ as
	\begin{align*}
		X\to X\times _ZY\to Y,
	\end{align*} 
	where the projection to $Y$ is proper because it is the base change of a proper morphism. The map $X\to X\times _ZY$ is the graph of $f$, which is a closed immersion because $g$ is separated by assumption; hence, it is also proper, so $f$ is proper. The rest of the conditions on $\mathsf{Sch}_{\text{aft}}$ required in \cite[Theorem 5.14]{CLL_six-functor-formalisms} hold: Every equivalence is truncated, and the product functor $\mathsf{Sch}_{\text{aft}}\times \mathsf{Sch}_{\text{aft}}\to \mathsf{Sch}_{\text{aft}}$ preserves pullbacks because it is right adjoint to the diagonal functor.

	If $i:X\to Y$ is a map in $I$, then by \cite[Chapter 4 Corollary 1.3.6]{GR17II}, $i_{\text{dR}}$ is inf-schematic, hence ind-inf-schematic. Similarly, $i_{\text{dR}}$ is proper, hence ind-proper, because its underlying map on reduced stacks is $i_{\text{red}}$.\footnote{Recall that $(X_{\text{dR}})_{\text{red}}=X_{\text{red}}$.} It remains to check that the conditions of \cite[Definition 5.13]{CLL_six-functor-formalisms} are satisfied. As before, by Example 2.7 of \textit{op. cit.}, it suffices to check conditions (1) and (2a). Since $i$ is proper, $i_*\dashv i^!$. Suppose given a Cartesian diagram
	\[ 
	\begin{tikzcd}
		X'\arrow[r,"i'"]\arrow[d,"g"]&Y'\arrow[d,"f"]\\
		X\arrow[r,"i"]&Y
	\end{tikzcd}.
	\] 
	Recalling that $(-)_{\text{dR}}$ preserves limits by \cite[Lemma 1.1.4]{GRCrystals}, the required left mate for condition (1) is
	\begin{align*}
		(i'_{\text{dR}})_*g_{\text{dR}}^!\to f_{\text{dR}}^!(i_{\text{dR}})_*.
	\end{align*} 
	Indeed, this mate is invertible by \cite[Chapter 4 Section 2.1.3]{GR17II}. We now check condition (2a). Since $(-)_{\text{dR}}$ preserves limits, this condition follows using the same proof as in \zcref{prop: indcoh is a 6-functor formalism}.
\end{proof} 
\subsection{Relating potent $\mathcal{D}$-modules to crystals}
\label{ssec: potent D-modules=crystals}
By \zcref{prop: formal completions are ind-ft}, $\mathsf{Sch}_{\text{aft}}$ is a full subcategory of $\mathsf{indGStk}_{\text{ft}}$. Since ind-inf-schematic ind-proper maps between almost finite type schemes are proper, the six-functor formalism $\Pot$ constructed in \zcref{thm:potent-D-modules-is-a-6FF} restricts to a six-functor formalism
\begin{align*}
	\Pot_{\text{sch}}\coloneq \Pot|_{\bSpan(\mathsf{Sch}_{\mathrm{aft}},\mathrm{proper})_{\mathrm{equiv},\mathrm{proper}}}:\bSpan(\mathsf{Sch}_{\text{aft}},\text{proper})_{\text{equiv},\text{proper}}\to \bDGCat_{\ku}.
\end{align*} 
\begin{lem}\label{lem: nat transformation of constant loops}
	There is a 1-morphism in
	\begin{align*}
		\Fun(\mathsf{Sch}_{\mathrm{aft}},\Fun(BS^1,\mathsf{PreStk}_{\mathrm{laft}}))
	\end{align*} 
	whose value on a scheme $X$ is the $S^1$-equivariant constant loop morphism
	\begin{align*}
		\triv_{\mathrm{sch}}^{S^1}(X)\to \mathcal{L}X,
	\end{align*} 
	where $\triv_{\mathrm{sch}}^{S^1}: \mathsf{Sch}_{\mathrm{aft}}\to \Fun(BS^1,\mathsf{PreStk}_{\mathrm{laft}})$ equips schemes with the trivial $S^1$-action.
\end{lem} 
\begin{proof}
	There is a functor $P:BS^1\times \Delta^1\to \mathsf{Spc}$ whose restriction to $0$ is $S^1$ and restriction to $1$ is $*$. Indeed, loop rotation defines a functor $R:BS^1\to \mathsf{Spc}$ whose value at $*\in BS^1$ is $S^1$. Since the point space equipped with the trivial $S^1$-action $*$ is terminal in $\Fun(BS^1,\mathsf{Spc})$, there is a map $R\to *$ in that functor category. Uncurrying this map gives $P$. Consider the composite
	\begin{align*}
		(BS^1)^{\text{op}}\times (\Delta^1)^{\text{op}}\times \mathsf{Sch}_{\text{aft}}\xrightarrow{P^{\text{op}}\times \id}\mathsf{Spc}^{\text{op}}\times \mathsf{Sch}_{\text{aft}}\xrightarrow{(K,X)\mapsto \Map(K,X)}\mathsf{PreStk}.
	\end{align*} 
	Currying the variables besides the $(\Delta^1)^{\text{op}}$ variable and using the equivalence $(BS^1)^{\text{op}}\simeq BS^1$ we obtain a 1-morphism in
	\begin{align*}
		\Fun(\mathsf{Sch}_{\text{aft}},\Fun(BS^1,\mathsf{PreStk}_\text{laft})).
	\end{align*} 
	The value of this 1-morphism at $X\in \mathsf{Sch}_{\text{aft}}$ is then the constant loop map $\triv_{\mathrm{sch}}^{S^1}(X)\to \mathcal{L}X$.
\end{proof} 
\begin{lem}\label{lem: nat transformation between loops and de Rham}
	There is a 1-morphism $\pi$ in
	\begin{align*}
		\Fun(\mathsf{Sch}_{\mathrm{aft}},\Fun(BS^1,\mathsf{PreStk}_\mathrm{laft}))
	\end{align*} 
	whose value on a scheme $X$ is the map
	\begin{align*}
		\pi_X: \mathcal{L}X\xrightarrow{\eta_{\mathcal{L}X}}(\mathcal{L}X)_{\mathrm{dR}}\xrightarrow{\simeq}X_{\mathrm{dR}}
	\end{align*} 
	of \cite[Section 1.2.2]{Preygel}, where $\eta$ denotes the unit of the adjunction
	\[
	\begin{tikzcd}
		(-)_{\mathrm{red}}:\mathsf{PreStk}\arrow[r,shift left]& \mathsf{PreStk}_{\mathrm{cl,red}}:(-)_{\mathrm{dR}}\arrow[l,shift left].
	\end{tikzcd}
	\] 
\end{lem} 
\begin{proof}
	By \zcref{lem: nat transformation of constant loops} and functoriality of $(-)_{\mathrm{dR}}$, there is a natural transformation
	\begin{align*}
		(-)_{\mathrm{dR}}\circ\triv_{\mathrm{sch}}^{S^1}\Rightarrow (-)_{\mathrm{dR}}\circ \mathcal{L}
	\end{align*} 
	of functors in $\Fun(\mathsf{Sch}_{\mathrm{aft}},\Fun(BS^1,\mathsf{PreStk}_\text{laft}))$. For $X\in\mathsf{Sch}_{\text{aft}}$, $\mathcal{L}X$ and $X$ have the same underlying reduced, classical prestack; hence, this natural transformation is an equivalence and admits an inverse. Since $\eta$ is also a natural transformation, we see that its composition with this inverse defines a natural transformation $\pi$ as desired.
\end{proof} 
 
\begin{prop}\label{prop: potent D-modules are crystals}
	The natural transformation $\pi$ of \zcref{lem: nat transformation between loops and de Rham} induces an equivalence
	\begin{align*}
		\Phi:\mathsf{Crys}^r_{\mathbb{Z}/2}\xRightarrow{\simeq}\Pot_{\mathrm{sch}}.
	\end{align*} 
\end{prop} 
\begin{proof}
	We begin by observing that the natural transformation $\pi$ of \zcref{lem: nat transformation between loops and de Rham} canonically lifts to a symmetric monoidal natural transformation, because the functors $\mathcal{L}$, $\triv_{\mathrm{sch}}^{S^1}$, and $(-)_{\mathrm{dR}}$ all commute with products of schemes. Thus, $\pi^{!,S^1}\coloneq \mathsf{IndCoh}(\pi)^{S^1}$ is a lax symmetric monoidal natural transformation, by lax symmetric monoidality of $\mathsf{IndCoh}$ and $S^1$-invariants (\zcref{prop: indcoh is a 6-functor formalism} and \zcref{lem: S^1-invariants is lax symmetric monoidal}). By \cite[Theorem 5.1.3]{Preygel} (cf. \cite[Example 5.2.1]{Preygel}) and \zcref{lem: crystals factor through Coh_t}, $\pi^!_X$ is left $t$-exact up to a shift. Since the forgetful functors from $S^1$-invariants detect lower bounds and intertwine $\pi_X^!$ and $\pi_X^{!,S^1}$, we see that $\pi_X^{!,S^1}$ is also left $t$-exact up to a shift. The functor $\mathsf{IndCoh}(\mathcal{L}-)^{S^1}$ factors through $\Coh_{t,k[u]}^\otimes $ by \zcref{prop: IndCoh(L-)^S^1 as a map of operads}; similarly, the functor $\mathsf{Crys}^r(-)^{S^1}$ does as well by \zcref{lem: crystals factor through Coh_t}. Thus, $\pi_X^{!,S^1}$ is a unary operation in $\mathsf{Coh}_{t,k[u]}^\otimes $ for any $X\in \mathsf{Sch}_{\mathrm{aft}}$. Moreover, the source and target functors of $\pi^{!,S^1}$ factor through this operad; hence, $\pi^{!,S^1}$ defines a morphism
	\begin{align*}
		\mathsf{Crys}^r(-)^{S^1}\Rightarrow  \mathsf{IndCoh}(\mathcal{L}-)^{S^1}
	\end{align*} 
	in $\mathsf{Alg}_{\mathsf{Sch}_{\mathrm{aft}}^{\mathrm{op}}}(\mathsf{Coh}_{t,k[u]})$. Postcomposition with the functor $\oblv^\otimes \circ q \circ p^{-1}$ of \zcref[noname]{eq: k-linear D^pot} defines a functor
	\begin{align*}
		\mathsf{Alg}_{\mathsf{Sch}_{\mathrm{aft}}^{\mathrm{op}}}(\mathsf{Coh}_{t,k[u]})\to \mathsf{Alg}_{\mathsf{Sch}_{\mathrm{aft}}^{\mathrm{op}}}(\DGCat_{k[u]})
	\end{align*} 
	whose value on $\pi^{!,S^1}$ is thus a lax symmetric monoidal natural transformation
	\begin{align*}
		\mathscr{R}\left(\mathsf{IndCoh}((-)_{\mathrm{dR}})^{S^1}\right)\Rightarrow  \mathscr{R}\left(\mathsf{IndCoh}(\mathcal{L}-)^{S^1}\right)
	\end{align*} 
	with component $\mathscr{R}(\pi^{!,S^1}_X)$ at $X\in \mathsf{Sch}_{\mathrm{aft}}$. Finally, extending scalars along $k[u]\to\ku$ gives a lax symmetric monoidal natural transformation
	\begin{align*}
		\pi^{!,\tT}: \mathsf{Crys}^r(-)^{\tT}\Rightarrow \Pot_{\mathrm{sch}}.
	\end{align*} 
	Now, define $\Phi: \mathsf{Crys}^r_{\mathbb{Z}/2}\Rightarrow \Pot_{\mathrm{sch}}$ as the composition
	\begin{align*}
		\mathsf{Crys}^r_{\mathbb{Z}/2}\xRightarrow{\xi}\mathsf{Crys}^r(-)^{\tT}\xRightarrow{\pi^{!,\tT}}\Pot_{\mathrm{sch}}
	\end{align*} 
	where $\xi$ is the natural equivalence of \zcref{lem: periodic crystals are tate invariants of crystals}. We now claim that $\Phi$ is a map of six-functor formalisms, which means that it is a 1-morphism in the image of the functor $h^*$ of \cite[Theorem 5.14]{CLL_six-functor-formalisms}. Thus, we need to show 
	\begin{enumerate}[label=\arabic*)]
		\item $\Phi$ is a natural transformation of lax symmetric monoidal $(\infty,2)$-functors;
		\item The natural transformation of functors underlying $\Phi$ is $(\mathrm{proper},\mathrm{equiv})$-biadjointable.
	\end{enumerate} 
	Item (1) follows because $\Phi$ is the composite of the lax symmetric monoidal natural transformations $\xi$ and $\pi^{!,\tT}$. For item (2), we observe that for fixed $X\in \mathsf{Sch}_{\mathrm{aft}}$,  
	\begin{align*}
		\Phi_X:\mathsf{Crys}^r_{\mathbb{Z}/2}(X)\to \Pot(X)
	\end{align*} 
	is an equivalence by \cite[Theorem 1.3.5]{Preygel}. Thus, $\Phi$ is a natural equivalence and hence is automatically $(\mathrm{proper},\mathrm{equiv})$-biadjointable, proving item (2). \cite[Theorem 5.14]{CLL_six-functor-formalisms} thus extends $\Phi$ to a natural equivalence of six-functor formalisms.
\end{proof}

%--------------------------------------------------------------------------------------------%
%--------------------------------------BIBLIOGRAPHY-------------------------------------%
%--------------------------------------------------------------------------------------------%

\printbibliography

\end{document}